\documentclass{coulonpaper}

\title{Small cancellation theory and Burnside problem.}
\author{R\'emi Coulon}

\begin{document} 

\maketitle

\begin{abstract}
In these notes we detail the geometrical approach of small cancellation theory used by T.~Delzant and M.~Gromov to provide a new proof of the infiniteness of free Burnside groups and periodic quotients of torsion-free hyperbolic groups.
\end{abstract}

\tableofcontents 

% !TEX root = notes.tex
\section{Introduction}

\paragraph{} Let $n$ be an integer.
A group $G$ has exponent $n$ if for all $g\in G$, $g^n=1$.
In 1902 W.~Burnside asked whether a finitely generated group with finite exponent is necessarily finite or not \cite{Bur02}.
To study this question it is natural to look at the free Burnside group $\burn rn = \free r/ {\free r}^n$, which is the quotient of the free group of rank $r$ by the (normal) subgroup generated by the $n$-th powers of all elements.
It is indeed the largest group of rank $r$ and exponent $n$.
In 1968, P. S.~Novikov and S. I.~Adian answered negatively the Burnside Problem.
They proved that if $r\geq 2$ and $n$ is odd larger than 4381, then $\burn rn$ is infinite  \cite{NovAdj68c}.
Since, this result has been extended in many directions (even exponents, periodic quotient of hyperbolic groups, etc)  \cite{Olc82,Olc91,Iva94,Lys96}.
More recently T.~Delzant and M.~Gromov provided a new approach of the Burnside problem \cite{DelGro08}.
In particular, they give an alternative proof of the following theorem.

\begin{theo}[{\cite[Th. 6.2.2]{DelGro08}}, see also \cite{Olc91}]
\label{res: Burnside problem hyperbolic groups}
	Let $G$ be a non-cyclic, torsion-free hyperbolic group.
	There exists an integer $n_0$ such that for all odd exponents $n \geq n_0$, the quotient $G/G^n$ is infinite.
\end{theo}

The aim of these notes is to give a comprehensive presentation of their method.
Actually a more general statement holds for an arbitrary hyperbolic group.

\begin{theo}[{\cite[Th. A]{IvaOlc96}}]
	Let $G$ be a non-virtually cyclic hyperbolic group.
	For every integer $n_0$ there exists an exponent $n \geq n_0$ such that the quotient $G/G^n$ is infinite.
\end{theo}

However this second result require to control the even torsion, which is much harder.
The notes do not cover this theorem.

\paragraph{}
From a geometrical point of view the difficulty to study Burnside groups comes from the fact that we have no ``nice'' metric space on which $\burn rn$ is acting by isometries.
Therefore the idea is to study $\burn rn$ as the direct limit of a sequence of groups,
\begin{displaymath}
	\free r = G_0 \twoheadrightarrow G_1 \twoheadrightarrow G_3 \twoheadrightarrow G_4 \twoheadrightarrow \dots \twoheadrightarrow G_k \twoheadrightarrow G_{k+1} \twoheadrightarrow \dots 
\end{displaymath}
where each $G_k$ is easier to understand. 
In this construction, T.~Delzant and M.~Gromov obtain $G_{k+1}$ as a quotient of $G_k$ using a geometrical version of the small cancellation theory.
It forces the groups $G_k$ to be non-virtually cyclic and hyperbolic.
Thus their limit $\burn rn$ cannot be finite.
All the other known strategies also use such a sequence of groups.
Before explaining the construction of T.~Delzant and M.~Gromov let us recall a few ideas about the usual small cancellation theory and its geometric generalization given by M.~Gromov \cite{Gro01b}.

\subsection{Usual small cancellation theory.} 
For more details about the usual small cancellation theory, we refer the reader to  \cite{LynSch77}.
Let $\mathbf F(S)$ be the free group generated by a finite set $S$. 
Let $R$ be a set of words over the alphabet $S \cup S^{-1}$.
The goal is to study the group $\bar G = \mathbf F(S) /\lnormal R \rnormal$, where $\lnormal R \rnormal$ stands for the normal subgroup of $\mathbf F(S)$ generated by $R$.
We assume that the elements of $R$ are non-trivial and  cyclically reduced.
We denote by $R^*$ the set of all cyclic conjugates of elements of $R\cup R^{-1}$.
A \emph{piece} is a common prefix of two distinct elements of $R^*$.
In  other words a piece is a subword that could cancel in the product $rs$ where $r,s \in R^*$.
Let $\lambda > 0$.
One says that $R$ satisfies the small cancellation assumption $C'(\lambda)$ if for all pieces $u$, for all relations $r \in R$ containing $u$, $|u| \leq \lambda |r|$.
Let us mention one important theorem.
In the next paragraph we will provide an extension of it.

\begin{theo}[see {\cite{GhyHar90}}]
\label{res: usual small cancellation theorem}
	Let $\lambda \in \left(0,  1/6\right)$.
	Let $R$ be a set of non-trivial, cyclically reduced words over the alphabet $S \cup S^{-1}$.
	If $R$ satisfies the condition $C'(\lambda)$ then the quotient $\mathbf F(S)/\lnormal R \rnormal$ is a hyperbolic group.
\end{theo}

Four our purpose we are going to consider a more stronger condition called $C''(\lambda)$.
It requires that for all pieces $u$, for all relations $r \in R$ (not necessarily containing $u$), $|u| \leq \lambda |r|$.
This assumption can be reformulated in a more geometrical way.
To that end we consider the Cayley graph of $\mathbf F(S)$ with respect to $S$, denoted by $X$ which is a tree.
Let $r \in R^*$.
It fixes two points $r^-$ and $r^+$ of the boundary at infinity $\partial X$ of $X$.
These points are joined by a bi-infinite geodesic called the \emph{axis} of $r$.
We denote it by $Y_r$.
The isometry $r$ acts on $Y_r$ by translation of length $\len r = \inf_{x \in X} \dist {rx}x$, where $\dist xy$ denotes the distance between $x$ and $y$.
Since $r$ is a cyclically reduced word, $\len r$ is in fact the length of the relation $r$.
Consider now two relations $r, s \in R^*$.
The length of a common piece of $r$ and $s$ is exactly $\diam\left(Y_r\cap Y_s\right)$.
Thus the $C''(\lambda)$ condition can be stated in this way:
\begin{displaymath}
	\sup_{r\neq s} \diam\left(Y_r\cap Y_s\right) \leq \lambda \inf_r\len r.
\end{displaymath}
With this idea in mind we provide in the next paragraph a larger framework to the small cancellation theory.

\subsection{Small cancellation theory in hyperbolic groups}

\paragraph{} From now on $X$ is a proper geodesic $\delta$-hyperbolic space.
Let $G$ be a group acting properly co-compactly by isometries on $X$.
The translation length of an isometry $g\in G$, denoted by $\len g$ is the quantity $\len g = \inf_{x \in X}\dist {gx}x$.
We consider a collection $\mathcal Q$ of pairs $(H,Y)$ satisfying the following properties
\begin{enumerate}
	\item $Y$ is a $2\delta$-quasi-convex subset of $X$ and $H$ a subgroup of $\stab{Y}$, the stabilizer in $G$ of $Y$, acting co-compactly on $Y$.
	\item the family $\mathcal Q$ is invariant under the action of  $G$ i.e., for all $g \in G$, for all $(H,Y) \in \mathcal Q$, $(gHg^{-1},gY)$ is an element of $\mathcal Q$,
	\item the set $\mathcal Q/G$ is finite.
\end{enumerate}
Our goal is to study the group $\bar G = G/K$, where $K$ is the (normal) subgroup of $G$ generated by the $H$'s.
For this purpose, we define two quantities which respectively play the role of the length of the largest piece and the length of the smallest relation.
\begin{itemize}
	\item
	The maximal overlap between two quasi-convex of $\mathcal Q$ is measured by the quantity
	\begin{displaymath}
		\Delta(\mathcal Q) = \sup_{(H_1,Y_1)\neq (H_2,Y_2)}\diam \left(Y_1^{+5\delta}\cap Y_2^{+5\delta}\right).
	\end{displaymath}
	\item
	The smallest translation length of $\mathcal Q$ is defined as
	\begin{displaymath}
		T(\mathcal Q) = \inf\left\{\len h \mid h\in H\setminus \{1\}, \; (H,Y) \in \mathcal Q \right\}.
	\end{displaymath}
\end{itemize}
Note that if $\Delta(\mathcal Q)$ is finite then two distinct groups $H$ cannot be associated in $\mathcal Q$ to the same quasi-convex $Y$.
In particular $H$ is a normal subgroup of $\stab Y$.
We are interested in situations where the ratios $\Delta(\mathcal Q)/T(\mathcal Q)$ and $\delta/T(\mathcal Q)$ are very small (see Theorem~\ref{res: rotation small cancellation theorem rescaled}).
In particular it forces every non-trivial element of every subgroup $H$ to be hyperbolic.

\exs 
\begin{enumerate}
	\item \label{exa: usual small cancellation}
	Let $R$ be a finite set of elements of $\mathbf F(S)$.
	Consider the collection 
	\begin{displaymath}
		\mathcal Q = \set{\left( \left< uru^{-1}\right>,uY_r\right)}{r \in R,\ u \in \mathbf F(S) }.
	\end{displaymath}
	In this case we recover the usual small cancellation theory.
	Since the the Cayley graph of $\mathbf F(S)$ is 0-hyperbolic, the assumption $C''(\lambda)$ is equivalent to $\Delta(\mathcal Q) \leq \lambda T(\mathcal Q)$
	
	\item \label{exa: small cancellation with graphs}
	The next example comes from small cancellation with graphs (see \cite{Oll06}).
	Let $\Gamma$ be a finite connected graph whose edges are labelled with elements of $S \cup S^{-1}$.
	Given a vertex $v_0$ of $\Gamma$, the labeling provides a homomorphism from the fundamental group of $\left(\Gamma, v_0\right)$ onto a subgroup $H$ of $\mathbf F(S)$.
	Moreover, we have a map from $\tilde \Gamma$, the universal cover of $\Gamma$, onto a subtree $T$ of the Cayley graph of $\mathbf F(S)$.
	We can consider the collection
	\begin{displaymath}
		\mathcal Q = \set{\left(gHg^{-1},gT\right)}{g\in G/H}.
	\end{displaymath}
	In this situation a \emph{piece} is, by definition, a  word which labels two distinct simple paths of $\Gamma$.
	Then $\Delta(\mathcal Q)$ is the length of the largest piece whereas $T(\mathcal Q)$ is the girth of the graph $\Gamma$ i.e., the length of the smallest embedded loop in $\Gamma$.
	For further detail we refer the reader to \cite[1/6-theorem]{Gro03} and \cite[Th.1]{Oll06}.
	
	\item \label{exa: small cancellation with large powers}
	Let $G$ be a group acting properly, co-compactly by isometries on a geodesic, $\delta$-hyperbolic space $X$.
	If $r \in G$ is a hyperbolic isometry it fixes two points $r^-$ and $r^+$ of $\partial X$, the boundary at infinity of $X$.
	Following the case of free groups we define $Y_r$ to be the set of all points of $X$ which are $10\delta$-close to some bi-infinite geodesic joining $r^-$ and $r^+$.
	This set is $2\delta$-quasi-convex.
	We consider the set $P$ of all hyperbolic elements $r \in G$ which are not proper powers such that $\len r \leq 1000\delta$.
	To every integer $n$ we associate the collection
	\begin{displaymath}
		\mathcal Q_n=\set{\left(\left<r^n\right>,Y_r\right)}{r\in P}
	\end{displaymath}
	Assume now that $G$ is non-virtually cyclic.
	Since $G$ acts properly co-compactly on $X$, the overlap $\Delta(\mathcal Q_n)$ is uniformly bounded.
	Therefore, by choosing an integer $n$ large enough, one can obtain a ratio $\Delta(\mathcal Q_n)/T(\mathcal Q_n)$ as small as desired.
\end{enumerate}

\paragraph{} Consider now an arbitrary family $\mathcal Q$ as defined previously.
In order to study the quotient $\bar G = G/K$ we construct a space $\bar X$ on which $\bar G$ acts properly, co-compactly, by isometries.
Roughly speaking, $\bar X$ is obtained as follows.
We fix a large real number $\rho$.
Its value will be made precise later.
For each $(H,Y) \in \mathcal Q$, we start by attaching on $X$ along $Y$ a cone of base $Y$ and radius $\rho$.
The cones are endowed with a negatively curved metric modeled on the distance of the hyperbolic plane $\H_2$.
This new space $\dot X$ (called cone-off of $X$)  inherits a metric coming from the distance on $X$ and the ones on the cones.
Since $\mathcal Q$ is endowed with an action of $G$, the action of $G$ on $X$ extends by homogeneity into an action of $G$ on $\dot X$.
On this space every group $H$ acts as a rotation of very large angle that fixes the apex of the cone of base $Y$.
One speaks of \emph{rotation family} (see Section~\ref{sec: rotation family}).
The space $\bar X$ is then the quotient of $\dot X$ by $K$.
By construction $\bar X$ is proper and geodesic.
Moreover, $\bar G$ acts properly co-compactly by isometries on it.
The study of the metric on $\dot X$ and $\bar X$ is a key question.
This will be done in Sections~\ref{sec: cone}-\ref{sec: small cancellation theory}.
One of the main features of the space $\bar X$ is the following.
Assume that the ratio $\Delta(\mathcal Q)/T(\mathcal Q)$ is very small, then every ball of $\bar X$ of radius $\rho/50$ is roughly $\boldsymbol \delta$-hyperbolic where $\boldsymbol \delta$ is the hyperbolicity constant of the hyperbolic plane $\H_2$.
If $\rho$ has been chosen large enough, we can then apply a Cartan-Hadamard Theorem (see Theorem~\ref{res: cartan hadamard}).
It leads to the following analogue of Theorem~\ref{res: usual small cancellation theorem}.

\begin{theo}
\label{res: rotation small cancellation theorem rescaled}
	There exist positive constants $\delta_0$ and $\lambda_0$ (which do not depend on $X$, $G$ or $\mathcal Q$) with the following property.
	Assume that 
	\begin{displaymath}
		\frac \delta{T(\mathcal Q)} \leq \delta_0 \quad \text{ and }\quad  \frac {\Delta(\mathcal Q)}{T(\mathcal Q)} \leq \lambda_0.
	\end{displaymath}
	Then  $\bar X$ is a hyperbolic space endowed with a proper co-compact action of $\bar G$.
\end{theo}

\paragraph{} Some hypotheses of the theorem can be weaken.
For instance we do not really need that the action of $G$ (\resp $H$) on $X$ (\resp $Y$) to be co-compact.
Without these assumptions, the space $\bar X$ remains hyperbolic nevertheless the action of $\bar G$ on $\bar X$ is no more proper or co-compact.

\subsection{Iterating small cancellation}
Theorem~\ref{res: rotation small cancellation theorem rescaled} is an important step of the induction process to prove the infiniteness of the Burnside groups.
However this is definitely not sufficient. 
The space $\bar X$ has in fact many other properties that we will enlighten in Section~\ref{sec: small cancellation theory}.
Let us explain briefly where the difficulty does come from.
Let $G_0$ be a non-virtually cyclic torsion-free group acting properly co-compactly on a $\delta$-hyperbolic space $X_0$ and $n_0$ an integer.
Following Example~\ref{exa: small cancellation with large powers}, we consider the quotient $G_1 = G_0/K_0$ where $K_0$ is the normal subgroup of $G_0$ generated by 
\begin{displaymath}
	\left\{r^{n_0}/ r\in G, \text{ hyperbolic, not a proper power, } \len r \leq 1000 \delta\right\}
\end{displaymath}
If $n_0$ is large enough, Theorem~\ref{res: rotation small cancellation theorem rescaled} applies. 
The quotient $G_1$ is also a hyperbolic group.
Thus we would like to iterate the process and kill some large powers in $G_1$.
According to Example~\ref{exa: small cancellation with large powers} this can be done, but there is no reason that the exponent $n_1$ that works for the second step should be the same as $n_0$.
In this way one can construct by induction an infinite torsion group.
Nevertheless to study Burnside groups we need to follow during the construction some parameters that allow us to use at each step the same exponent $n$.
Note already that the constants $\delta_0$ and $\lambda_0$ in Theorem~\ref{res: rotation small cancellation theorem rescaled} do not depend on $G$, $X$ or $\mathcal Q$.
This is an important but not sufficient fact.
To control the ratio $\Delta(\mathcal Q)/T(\mathcal Q)$ we will take care of two quantities associated to the action of $G$ on $X$. 
The radius of injectivity which represent the smallest asymptotic translation length of a hyperbolic element of $G$ and the invariant $A(G,X)$ which measures the maximal overlap between the axes of two small isometries of $G$ (see Section~\ref{sec: group acting on a hyperbolic space} for the precise definitions).
During this construction we will require the exponent $n$ to be odd. 
With this assumption we can prove that every subgroup of any $G_k$ is either cyclic or contains a free group. 
It will help to control $\Delta(\mathcal Q)$ in terms of $A(G,X)$.
Burnside groups of large even exponents are known to be infinite \cite{Iva94,Lys96}.
However the proof is much harder.

\paragraph{}
The hypotheses required to apply small cancellation are not very restricting. 
See for instance \cite{Olc92}.
Therefore iterating small cancellation can be exploited to build many groups with pathological properties.
It is one of the ingredients involved in the construction of Gromov's monster group \cite{Gro03,ArzDel08}.
This object is a limit of groups $G_k$ where $G_{k+1}$ is obtained from $G_k$ using graphical small cancellation (see Example~\ref{exa: small cancellation with graphs}).
One remarkable property of this group is that it does not coarsely embed into a Hilbert space and thus does not satisfy the Baum-Connes conjecture with coefficients.

\paragraph{}
For other approaches of iterated small cancellation theory we refer the reader to \cite{Rips:1982ge,Olc91a,McCammond:2000he}

\subsection{Outline of the paper}
In the Section~\ref{sec: hyperbolic geometry} we recall the main features of hyperbolic spaces in the sense of Gromov and groups acting on a hyperbolic space.
Section~\ref{sec: rotation family} is dedicated to the study of rotation families. 
This section differs from the tools used by T.~Delzant and M.~Gromov.
In this way we do not need to use orbifolds as in \cite{DelGro08}.
Sections~\ref{sec: cone} and \ref{sec: cone-off construction} detail the cone-off construction.
The core of the proof is contained in Section~\ref{sec: small cancellation theory}.
We prove among others the small cancellation theorem (see Theorem~\ref{res: SC - small cancellation theorem}) and investigate the properties of the space $\bar X$.
Thanks to the small cancellation we state at the beginning of Section~\ref{sec: application to periodic groups} an induction lemma which explains how the group $G_{k+1}$ can be obtained from $G_k$.
In particular it provides a set of assumptions which if they are satisfied for $G_k$ also hold for $G_{k+1}$.
This allow us to iterate the construction and prove the main theorem (see Theorem~\ref{res : SC - periodic quotient}).
Our approach of the small cancellation theory is very broad. 
In particular it covers all the main results that are needed to construct Gromov's monster group. 
In Section~\ref{sec: gromov monster} we also give some keys to understand this construction as explained in \cite{ArzDel08}.
In the appendix, we propose an alternative proof of the Cartan-Hadamard Theorem based on the ideas given in \cite{DelGro08}.

\paragraph{Acknowledgments.} These notes follow a series of lectures given at the Vanderbilt University during Fall 2011. The author wishes to express his gratitude to M.~Mihalik, A.~Ol'shanski\u\i, D.~Osin and M.~Sapir for their patient and active interest in this work and also for encouraging him to write this text.
He is also greatly thankful to T.~Delzant who guides his first steps in this field.

% !TEX root = notes.tex

\section{Hyperbolic geometry}
\label{sec: hyperbolic geometry}

	\paragraph{} In this section we some of the standard facts on hyperbolic spaces in the sense of M.~Gromov.
	We only give the proofs of quantitative results, which are not a straight forward application of the four points condition. 
	For more details we refer the reader to the original paper of M.~Gromov \cite{Gro87} or to \cite{CooDelPap90,GhyHar90}.
	
	\paragraph{}Let $X$ be a metric length space.
	Unless otherwise stated a path is a rectifiable path parametrized by arclength.
	Given two points $x$ and $x'$ of $X$, we denote by $\dist[X]x{x'}$ (or simply $\dist x{x'}$) the distance between them.
	We write $B(x,r)$ for the closed ball of center $x$ and radius $r$ i.e., the set of points $y \in X$ such that $\dist xy \leq r$.
	Let $Y$ be a subset of $X$.
	The distance between a point $x$ of $X$ and $Y$ is denoted by $d(x,Y)$.
	We write $Y^{+\alpha}$ for the $\alpha$-neighborhood of $Y$ i.e., the set of points $x \in X$ such that $d(x,Y) \leq \alpha$.
	Let $\eta \geq 0$.
	A point $p$ of $Y$ is an \emph{$\eta$-projection} of $x \in X$ on $Y$ if $\dist xp \leq d(x,Y) +\eta$.
	A 0-projection is simply called a \emph{projection}.

\subsection{The four points inequality.}

	\paragraph{}
	The \emph{Gromov product} of three points $x,y,z \in X$  is defined by 
	\begin{displaymath}
		\gro xyz = \frac 12 \left\{ \fantomB \dist xz + \dist yz - \dist xy \right\}.
	\end{displaymath}
	The space $X$ is \emph{$\delta$-hyperbolic} if for every $x,y,z,t \in X$
	\begin{equation}
	\label{eqn: hyperbolicity condition 1}
		\gro xzt \geq \min\left\{\fantomB \gro xyt, \gro yzt \right\} - \delta,
	\end{equation}
	or equivalently
	\begin{equation}
	\label{eqn: hyperbolicity condition 2}
		\dist xz + \dist yt \leq \max\left\{ \fantomB \dist xy + \dist zt, \dist yz + \dist xt \right\} +2\delta.
	\end{equation}
	To show that a space is hyperbolic, it is actually sufficient to prove (\ref{eqn: hyperbolicity condition 1}) with a fixed base point.
	\begin{lemm}[see {\cite[Chap.1 Prop. 1.2]{CooDelPap90}}]
	\label{res: criterion for hyperbolicity}
		Let $\delta \geq 0$. 
		Let $(X,t)$ be a pointed metric space.
		If for every $x,y,z \in X$, we have $\gro xzt \geq \min\left\{\gro xyt, \gro yzt \right\} - \delta$, then $X$ is $2\delta$-hyperbolic.
	\end{lemm}
	
	\rems Note that in the definition of hyperbolicity we do not assume that $X$ is a geodesic.
	In Section~\ref{sec: cone-off construction} we construct a length space (the cone-off) for which it is not obvious to prove that it is geodesic.
	However it satisfies (\ref{eqn: hyperbolicity condition 1}).
	Therefore we prefer to state all the results concerning hyperbolicity for the class of length spaces.
	To compensate for the absence of geodesics we use the following property.
	For every $x,x' \in X$, for every $\eta >0$, for every $ 0 \leq l \leq \dist x{x'}$ there exists $y \in X$ such that $\dist xy=l$ and $\gro x{x'}y \leq \eta$.
	In this context, the Gromov product $\gro x{x'}y$ should be thought as an analogue of the distance between $y$ and a ``geodesic'' joining $x$ and $x'$.
	
	\paragraph{}
	If $X$ is 0-hyperbolic, then it can be isometrically embedded in an $\R$-tree \cite[Chap. 2, Prop. 6]{GhyHar90}.
	However we will always assume that the hyperbolicity constant $\delta$ is positive.
	Most of the results hold for $\delta = 0$.
	But this is more convenient to define particular subsets (see for instance Definition~\ref{def: hull} or Lemmas~\ref{res: intersection of quasi-convex} and \ref{res: intersection of thickened quasi-convex}) without introducing other auxiliary positive parameters.

	 \paragraph{}
	The hyperbolicity constant of the hyperbolic plane $\H_2$ will play a particular role.
	Therefore we denote it by $\boldsymbol \delta$ (bold delta).
		
	\paragraph{}
	From now on we assume that the space $X$ is $\delta$-hyperbolic.
	It is known that triangles in a hyperbolic geodesic space are thin (every side lies in a uniform neighborhood of the union of the two others). 
	Since our space is not necessarily geodesic, we use instead the following metric inequalities.
	They will be used in situations where the Gromov products $\gro xys$, $\gro xyt$ or $\gro xzt$ are ``small''.
	Their proof is left to the reader.
	\begin{lemm}
	\label{res: metric inequalities}
		Let $x$, $y$, $z$, $s$ and $t$ be five points of $X$.
		\begin{enumerate}
			\item \label{enu: metric inequalities - thin triangle}
			\begin{math}
				\gro xyt \leq \max \left\{\fantomB \dist xt - \gro yzx , \gro xzt  \right\} + \delta,
			\end{math}
			\item \label{enu: metric inequalities - two points close to a geodesic}
			\begin{math}
				\dist st \leq \left|\fantomB \dist xs - \dist xt \right| + 2\max\left\{\fantomB \gro xys, \gro xyt\right\} + 2\delta,
			\end{math}
			\item \label{enu: metric inequalities - comparison tripod} The distance $\dist st$ is bounded above by 
			\begin{displaymath}
				\max\left\{\fantomB \dist{\fantomB \dist xs}{\dist xt} + 2\max \left\{ \gro xys,\gro xzt\right\}, \dist xs +\dist xt - 2 \gro yzx \right\} + 4\delta.
			\end{displaymath}
		\end{enumerate}
	\end{lemm}
		
%	\begin{proof} Let $x$, $y$, $z$, $s$ and $t$ be points of $X$.
%		\paragraph{Point~\ref{enu: metric inequalities - thin triangle}} By hyperbolicity we have $\min\left\{\fantomB \gro xyt, \gro yzt \right\} \leq \gro xzt + \delta$.
%		However $\gro zyt = \gro xyt - \dist xt + \gro yzx + \gro xzt$.
%		Therefore
%		\begin{displaymath}
%			\min\left\{\fantomB \gro xyt,\gro xyt - \dist xt + \gro yzx + \gro xzt\right\} \leq \gro xzt + \delta,
%		\end{displaymath}
%		which is equivalent to Point~\ref{enu: metric inequalities - thin triangle}.
%		
%		\paragraph{Point~\ref{enu: metric inequalities - two points close to a geodesic}}
%		By hyperbolicity one has
%		\begin{eqnarray*}
%			\dist st & = & \dist xs + \dist xt - 2 \gro stx  \\
%				& \leq & \dist xs + \dist xt - 2 \min \left\{ \fantomB \gro syx,  \gro ytx \right\} + 2\delta.
%		\end{eqnarray*}
%		However $\gro syx = \dist xs - \gro xys$ and $\gro ytx = \dist xt - \gro xyt$ which leads to Point~\ref{enu: metric inequalities - two points close to a geodesic}.
%		\paragraph{Point~\ref{enu: metric inequalities - comparison tripod}}
%		By hyperbolicity one has
%		\begin{eqnarray*}
%			\dist st & = & \dist xs + \dist xt - 2 \gro stx  \\
%				& \leq & \dist xs + \dist xt - 2 \min \left\{ \fantomB \gro syx, \gro yzx, \gro ztx \right\} + 4\delta.
%		\end{eqnarray*}
%		However $\gro syx = \dist xs - \gro xys$ and $\gro tzx = \dist xt - \gro xzt$ which leads to Point~\ref{enu: metric inequalities - comparison tripod}.
%	\end{proof}

\subsection{Quasi-geodesics}
\label{sec: quasi-geodesics}
	
	\begin{defi}
	\label{def: quasi-geodesic}
		Let $l \geq 0$, $k \geq 1$ and $L\geq 0$.
		Let $I$ be an interval of $\R$.
		A path $\gamma : I \rightarrow X$ is
		\begin{enumerate}
			\item
			a \emph{$(k,l)$-quasi-geodesic} if for all $t,t' \in I$
			\begin{displaymath}
				\dist{\gamma(t)}{\gamma(t')} \leq \dist t{t'} \leq k\dist{\gamma(t)}{\gamma(t')}+l.
			\end{displaymath}
			\item
			an \emph{$L$-local $(k,l)$-quasi-geodesic} if the restriction of $\gamma$ to any interval of length $L$ is a $(k,l)$-quasi-geodesic.
		\end{enumerate}
	\end{defi}
	
	\rems The first inequality in the definition just follows from the fact that $\gamma$ is parametrized by arc length.
	Since $X$ is a length space, for every $x,x' \in X$, for every $l>0$, there exists a $(1,l)$-quasi-geodesic joining $x$ and $x'$ (take a rectifiable path between them whose length is shorter than $\dist x{x'} +l$).
	Let $\gamma : I \rightarrow X$ be a $(1,l)$-quasi-geodesic.
	Let $x = \gamma(t)$, $x'=\gamma(t')$ and $y=\gamma(s)$ be three points on $\gamma$.
	If $t\leq s\leq t'$ then $\gro x{x'}y\leq l/2$.
	
	\begin{prop}
	\label{res: quasi-convexity of quasi-geodesics}
		Let  $\gamma : I \rightarrow X$ be a $(1,l)$-quasi-geodesic of $X$.
		\begin{enumerate}
			\item \label{enu: quasi-convexity of quasi-geodesics - gromov product}
			Let $x$ be a point of $X$ and $p$ an $\eta$-projection of $x$ on $\gamma$.
			For all $y \in \gamma$, $\gro xyp \leq l + \eta + 2\delta$.
			\item \label{enu: quasi-convexity of quasi-geodesics - quasi-convex}
			For every $x \in X$ for every $y,y' \in \gamma$ we have $\gro y{y'}x -l \leq d(x,\gamma) \leq \gro y{y'}x +l + 3\delta$.
		\end{enumerate}
	\end{prop}
	
	\begin{proof}
		Let $x$ be a point of $X$.
		We denote by $p=\gamma(s)$ an $\eta$-projection of $x$ on $\gamma$ and $y = \gamma(t)$ a point of $\gamma$.
		By reversing if necessary the parametrization of $\gamma$ we can assume that $s \leq t$.
		Note that $\dist yp \geq \gro xyp$.
		Therefore there exists $r \in \intval st$ such that the point $z = \gamma (r)$ satisfies  $\dist zp = \gro xyp$ and $\gro pyz \leq l/2$.
		By Lemma~\ref{res: metric inequalities}-\ref{enu: metric inequalities - thin triangle} we have $\gro xpz \leq l/2 + \delta$ i.e.,
		\begin{displaymath}
			\gro xyp = \dist pz  \leq \dist xp -   \dist xz + l + 2 \delta.
		\end{displaymath}
		Nevertheless, $p$ is an $\eta$-projection of $x$ on $\gamma$. 
		Thus $ \dist xp -   \dist xz \leq \eta$.
		It follows that $\gro xyp \leq l + \eta + 2\delta$, which proves \ref{enu: quasi-convexity of quasi-geodesics - gromov product}.
		
		\paragraph{} Let $x \in X$. 
		Let $y$ and $y'$ be two points of $\gamma$.
		The first inequality of \ref{enu: quasi-convexity of quasi-geodesics - quasi-convex} follows from the triangle inequality.
		We denote by $p$ a projection of $x$ on $\gamma$.
		It follows from Point~\ref{enu: quasi-convexity of quasi-geodesics - gromov product} that $d(x,\gamma) =\gro ypx + \gro xyp$ is bounded above by  $\gro ypx + l + 2 \delta$.
		In the same way $d(x,\gamma) \leq \gro {y'}px + l +2 \delta$.
		Consequently, the hyperbolicity condition~(\ref{eqn: hyperbolicity condition 1}) gives
		\begin{displaymath}
			d(x,\gamma) \leq \min \left\{ \gro ypx, \gro {y'}px \right\} + l + 2 \delta \leq \gro y{y'}x + l + 3 \delta. \qedhere
		\end{displaymath}
	\end{proof}
	
	\begin{prop}[Stability of quasi-geodesics {\cite[Chap. 3, Th. 1.2-1.4]{CooDelPap90}}]
	\label{res: stability quasi-geodesic}
		Let $k \geq 1$, $k' >k$ and $l \geq 0$.
		There exists $L\geq 0$, $D \geq 0$ which only depend on $\delta$, $k$, $k'$ and $l$ such that the Hausdorff distance between two $L$-local $(k,l)$-quasi-geodesics joining the same endpoints is at most  $D$.
		Moreover, every $L$-local $(k,l)$-quasi-geodesic is also a (global) $(k',l)$-quasi-geodesic.
	\end{prop}
	
	\paragraph{}
	In this paper we are mostly using $L$-local $(1,l)$-quasi-geodesics.
	For these paths one can easily provide a precise value for $D$ (see Corollary~\ref{res: stability (1,l)-quasi-geodesic}).
	This is not really crucial but it will decrease the number of parameters that we have to deal with in all the proofs.
	
	\begin{coro}
	\label{res: stability (1,l)-quasi-geodesic}
		Let $l,l' \geq 0$.
		There exists $L=L(l,l',\delta)$ which only depends on $\delta$, $l$ and $l'$ with the following property.
		Let $\gamma : [a,b]\rightarrow X$ be a $L$-local $(1,l)$-quasi-geodesic and $\gamma : [a',b']\rightarrow X$ a $L$-local $(1,l')$-quasi-geodesic.
		If they join the same extremities then $\gamma$ lies in the $(l/2+ l' + 5 \delta)$-neighborhood of $\gamma'$.		
	\end{coro}
	
	\begin{proof}
		Let $l \geq 0$.
		The constants $L$, $D$ and $k$ are chosen so that we can apply Proposition~\ref{res: stability quasi-geodesic} to $\gamma$ and $\gamma'$:
		the Hausdorff distance between $\gamma$ and $\gamma'$ is at most $D$ and $\gamma'$ is a $(k,l')$-global quasi-geodesic. 
		Note that $k$, $D$ and $L$ only depend on $l$, $l'$ and $\delta$.
		Without loss of generality we may assume that $L > k(4D + 3l+ 6\delta) +l'$.
		Let $x = \gamma(t)$ be a point on $\gamma$.
		We assume that $\dist at$ and $\dist bt$ are bounded below by $D + 3l/2 + 3\delta$. 
		The other cases are similar.
		Let us put $t_\pm = t \pm(D + 3l/2 + 3\delta)$ and $x_\pm = \gamma(t_\pm)$.
		By the stability of quasi-geodesics there exist points $y_\pm = \gamma'(s_\pm)$ such that $\dist{x_\pm}{y_\pm} \leq D$.
		Without loss of generality we can assume that $s_- \leq s_+$.
		By hyperbolicity we obtain
		\begin{displaymath}
			\min\left\{ \dist {x_-}x - \dist {x_-}{y_-}, \gro {y_-}{y_+}x, \dist {x_+}x - \dist {x_+}{y_+}\right\} \leq \gro {x_-}{x_+}x + 2 \delta \leq l/2 + 2 \delta.
		\end{displaymath}
		However by construction $\dist{x_\pm}x$ is bounded below by $D + l/2 + 3\delta$ and $\dist{x_\pm}{y_\pm}$ above by $D$.
		Therefore, we necessarily have $\gro {y_-}{y_+}x \leq l/2 + 2 \delta$.
		We now claim that $\dist{s_+}{s_-} \leq L$.
		Indeed $\gamma'$ is a $(k,l')$-quasi-geodesic.
		Therefore the triangle inequality leads to 
		\begin{displaymath}
			k^{-1}\left(\fantomB\dist{s_+}{s_-} - l' \right)\leq \dist {y_+}{y_-} \leq \dist {x_+}{x_-} + 2D \leq 4D + 3l + 6\delta
		\end{displaymath}
		Our claim follows from the assumption on $L$.
		In particular $\gamma'$ restricted to $\intval{s_-}{s_+}$ is a (global) $(1,l)$-quasi-geodesic.
		By Proposition~\ref{res: quasi-convexity of quasi-geodesics}-\ref{enu: quasi-convexity of quasi-geodesics - quasi-convex}, the point $x$ is  $(l/2 + l'+ 5\delta)$-close to $\gamma'$.
	\end{proof}

	\rems We keep the notations of the corollary.
	Let $x = \gamma(t)$, $x'=\gamma(t')$ and $y=\gamma(s)$ be three points on $\gamma$.
	If $t \leq s \leq t'$ then $\gro x{x'}y \leq l/2 + 5\delta$.
	Moreover, for every $p \in X$ we have $d(p, \gamma) \leq \gro x{x'}p + l+8\delta$.
	
	\begin{prop}[Stability of discrete quasi-geodesics]
	\label{res: stability discrete quasi-geodesic}
		Let $l > 0$.
		There exists $L = L(l,\delta)$ which only depends on $\delta$ and $l$ with the following property.
		Let $x_0, \dots, x_m$ be a sequence of points of $X$ such that
		\begin{enumerate}
			\item for every $i \in \intvald 1{m-1}$, $\gro {x_{i-1}}{x_{i+1}}{x_i} \leq l$,
			\item for every $i \in \intvald 1{m-2}$, $\dist {x_{i+1}}{x_i} \geq L$.
		\end{enumerate}
		Then for all $i \in \intvald 0m$, $\gro x{x'}{x_i} \leq l+ 5\delta$.
		Moreover, for all $p \in X$ there exists $i \in \intvald 0{m-1}$ such that $\gro{x_{i+1}}{x_i}p \leq \gro {x_0}{x_m}p + 2l + 8\delta$.
	\end{prop}
	
	\begin{proof}
		Let $\eta > 0$.
		For every $i \in \intvald 0{m-1}$ we choose a $(1,\eta)$-quasi-geodesic $\gamma_i$ joining $x_i$ to $x_{i+1}$.
		We denote by $\gamma$ the concatenation $\gamma_0 \dots \gamma_{m-1}$.
		According to our assumptions this is a $L$-local $(1,2l+4\eta)$-quasi-geodesic.
		We can therefore apply Proposition~\ref{res: stability quasi-geodesic}.
		If $L$ is sufficiently large then we obtain the followings.
		\begin{enumerate}
			\item For every $i \in \intvald 0m$ we have $\gro x{x'}{x_i} \leq  l+2\eta + 5\delta$.
			\item For every $p \in X$ we have $d(p,\gamma) \leq \gro{x_0}{x_m}p  + 2l + 4 \eta + 8 \delta$.
			However there exists $i \in \intvald 0{m-1}$ such that the distance between $p$ and $\gamma$ is exactly $d(p, \gamma_i)$.
			The path $\gamma_i$ being a $(1, \eta)$-quasi-geodesic, we get $\gro {x_i}{x_{i+1}}p \leq \gro{x_0}{x_m}p + 2l + 5\eta + 8 \delta$.
		\end{enumerate}
		The inequalities that we obtained hold for every sufficiently small $\eta>0$, which gives the desired conclusion.
	\end{proof}

\subsection{Quasi-convex subsets}

	\begin{defi}
	\label{def: quasi-convex}
		Let $\alpha \geq 0$.
		A subset $Y$ of $X$ is \emph{$\alpha$-quasi-convex} if for every $x \in X$, for every $y,y' \in Y$, $d(x,Y) \leq \gro y{y'}x + \alpha$.
	\end{defi}
	
	\rem Since $X$ is not a geodesic space our definition of quasi-convex slightly differs from the usual one (every geodesic joining two points of $Y$ remains in the $\alpha$-neighborhood of $Y$).
	However if $X$ is geodesic, an $\alpha$-quasi-convex subset in the usual sense is $(\alpha + 3\delta)$-quasi-convex in our sense and conversely.
	According to Proposition~\ref{res: quasi-convexity of quasi-geodesics} every $(1,l)$-quasi-geodesic is $(l+3\delta)$-quasi-convex.
	If $L$ is sufficiently large then every $L$-local $(1,l)$-quasi-geodesic is $(l + 8 \delta)$-quasi-convex (see Corollary~\ref{res: stability (1,l)-quasi-geodesic}).
	For our purpose we will also need a slightly stronger version of quasi-convexity.
	
	\begin{defi}
		Let $Y$ be a subset of $X$ connected by rectifiable paths.
		The length metric on $Y$ induced by the restriction of $\distV[X]$ to $Y$ is denoted by $\distV[Y]$.
		We say that $Y$ is \emph{strongly quasi-convex} if $Y$ is $2\delta$-quasi-convex and for every $y,y' \in Y$,
		\begin{displaymath}
			\dist[X]y{y'} \leq \dist[Y]y{y'} \leq \dist[X]y{y'} + 8\delta.
		\end{displaymath}
	\end{defi}
	
	\rem The first inequality is just a consequence of the definition of $\distV[Y]$.
		
%	\begin{prop}
%	\label{res: strong quasi-convex implies quasi-convex}
%		Let $\alpha \geq 0$.
%		Let $Y$ be a subset of $X$.
%		If $Y$ is $\alpha$-strongly quasi-convex then it is also $(\alpha +3\delta)$-quasi-convex.
%	\end{prop}
%	
%	\begin{proof}
%		Let $x \in X$ and $y,y' \in Y$.
%		Let $\eta >0$.
%		Since $Y$ is $\alpha$-strongly convex there is a path $\gamma$ contained in $Y$ joining $y$ to $y'$ whose length is bounded above by $\dist[X]y{y'} + \alpha +\eta$.
%		In particular this path is an $(1,\alpha +\eta)$-quasi-geodesic.
%		Therefore by Proposition~\ref{res: quasi-convexity of quasi-geodesics}
%		\begin{displaymath}
%			d(x,Y) \leq d(x,\gamma) \leq \gro y{y'}x +\alpha + 3\delta + \eta.
%		\end{displaymath}
%		This inequality holds for every $\eta>0$, hence $Y$ is $(\alpha +3\delta)$-quasi-convex.
%	\end{proof}
	
	\begin{lemm}[Projection on a quasi-convex. Compare {\cite[Chap. 10, Prop. 2.1]{CooDelPap90}}]
	\label{res: proj quasi-convex}
		Let $Y$ be an $\alpha$-quasi-convex subset of $X$. 
		Let $x, x' \in X$.
		\begin{enumerate}
			\item \label{enu: proj quasi-convex - gromov product}
			If $p$ is an $\eta$-projection of $x$ on $Y$, then for all $y \in Y$, $\gro xyp \leq \alpha + \eta$.
			\item \label{enu: proj quasi-convex - distance two points }
			If $p$ and $p'$ are respective $\eta$- and $\eta'$-projections of $x$ and $x'$ on $Y$, then 
			\begin{displaymath}
				\dist p{p'} \leq \max \left\{\fantomB \dist x{x'}-\dist xp - \dist {x'}{p'} +2\epsilon, \epsilon \right\},
			\end{displaymath}
			where $\epsilon = 2 \alpha + \eta + \eta' + \delta$.
		\end{enumerate}
	\end{lemm}
	
%	\begin{proof}
%		Let $y \in X$ and $p$ be an $\eta$-projections of $x$ on $Y$.
%		Applying the definition of a quasi-convex, $\dist xp \leq d(x,Y)+\eta \leq \gro ypx + \alpha + \eta$.
%		Therefore $\gro xyp = \dist xp - \gro ypx \leq \alpha + \eta$, which proves \ref{enu: proj quasi-convex - gromov product}.
%		Let us now consider $p$ and $p'$ respective $\eta$- and $\eta'$ projections of $x$ and $x'$ on $Y$.
%		Applying twice \ref{enu: proj quasi-convex - gromov product}, we get
%		\begin{displaymath}
%		\dist xp + \dist {x'}{p'} + 2\dist p{p'}  \leq \dist x{p'} + \dist {x'}p +4\alpha +2\eta +2\eta'
%		\end{displaymath}
%		Therefore the hyperbolicity condition (\ref{eqn: hyperbolicity condition 2}) leads to
%		\begin{displaymath}
%			\dist xp + \dist {x'}{p'} + 2\dist p{p'}
%			\leq \max\left\{\fantomB \dist xp + \dist {x'}{p'} , \dist x{x'} +\dist p{p'} \right\} + 2\epsilon,
%		\end{displaymath}
%	 	where $\epsilon = 2 \alpha + \eta + \eta' + \delta$.
%		Point \ref{enu: proj quasi-convex - distance two points } follows from this inequality.
%	\end{proof}

	\begin{lemm}[Neighborhood of a quasi-convex. Compare {\cite[Chap. 10, Prop. 1.2]{CooDelPap90}}]
	\label{res: neighborhood of a quasi-convex}
		Let $\alpha \geq 0$.
		Let $Y$ be an $\alpha$-quasi-convex subset of $X$.
		For every $A \geq \alpha$, the $A$-neighborhood of $Y$ is $2 \delta$-quasi-convex.
	\end{lemm} 
	
%	\begin{proof}
%		Let $x$ be a point of $X$ and $y,y'$ points in the $A$-neighborhood of $Y$.
%		Let $\eta >0$.
%		We denote by $p$ and $p'$ respective $\eta$-projection of $y$ and $y'$ on $Y$.
%		By hyperbolicity 
%		\begin{displaymath}
%			\min\left\{\fantomB \gro ypx, \gro p{p'}x, \gro {p'}{y'}x \right\} \leq \gro y{y'}x + 2 \delta.
%		\end{displaymath}
%		Since $Y$ is $\alpha$-quasi-convex, $\gro p{p'}x \geq d(x,Y) -\alpha \geq  d(x,Y) - A$.
%		On the other hand, the triangle inequality gives, 
%		\begin{displaymath}
%			\gro ypx \geq \dist xp - \dist yp \geq d(x,Y) - A - \eta .
%		\end{displaymath}
%		In the same way $\gro {y'}{p'}x \geq d(x,Y) - A - \eta$.
%		Hence $d(x,Y) \leq \gro y{y'}x + A + \eta + 2\delta$.
%		This inequality holds for every $\eta >0$.
%		Consequently $d(x,Y) \leq \gro y{y'}x + A + 2\delta$.
%		However $X$ is a length-space
%		Thus the distance of $x$ to the $A$-neighborhood of $Y$ is at most $\gro y{y'}x + 2\delta$.
%	\end{proof}
		
	\begin{lemm}
	\label{res: intersection of quasi-convex}
		Let $Y_1$ and $Y_2$ be respectively $\alpha_1$- and $\alpha_2$-quasi-convex subsets of $X$.
		The subset
		\begin{displaymath}
			 Z = Y_1^{+ \alpha_1 +3\delta} \cap Y_2^{+ \alpha_2 + 3\delta}
		\end{displaymath}
		is $7\delta$-quasi-convex. 
	\end{lemm}
	
	\begin{proof}
		Let $x \in X$ and $t,t' \in Z$.
		Let $\eta \in (0, 2\delta)$.
		We denote by $\gamma$ a $(1,\eta)$-quasi-geodesic joining $t$ and $t'$.
		Let $u$ be a point of $\gamma$.
		In particular $\gro t{'t}u \leq \eta/2$.
		We assume first that $\dist tu \geq 4\delta + 2 \eta$ and $\dist {t'}u \geq  4\delta + 2 \eta$.
		We denote by $y$ and $y'$ respective $\eta$-projections of $t$ and $t'$ on $Y_i$.
		By hyperbolicity 
		\begin{equation}
		\label{eqn: intersection of quasi-convex}
			\min \left\{\fantomB \gro tyu, \gro y{y'}u, \gro{y'}{t'}u \right\} \leq \gro t{t'}u + 2 \delta \leq \frac  12\eta + 2\delta.
		\end{equation}
		Assume that the minimum is achieved by $\gro tyu$, it gives
		\begin{displaymath}
			\dist yu = \dist ty - \dist tu + 2 \gro tyu \leq \alpha_i + 3\delta
		\end{displaymath}
		Thus $u$ belongs to the $(\alpha_i +3\delta)$-neighborhood of $Y_i$.
		The same holds if the minimum in (\ref{eqn: intersection of quasi-convex}) is achieved by $\gro {y'}{t'}u$.
		Suppose now that the minimum is achieved by $\gro y{y'}u$.
		The set $Y_i$ being $\alpha$-quasi-convex we get $d(u,Y_i) \leq \gro y{y'}u + \alpha_i \leq \alpha_i + 3 \delta$.
		In all cases $u$ lies in the $(\alpha_i + 3\delta)$-neighborhood of $Y_i$.
		It follows that any point of $\gamma$ is in the $(4\delta + 2\eta)$-neighborhood of $Z$.
		
		\paragraph{}
		Let $p$ be a projection of $x$ on $\gamma$.
		The path $\gamma$ is $(\eta+ 3\delta)$-quasi-convex (Proposition~\ref{res: quasi-convexity of quasi-geodesics}) hence $\dist xp \leq \gro t{t'}x + \eta + 3\delta$.
		According to the previous remark $p$ lies in the $(4\delta + 2 \eta)$-neighborhood of $Z$ thus $d(x,Z) \leq \gro t{t'}x + 3\eta + 7\delta$.
		This inequality holds for every sufficiently small $\eta$. 
		Thus $Z$ is $7\delta$-quasi-convex.
	\end{proof}
	
	\begin{lemm}
	\label{res: intersection of thickened quasi-convex}
		Let $Y$ and $Z$ be respectively $\alpha$- and $\beta$-quasi-convex subsets of $X$.
		For all $A \geq 0$ we have
		\begin{displaymath}
			\diam \left( Y^{+A} \cap Z^{+A} \right) 
			\leq \diam \left( Y^{+\alpha+3\delta} \cap Z^{+\beta+3\delta} \right) +2A + 4\delta.
		\end{displaymath}
	\end{lemm}
	
	\begin{proof}
		Let $x$ and $x'$ be two points of $Y^{+A} \cap Z^{+A}$.
		We assume that $\dist x{x'} > 2A + 4\delta$.
		Let $\eta \in (0, \delta)$ such that $\dist x{x'} > 2A + 4\delta+ 6\eta$.
		There exist $t,t' \in X$ such that $\dist xt = \dist {x'}{t'} = A + 2\delta + 3\eta$ and $\gro x{x'}t, \gro x{x'}{t'} \leq \eta$.
		Note that $\dist x{t'} , \dist {x'}t \geq A + 2\delta + 3\eta$.
		We claim that $t$ and $t'$ belong to the $(\alpha+3\delta)$-neighborhood of $Y$.
		Let us denote by $y$ and $y'$ respective $\eta$-projection of $x$ and $x'$ on $Y$.
		By hyperbolicity 
		\begin{displaymath}
		\min\left\{ \dist xt - \dist xy , \gro y{y'}t, \dist{x'}t - \dist{x'}{y'} \right\} \leq \gro x{x'}t + 2 \delta \leq  2\delta + \eta.
		\end{displaymath}
		It follows that $\gro y{y'}t \leq 2\delta+ \eta \leq 3 \delta$.
		The subset $Y$ being $\alpha$-quasi-convex we get $d(t, Y) \leq \alpha + 3\delta$.
		The same  holds for $t'$, which proves our claim.
		Similarly $t$ and $t'$ lie in the $(\beta + 3\delta)$-neighborhood of $Z$.
		Consequently, 
		\begin{displaymath}
			\dist x{x'} 
			\leq \dist t{t'} + 2A + 4\delta + 6\eta
			\leq \diam \left( Y^{+\alpha+3\delta} \cap Z^{+\alpha+3\delta} \right) +2A + 4\delta + 6\eta.
		\end{displaymath}
		This last inequality actually holds for every sufficiently small $\eta$ and $x$, $x'$ in $Y^{+A} \cap Z^{+A}$, which leads to the conclusion.
	\end{proof}
		
	\begin{defi}
	\label{def: hull}
		Let $Y$ be a subset of $X$.
		The \emph{hull} of $Y$ denoted by $\hull Y$ is the union of all $(1, \delta)$-quasi-geodesics joining two points of $Y$.
	\end{defi}

	\begin{lemm}
	\label{res: hull quasi-convex}
		Let $Y$ be a subset of $X$. 
		The hull of $Y$ is $6\delta$-quasi-convex.
	\end{lemm}
	
	\begin{proof}
		Let $x \in X$ and $y,y' \in \hull Y$.
		By definition there exist $\gamma : \intval ab \rightarrow X$ and $\gamma' : \intval {a'}{b'} \rightarrow X$ two $(1, \delta)$-quasi-geodesics joining points of $Y$ such that $y$ and $y'$ respectively lie on $\gamma$ and $\gamma'$.
		Since $X$ is a length space, there exists a $(1,\delta)$-quasi-geodesic $\gamma_0$ between $\gamma(a)$ and $\gamma'(a')$.
		In particular $\gamma_0 \subset \hull Y$.
		By hyperbolicity 
		\begin{displaymath}
			\min \left\{ \fantomB \gro y{\gamma(a)}x, \gro {\gamma(a)}{\gamma(a')}x, \gro{\gamma(a')}{y'}x \right\} \leq \gro y{y'}x + 2 \delta.
		\end{displaymath}
		However $\gamma$ is $4\delta$-quasi-convex (Proposition~\ref{res: quasi-convexity of quasi-geodesics}), thus $d(x, \hull Y) \leq d(x,\gamma) \leq \gro y{\gamma(a)}x + 4\delta$.
		We have similar inequalities for $\gamma_0$ and $\gamma'$.
		Hence $d(x,\hull Y)  \leq \gro y{y'}x + 6 \delta$.
	\end{proof}
	
	\begin{lemm}
	\label{res: gromov product and hull}
		Let $Y$ and $Z$ be two subsets of $X$.
		Let $x$ be a point of $X$.
		Assume that for all $y \in Y$, for all $z \in Z$, $\gro yzx \leq \alpha$.
		Then for all $y \in \hull Y $, for all $z \in \hull Z$, $\gro yzx \leq \alpha + 3\delta$.
	\end{lemm}
	
	\begin{proof}
		Let $y \in \hull Y$ and $z \in \hull Z$.
		By definition there exists $y_1, y_2 \in Y$ (\resp $z_1,z_2 \in Z$) such that $y$ (\resp $z$) lies on a $(1,\delta)$-quasi-geodesic between $y_1$ and $y_2$ (\resp $z_1$ and $z_2$).
		By hyperbolicity
		\begin{displaymath}
			\min\left\{ \gro {y_1}xy, \gro{y_2}xy \right\} \leq \gro {y_1}{y_2}y + \delta \leq \frac32 \delta.
		\end{displaymath}
		In particular there is $i \in \{1,2\}$ such that $\gro {y_i}xy \leq 3\delta/2$.
		In the same way there is $j \in \{1,2\}$ such that $\gro {z_j}xz \leq  3\delta/2$.
		By triangle inequality we obtain
		\begin{displaymath}
			\gro yzx \leq \gro {y_i}{z_j}x + \gro {y_i}xy +  \gro {z_j}xz \leq \alpha + 3\delta. \qedhere
		\end{displaymath}
	\end{proof}

\subsection{Ultra-limit of hyperbolic spaces}

	\paragraph{} Let us first recall the definition of the ultra-limit of a sequence of metric spaces and some related notations.
	
	\paragraph{} A non-principal ultra-filter is a finite additive map $\omega : \mathcal P(\N) \rightarrow \{0,1\}$ such that $\omega(\N) = 1$ and which vanishes on every finite subset of $\N$.
	A property $P_n$ is true \oas if $\omega \left( \left\{ n \in \N | P_n \text{ is true} \right\} \right) = 1$.
	A real sequence $\left(u_n\right)$ is \oeb if there exists $M$ such that $\left|u_n\right| \leq M$ \oas.
	Given $l \in \R$, we say that the $\omega$-limit of $\left(u_n\right)$ is $l$ and write $\limo u_n = l$ if for all $\epsilon >0$, $\dist{u_n}l \leq \epsilon$ \oas.
	In particular, any sequence which is \oeb admits a $\omega$-limit \cite{Bou71}.
	
	\paragraph{} Let $\left(X_n, x_n^0 \right)$ be a sequence of pointed metric spaces. We define the following set
	\begin{displaymath}
		\Pi_\omega X_n =\set{\left(x_n\right) \in \Pi_{n\in \N} X_n }{ \left(\dist{x_n}{x_n^0}\right) \text{ is \oeb} }.
	\end{displaymath}
	The space $\Pi_\omega X_n$ is endowed with a pseudo-metric defined in the following way:
	\begin{math}
		\dist{(x_n)}{(y_n)} = \limo \dist {x_n}{y_n}.
	\end{math}
	
	\begin{defi}
		The \emph{$\omega$-limit} of $\left(X_n, x_n^0\right)$, denoted by $\limo \left(X_n, x_n^0\right)$ or simply $\limo X_n$, is the quotient of the space $\Pi_\omega X_n$ by the equivalence relation which identifies two points at distance zero.
		The pseudo-distance on $\Pi_\omega X_n$ induces a distance on $\limo X_n$.
	\end{defi}
	
	\rem If the diameter of $X_n$ is uniformly bounded, then $\limo \left(X_n, x_n^0\right)$ does not depend on the choice of a base point $x_n^0$. 
	
	\notas
	\begin{itemize}
		\item Given a sequence $\left(x_n \right) \in \Pi_\omega X_n$ we write $\limo x_n$ for its equivalence class in $\limo X_n$.
		\item For all $n \in \N$, let $Y_n$ be a subset of $X_n$. The set $\limo Y_n$ is defined by 
		\begin{displaymath}
			\limo Y_n = \set{ \limo y_n }{ \left(y_n\right) \in \Pi_\omega X_n \text{ and } y_n \in Y_n \text{ \oas}}.
		\end{displaymath}
	\end{itemize}
	
	\begin{prop}[Ultra-limit of hyperbolic spaces [{\cite[Prop. 1.1.2]{Coulon:il}}]]
	\label{res: ultra-limit of hyperbolic spaces} 
		Let $\omega$ be a non-principal ultra-filter.
		Let $\left(X_n, x_n^0\right)$ be a sequence of pointed $\delta_n$-hyperbolic length spaces such that $\delta = \limo \delta_n$.
		Then $\limo X_n$ is a $\delta$-hyperbolic geodesic space. 
		In particular, if $\delta = 0$, $\limo X_n$ is an $\R$-tree.
	\end{prop}
	
	\begin{prop}[{\cite[Prop. 1.1.4]{Coulon:il}}]
	\label{res: ultra limit of spaces which is hyp}
		Let $\delta \geq 0$.
		Let $\omega$ be a non-principal ultra-filter.
		Let $\left(X_n, x_n^0\right)$ be a sequence of pointed length spaces.
		Assume that $\limo X_n$ is a $\delta$-hyperbolic space. 
		Then for every $\eta>0$ for every $r>0$, every ball of radius of $r$ of $X_n$ is $(\delta+\eta)$-hyperbolic \oas.
	\end{prop}

	\begin{prop}
	\label{res: ultra-limit and diaminter}
		Let $\omega$ be a non-principal ultra-filter.
		Let $\left(X_n, x_n^0\right)$ be a sequence of pointed $\delta_n$-hyperbolic length spaces with $\limo \delta_n = 0$.
		For every $n \in \N$ let $Y_n$ and $Z_n$  be respectively $\alpha_n$- and $\beta_n$-quasi-convex subsets of $X_n$.
		We denote by $Y= \limo Y_n$ and $Z = \limo Z_n$ the corresponding limit subsets of $X = \limo X_n$.
		Then 
		\begin{displaymath}
			\diam\left(Y\cap Z\right) \leq \limo \diam\left(Y_n^{+\alpha_n + 3 \delta_n} \cap Z_n^{+ \beta_n + 3 \delta_n}\right).
		\end{displaymath}

	\end{prop}
	
	\begin{proof}
		Let $x = \limo x_n$ and $x' = \limo x'_n$ be two points of $Y \cap Z$.
		Let $A >0$.
		Since $x$ and $x'$ belong to both $Y$ and $Z$, $x_n$ and $x'_n$ belong to $Y_n^{+A}\cap Z_n^{+ A}$ \oas.
		Applying Proposition~\ref{res: intersection of thickened quasi-convex} we obtain
		\begin{displaymath}
			\dist{x'_n}{x_n} \leq  \diam\left(Y_n^{+\alpha_n + 3 \delta_n} \cap Z_n^{+ \beta_n + 3 \delta_n}\right) + 2A + 4\delta_n \text{ \oas}.
		\end{displaymath}
		After taking the $\omega$-limit, it gives
		\begin{displaymath}
			\dist {x'}x \leq \limo \diam\left(Y_n^{+\alpha_n + 3 \delta_n} \cap Z_n^{+ \beta_n + 3 \delta_n}\right) +2A.
		\end{displaymath}
		This last inequality holds for every $x,x' \in Y\cap Z$ and every $A>0$, which leads to the result.
	\end{proof}

\subsection{Isometries of a hyperbolic space}
\label{sec: isometries hyperbolic space}

	In this section we assume that the space $X$ is geodesic and proper.
	By \emph{proper} we mean that every closed ball of $X$ is compact.
	Although it is not necessarily unique, $\geo x{x'}$ stands for a geodesic between two points $x$ and $x'$ of $X$.
	We denote by $\partial X$ the boundary at infinity of $X$ (see \cite[Chap. 2]{CooDelPap90}).
	The space $X$ being proper any two distinct points of $\partial X$ are joined by a bi-infinite geodesic.
	In this situation one can precise the constants that appears in Corollary~\ref{res: stability (1,l)-quasi-geodesic} (see \cite[Chap. III.H, Th. 1.13]{BriHae99}):
	the Hausdorff distance between two $200\delta$-local $(1,0)$-quasi-geodesics of $X$ joining the same extremities (eventually in $\partial X$) is at most $5\delta$.
	In particular the Hausdorff distance between two bi-infinite geodesics joining the same points of $\partial X$ is at most $5\delta$.
	
	\begin{lemm}
	\label{res: strong neighborhood of a quasi-convex}
		Let $\alpha \geq 0$ and $Y$ be an $\alpha$-quasi-convex subset of $X$.
		For every $A \geq \alpha + 2\delta$, the $A$-neighborhood of $Y$ is strongly quasi-convex.
	\end{lemm}

	\rem An analog statement is true if the space $X$ is not proper of geodesic.
	However is requires to consider the \emph{open} $A$-neighborhood of $X$.
	That is why we preferred to state this result with this stronger assumption.
	
	\begin{proof}
		According to Lemma~\ref{res: neighborhood of a quasi-convex}, it is sufficient to prove that the $2\delta$-neighborhood of a closed $2\delta$-quasi-convex subset $Y$ of $X$ is $8\delta$-strongly quasi-convex.
		Let $x$ and $x'$ be two points of $X$ which are $2\delta$-close to $Y$.
		We denote by $p$ and $p$ their respective projections on $Y$.
		By construction the geodesics $\geo xp$ and $\geo {x'}{p'}$ lie in the $2\delta$-neighborhood of $Y$.
		By quasi-convexity the same holds for $\geo p{p'}$.
		Thus by concatenating the three geodesics we obtain a path contained in the $2\delta$-neighborhood of $Y$ joining $x$ to $x'$ whose length is at most $\dist x{x'} + 8\delta$.
		Consequently the $2\delta$-neighborhood of $Y$ is $8\delta$-strongly quasi-convex.
	\end{proof}

	\paragraph{}
	Let $x$ be a point of $X$.
	An isometry $g$ of $X$ is either
	\begin{itemize}
		\item \emph{elliptic} i.e., the orbit of $x$ under $g$ is bounded,
		\item \emph{parabolic} i.e., the orbit of $x$ under $g$ has exactly one accumulation point in $\partial X$.
		\item \emph{hyperbolic} i.e., the orbit of $x$ under $g$ has exactly two accumulation points in $\partial X$.
	\end{itemize}
	Note that these definitions do not depend on $x$.
	In order to measure the action of $g$ on $X$, we define two translation lengths.
	By the \emph{translation length} $\len[espace=X]g$ (or simply $\len g$) we mean 
	\begin{displaymath}
		\len[espace=X] g = \inf_{x \in X} \dist {gx}x.
	\end{displaymath}
	The \emph{asymptotic translation length} $\len[stable, espace=X] g$ (or simply $\len[stable]g$) is
	\begin{displaymath}
		\len[espace=X,stable] g = \lim_{n \rightarrow + \infty} \frac 1n \dist{g^nx}x.
	\end{displaymath}
	These two lengths satisfy the following inequality $\len[stable]g \leq \len g \leq \len[stable] g + 32\delta$ \cite[Chap. 10, Prop. 6.4]{CooDelPap90}.
	An isometry $g$ of $X$ is hyperbolic if and only if $\len[stable] g >0$ \cite[Chap. 10, Prop. 6.3]{CooDelPap90}.

	\begin{lemm}
	\label{res: quasi-convexity distance isometry}
		Let $x$, $x'$ and $y$ be three points of $X$.
		Let $g$ be an isometry of $X$.
		Then $\dist {gy}y \leq \max\left\{ \dist {gx}x, \dist {gx'}{x'} \right\} + 2 \gro x{x'}y  + 6 \delta$.
	\end{lemm}

	\begin{proof}
		By hyperbolicity
		\begin{equation}
		\label{eqn: quasi-convexity distance isometry - 1}
			\min \left\{ \fantomB \gro x{gx}y , \gro {gx}{gx'}y , \gro {gx'}{x'}y \right\} \leq \gro x{x'}y + 2 \delta.
		\end{equation}
		Assume that the minimum is achieved by $\gro x{gx}y$.
		Using the triangle inequality we obtain $\dist {gy}y \leq \dist {gx}x + 2  \gro x{gx}y \leq  \dist {gx}x + 2 \gro x{x'}y + 4 \delta$.
		A similar inequality holds if the minimum is achieved by $\gro {gx'}{x'}y$.
		Suppose now that the minimum in (\ref{eqn: quasi-convexity distance isometry - 1}) is achieved by $\gro {gx}{gx'}y$.
		Hence $\gro {gx}{gx'}y \leq \gro x{x'}y + 2 \delta$.
		Applying (\ref{eqn: hyperbolicity condition 2}) we obtain
		\begin{equation}
		\label{eqn: quasi-convexity distance isometry - 2}
			\dist {gy}y + \dist {gx}{gx'}\leq \max \left\{ \fantomB \dist {gx}{gy} + \dist {gx'}y,  \dist {gx'}{gy} + \dist {gx}y \right\} +2 \delta.
		\end{equation}
		However, by triangle inequality 
		\begin{eqnarray*}
			\dist {gx}{gy} + \dist {gx'}y 
			& \leq &  \dist{gx}x + \dist {gx'}{gx} + 2 \gro {gx}{gx'}y\\
			& \leq & \dist{gx}x + \dist {gx'}{gx}+ 2\gro x{x'}y +4 \delta.
		\end{eqnarray*}
		The same inequality holds after swapping $x$ and $x'$.
		Therefore (\ref{eqn: quasi-convexity distance isometry - 2}) leads to the desired result.
	\end{proof}

	\begin{defi}
	\label{def: axes}
		Let $g$ be an isometry of $X$.
		The \emph{axis} of $g$ denoted by $A_g$ is the set of points $x \in X$ such that  $\dist {gx}x  \leq \max\{ \len g, 8\delta\}$.
	\end{defi}
	
	\rems Note that we do not require $g$ to be hyperbolic.
	This definition works also for parabolic or elliptic isometries.
	This subset is not empty because $X$ is proper.
	It is also closed.
		
	\begin{prop}
	\label{res: axes is quasi-convex}
	Let $g$ be an isometry of $X$.
	Let $x$ be a point of $X$.
	\begin{enumerate}
		\item \label{enu: axes - min deplacement}
		$\dist {gx}x \geq 2 d(x, A_g) + \len g - 14\delta$,
		\item \label{enu: axes - deplacement donne distance a l'axe}
		if $\dist{gx}x \leq \len g + A$, then $d(x,A_g) \leq \frac 12 A+7\delta$,
		\item \label{enu: axes - quasi-convex}
		$A_g$ is $14\delta$-quasi-convex. 
	\end{enumerate}
	\end{prop}
	
	\begin{proof}
		Let $x \in X$. 
		Note that if $x$ belongs to $A_g$, Point \ref{enu: axes - min deplacement} is true.
		Therefore we can assume that $x \notin A_g$.
		We denote by $y$ a projection of $x$ on $A_g$. 
		Observe that any geodesic $\geo y{gy}$ is contained in $A_g$.
		Such a geodesic is $3\delta$-quasi-convex.
		Moreover, $y$ and $gy$ are respective projections of $x$ and $gx$ on it.
		Proposition~\ref{res: proj quasi-convex} gives	
		\begin{equation}
		\label{eqn: axes - maj projections}
			\dist {gy}y \leq \max\left\{ \fantomB \dist {gx}x -2 \dist xy + 14\delta, 7\delta \right\}.
		\end{equation}
		On the other hand we claim that $\dist {gy}y \geq 8\delta$.
		Recall that $x$ does not belong to $A_g$.
		By construction of $y$ any point $z$ on $\geo xy$ distinct from $y$ does not belong to $A_g$,
		Therefore $\dist {gz}z \geq \max\{\len g, 8 \delta\}$.
		Taking the limit as $z$ approaches $y$ leads to the claim.
		Hence (\ref{eqn: axes - maj projections}) gives
		\begin{equation}
		\label{eqn : axes conclusion projection}
			\dist {gx}x \geq \dist {gy}y + 2 \dist xy - 14 \delta \geq \len g + 2d(x,A_g) - 14 \delta,
		\end{equation}
		which proves Point~\ref{enu: axes - min deplacement}.
		Point~\ref{enu: axes - deplacement donne distance a l'axe} is a consequence of \ref{enu: axes - min deplacement}.
		Let us now prove Point~\ref{enu: axes - quasi-convex}.
		Let $y$ and $y'$ be two points of $A_g$.
		Let $x$ be a point of $X$.
		By Lemma~\ref{res: quasi-convexity distance isometry},
		\begin{displaymath}
			\dist {gx}x \leq \max \left\{ \fantomB \dist {gy}y, \dist{gy'}{y'}\right\} + 2\gro y{y'}x + 6 \delta \leq \len g + 2 \gro y{y'}x + 14\delta.
		\end{displaymath}
		It follows then from Point~\ref{enu: axes - deplacement donne distance a l'axe} that $d(x,A_g) \leq \gro y{y'}x + 14\delta$.
	\end{proof}	
	
	\paragraph{}
	Let $g$ be a hyperbolic isometry of $X$.
	We write $g^-$ and $g^+$ for the accumulation points in $\partial X$ of an orbit of $g$.
	They are the only points of $\partial X$ fixed by $g$.
	
	\begin{defi}
	\label{def: cylinder}
		Let $g$ be a hyperbolic isometry of $X$.
		The \emph{cylinder} of $g$ denoted by $Y_g$ is the set of points lying in the $10\delta$-neighborhood of a geodesic joining $g^-$ and $g^+$.
	\end{defi}
	
	\begin{lemm}
	\label{res: cylinder quasi-convex}
		Let $g$ be a hyperbolic isometry of $X$.
		The union $\Gamma$ of all geodesics joining $g^-$ to $g^+$ is $8\delta$-quasi-convex.
		The set $Y_g$ is strongly quasi-convex.
	\end{lemm}

	\begin{proof}
	According to Lemma~\ref{res: strong neighborhood of a quasi-convex} it is sufficient to prove that $\Gamma$ is $8\delta$-quasi-convex.
	Let $x \in X$ and $y,y' \in \Gamma$.
	There exist $\gamma$ and $\gamma'$ two geodesics joining $g^-$ and $g^+$ such that $y$ and $y'$ respectively lie on $\gamma$ and $\gamma'$.
	We denote by $p'$ a projection of $y'$ on $\gamma$.
	Since $\gamma$ and $\gamma'$ join the same extremities the Hausdorff distance between them is at most $5\delta$.
	Thus $\dist {y'}{p'} \leq 5\delta$.
	The path $\gamma$ being a bi-infinite geodesic, it follows that 
	\begin{displaymath}
		d(x, \Gamma) \leq d(x, \gamma) \leq \gro y{p'}x + 3 \delta \leq \gro y{y'}x + 8\delta. \qedhere
	\end{displaymath}
\end{proof}
	
	\begin{lemm}
	\label{res: cylinder contained in invariant subset}
		Let $g$ be a hyperbolic isometry of $X$.
		Let $Y$ be a $g$-invariant $\alpha$-quasi-convex subset of $X$.
		Then $Y_g$ lies in the $(\alpha + 22\delta)$-neighborhood of $Y$.
		In particular $Y_g$ is contained in the $36\delta$-neighborhood of $A_g$.
	\end{lemm}
	
	\begin{proof}
		It is sufficient to prove that every bi-infinite geodesic joining $g^-$ to $g^+$ lies in the $(\alpha +12\delta)$-neighborhood of $Y$.
		Let $\gamma$ be such a geodesic and $x$ a point of $\gamma$.
		Let $\eta >0$.
		We denote by $y$ an $\eta$-projection of $x$ on $Y$.
		Since $g$ is hyperbolic there exists $m \in \N$ such that $\dist {g^mx}{g^{-m}x} > 2\dist xy + 26\delta$.
		The geodesics $g^m\gamma$ and $g^{-m}\gamma$ also join $g^-$ and $g^+$.
		Therefore $g^m x$ and $g^{-m}x$ are $5\delta$-close to $\gamma$.
		We denote by $p_-$ and $p_+$ respective projections of these points on $\gamma$.
		Note that $x$ lies on the portion of $\gamma$ between $p_-$ and $p_+$.
		Indeed if it was not the case we would have
		\begin{displaymath}
			\dist {g^mx}{g^{-m}x} 
			\leq \dist{p_-}{p_+} + 10\delta 
			= \dist{\fantomB\dist x{p_-}}{\dist x{p_+}}+10\delta
			\leq 20\delta,
		\end{displaymath}
		which contradicts our assumption on $m$.
		In particular $\gro {g^mx}{g^{-m}x}x \leq 10\delta$.
		By hyperbolicity we get
		\begin{displaymath}
			\min\left\{\dist{g^mx}x - \dist xy ,\gro{g^my}{g^{-m}y}x\right\} \leq\gro{g^mx}{g^{-m}x}x + 2\delta \leq 12\delta
		\end{displaymath}
		By construction of $m$, $\dist{g^mx}x$ is bounded below by $\dist xy +13\delta$.
		Therefore the minimum in the previous inequality is achieved by $\gro {g^my}{g^{-m}y}x$.
		However $Y$ being $g$-invariant and $\alpha$-quasi-convex, $g^my$ and $g^{-m}y$ are two points of $Y$ and $d(x,Y) \leq \gro {g^my}{g^{-m}y}x +\alpha$.
		Consequently $x$ lies in the $(\alpha + 12\delta)$-neighborhood of $Y$. 
	\end{proof}

	\paragraph{}Let $g$ be an isometry of $X$ such that $\len g >200\delta$.
	(In particular, $g$ is hyperbolic.)
	Let $x$ be a point of $A_g$.
	We consider a geodesic $\gamma : J \rightarrow X$ between $x$ and $gx$ parametrized by arc length. 
	We extend $\gamma$ in a $g$-invariant path $\gamma : \R \rightarrow X$ in the following way: for all $t \in J$, for all $m \in \Z$, $\gamma\left(t + m \len g \right) = g^m\gamma(t)$. 
	This is a $\len g$-local (1,0)-quasi-geodesic contained in $A_g$ joining $g^-$ to $g^+$.
	By stability of quasi-geodesics $\gamma$ is actually $8\delta$-quasi-convex.
	We call such a path a \emph{nerve} of $g$.
	The Hausdorff distance between two nerves of $g$ is at most $5\delta$.	

	\begin{lemm}
	\label{res: axis in cylinder}
		Let $g \in G$ such that $\len g > 200\delta$.
		Let $\gamma$ be a $200\delta$-local $(1,0)$-quasi-geodesic joining $g^-$ to $g^+$
		Then $A_g$ is contained in the $5\delta$-neighborhood of $\gamma$.
		In particular $A_g$ is contained  in $Y_g$ and in the $5\delta$-neighborhood of any nerve of $g$.
	\end{lemm}
	
	\begin{proof}
		Let $x \in A_g$.
		There exists $\gamma'$ a nerve of $g$ going through $x$.
		Both $\gamma$ and $\gamma'$ are $200\delta$-local $(1,0)$-quasi-geodesic joining $g^-$ to $g^+$.
		By the stability of quasi-geodesics $\gamma'$ and thus $x$ lies in the $5\delta$-neighborhood of $\gamma$.
	\end{proof}

	The next lemma explains the following fact.
	Let $g$ be a hyperbolic isometry of $X$.
	A quasi-geodesic contained in the neighborhood of the axis of $g$ almost behaves like a nerve of $g$.

	\begin{lemm}
	\label{res : quasi-geodesic behaving like a nerve}
		Let $g \in G$ such that $\len g > 200 \delta$.
		Let $\gamma : \intval ab \rightarrow X$ be a $\len g$-local $(1,0)$-quasi-geodesic contained in the $A$-neighborhood of $A_g$.
		Then there exists $\epsilon \in \{ \pm 1 \}$ such that for every $s \in \intval ab$  if $s \leq b - \len g$ then 
		\begin{displaymath}
			\dist {g^\epsilon\gamma(s)}{\gamma(s+ \len g)} \leq 4A + 80\delta.
		\end{displaymath}
	\end{lemm}
	
	\begin{proof}
		We denote by $\gamma_g$ an nerve of $g$.
		Since $\len g > 200\delta$, the $5\delta$-neighborhood of $\gamma_g$ contains $A_g$.
		Thus $\gamma$ lies in the $(A + 5\delta)$-neighborhood of $\gamma_g$.
		In particular there exist $c, d \in \R$ such that $\dist {\gamma(a)}{\gamma_g(c)} \leq A + 5\delta$ and $\dist {\gamma(b)}{\gamma_g(d)}\leq A +5\delta$.
		By replacing if necessary $g$ by $g^{-1}$ we can assume that $c \leq d$.

		\paragraph{}
		Let $s \in \intval ab$ such that $s \leq t- \len g$.
		Using the stability of quasi-geodesics there exist $t \in \intval cd$ and $t' \in \intval td$ such that $\dist {\gamma(s)}{\gamma_g(t)} \leq A + 15\delta$ and $\dist {\gamma(s+ \len g)}{\gamma_g(t')} \leq A + 19\delta$.
		It follows that 
		\begin{displaymath}
			\dist{\fantomA \dist {\gamma_g(t)}{\gamma_g(t')}}{\dist {\gamma(s)}{\gamma(s+ \len g)}} \leq  2A + 34\delta.
		\end{displaymath}
		However $\dist {\gamma(s)}{\gamma(s+ \len g)}$ and $\dist {\gamma_g(t)}{\gamma_g(t+ \len g)}$ both equal $\len g$.
		Consequently
		\begin{displaymath}
			\dist{\fantomA \dist {\gamma_g(t)}{\gamma_g(t')}}{\dist {\gamma_g(t)}{\gamma_g(t+ \len g)}} \leq  2A + 34\delta.
		\end{displaymath}
		Since $t'$ and $t + \len g$ are larger than $t$ we get by Lemma~\ref{res: metric inequalities}-\ref{enu: metric inequalities - two points close to a geodesic}.
		\begin{displaymath}
			\dist{\gamma_g(t')}{g\gamma_g(t)} = \dist{\gamma_g(t')}{\gamma_g(t+\len g)} \leq 2A + 46\delta.
		\end{displaymath}
		It follows then from the triangle inequality that $\dist {g\gamma(s)}{\gamma(s + \len g)} \leq 4A  + 80\delta$.
	\end{proof}

\subsection{Group acting on a hyperbolic space}
\label{sec: group acting on a hyperbolic space}

	In this section $G$ denotes a group acting by isometries on $X$.
	We still assume that $X$ is geodesic and proper.
	Moreover, we require the action of $G$ on $X$ to be
	\begin{enumerate}
		\item 
		\emph{proper} i.e., for every $x \in X$, there exists $r >0$ such that the set of elements $g \in G$ satisfying $gB(x,r) \cap B(x,r) \neq \emptyset$ is finite.
		\item \emph{co-compact} i.e., the quotient $X/G$ endowed with the induced topology is compact.
	\end{enumerate}
	Since the space $X$ is proper, the properness of the action of $G$ implies this more general fact.
	Let $Y$ be a bounded subset of $X$. The set of elements $g \in G$ such that $gY$ intersects $Y$ is finite \cite[Chap. I.8, Remark 8.3]{BriHae99}.
	It follows from these assumptions that a subgroup of $G$ is either \emph{elementary} i.e., virtually cyclic or contains a free group of rank 2 \cite[Chap. 8, Th. 37]{GhyHar90}.
	
	\nota Given a subset $Y$ of $X$ we denote by $\stab Y$ the stabilizer of $Y$ i.e., the set of elements $g \in G$ such that $gY = Y$.
	
	\paragraph{Finite subgroups.}
	We start by studying some properties of the finite subgroups of $G$.
	To that end we associate to each such subgroup a particular subset of $X$.
	\begin{defi}
	\label{def: characteristic set of a finite subgroup}
		Given a finite subgroup $F$ of $G$ we denote by $C_F$ the set of points $x \in X$ such that for every $g \in F$, $\dist {gx}x \leq 10 \delta$.
	\end{defi}
	
	It follows from the definition that $C_F$ is an $F$-invariant subset of $X$.
		
	\begin{prop}
	\label{res: characteristic subset elliptic group - technical lemma}
		Let $F$ be a finite subgroup of $G$.
		Let $x$ be a point of $X$.
		Let $g \in F$ such that $\dist {gx}x$ is maximal.
		We denote by $m$ the midpoint of a geodesic $\geo x{gx}$.
		Then $m$ belongs to $C_F$.
		In particular $C_F$ is non-empty.
	\end{prop}
	
	\begin{proof}
	Let $h$ be an element of $F$.
	In order to simplify the notations we put $y = gx$, $z = hx$ and $t = hgx$.
	The point $p = hm$ is thus the midpoint of $\geo zt$.
	The only fact that we are going to use is that $\dist xy = \dist zt$ is the largest distance between any two points of $\{x,y,z,t\}$. 
	Using hyperbolicity condition~(\ref{eqn: hyperbolicity condition 2}) we have
	\begin{equation}
	\label{eqn: characteristic subset elliptic group - technical lemma}
		\dist xy + \dist zt \leq \max\left\{\dist xt +\dist yz, \dist yt + \dist xz \right\} +2\delta.
	\end{equation}
	Note that $z$ and $t$ play a symmetric role.
	Without loss of generality we can assume that the maximum is achieved by $\dist xt +\dist yz$.
	It follows that 
	\begin{displaymath}
		\left(\fantomB \dist xy- \dist xt \right) + \left(\fantomB \dist zt - \dist yz \right) \leq 2 \delta.
	\end{displaymath}
	Recall that $\dist xy$ and $\dist zt$ are respectively larger than or equal to $\dist xt$ and $ \dist yz$.
	Consequently $0 \leq \dist xy - \dist xt \leq 2 \delta$.
	Roughly speaking the triangle $[x,y,t]$ has two sides with approximatively the same length namely $\geo xy$ and $\geo xt$, this length being larger than the one of the last side.
	It follows that $\gro ytx \geq \dist xm - \delta$. 
	Similarly we have $\gro xzt \geq \dist tp - \delta$.
	Applying twice Lemma~\ref{res: metric inequalities}-\ref{enu: metric inequalities - thin triangle}, we obtain
	\begin{equationarray*}{lclcl}
		\gro xtm
		& \leq &  \max\left\{ \fantomB \dist xm - \gro ytx, \gro xym \right\} +\delta & \leq & 2\delta, \\
		\gro xtp
		& \leq  &   \max\left\{ \fantomB \dist tp - \gro xzt ,\gro ztp\right\} +\delta & \leq & 2\delta.
	\end{equationarray*}
	Lemma~\ref{res: metric inequalities}-\ref{enu: metric inequalities - two points close to a geodesic}, leads then to 
	\begin{displaymath}
		\dist pm
		\leq \left|\fantomB \dist xp- \dist xm \right| + 2\max\left\{\gro xtm, \gro xtp \right\} + 2\delta
		\leq \left|\fantomB \dist xp - \dist xm \right|+ 6\delta . 
	\end{displaymath}
	However $m$ and $p$ are respectively the midpoints of $\geo xy$ and $\geo zt$ which have the same length, thus
	\begin{displaymath}
		\dist xp - \dist xm = \left(\fantomB \dist xt - \dist xy\right) + 2 \gro xtp. 
	\end{displaymath}
	Therefore $\dist pm\leq 10\delta$ i.e., $\dist {hm}m \leq 10\delta$.
	In other words $m$ belongs to $C_F$.
	\end{proof}
	
	\begin{coro}
	\label{res: characteristic subset elliptic group - quasi convex}
		Let $F$ be a finite subgroup of $G$.
		The subset $C_F$ is $8\delta$-quasi-convex.
	\end{coro}
	
	\begin{proof}
		Let $x \in X$ and $y,y' \in C_F$.
		We denote by $g$ an element of $F$ such that  $\dist {gx}x$ is maximal and $m$ the midpoint of $\geo x{gx}$.
		According to Lemma~\ref{res: quasi-convexity distance isometry}, 
		\begin{displaymath}
			2 \dist xm = \dist {gx}x \leq \max\left\{ \dist{gy}y, \dist{gy'}{y'} \right\} + 2 \gro y{y'}x + 6 \delta.
		\end{displaymath}
		However $y$ and $y'$ belong to $C_F$, hence $ \dist xm \leq  \gro y{y'}x + 8\delta$.
		By Proposition~\ref{res: characteristic subset elliptic group - technical lemma}, $m$ belongs to $C_F$.
		Therefore $d(x,C_F) \leq \gro y{y'}x + 8\delta$
	\end{proof}

	\begin{coro}
	\label{res: characteristic subset elliptic group - non empty}
		Let $F$ be a finite subgroup of $G$.
		Let $Y$ be a non-empty, $F$-invariant, $\alpha$-quasi-convex subset of $X$.
		Then the $\alpha$-neighborhood of $Y$ intersects $C_F$.
	\end{coro}
	
	\begin{proof}
		Let $x$ be a point of $Y$.
		We denote by $g$ an element of $F$ such that  $\dist {gx}x$ is maximal and $m$ the midpoint of $\geo x{gx}$.
		According to Proposition~\ref{res: characteristic subset elliptic group - technical lemma}, $m$ lies in $C_F$.
		On the other hand, $Y$ is $F$-invariant, therefore $gx \in Y$.
		Since $Y$ is $\alpha$-quasi-convex, $d(y, Y) \leq \gro x{gx}m + \alpha \leq \alpha$.
		Consequently $m$ belongs to the $\alpha$-neighborhood of $Y$.
	\end{proof}

	\paragraph{Infinite elementary subgroups.}
		
	\paragraph{}
	Let $H$ be an infinite elementary subgroup of $G$.
	By definition $H$ contains a finite index subgroup isomorphic to $\Z$.
	The set of accumulation points in $\partial X$ of an orbit of $H$, that we denote $\partial H$, has exactly two points.
	There exists a subgroup $H^+$ of $H$ of index at most 2 which fixes pointwise $\partial H$.
	If $H^+ \neq H$ then  $H$ contains an element of order 2.
	A Schur Theorem (see \cite[Th. 5.32]{Rotman:1995ud}) implies that $H^+$ contains a unique maximal finite subgroup $F$. 
	This  group is actually a normal subgroup of $H^+$.
	Moreover, there exists a hyperbolic element $h \in H^+$ so that $H^+$ is isomorphic to $\sdp F{\Z}$ where $\Z$ is identified with the subgroup $\left< h \right>$ acting by conjugation on $F$.
	
	\paragraph{}Let $g$ be a hyperbolic element of $G$.
 	The subgroup $E$ of $G$ that stabilizes $\{g^-,g^+\}$ is the maximal elementary subgroup of $G$ containing $g$ \cite[Chap. 10, Prop. 7.1]{CooDelPap90}.
	The isometry $g$ is said to be a \emph{proper power} if there exist $h \in G$ and an integer $n\geq 2$ such that $g=h^n$.
	Any hyperbolic element of $G$ is a power of an isometry which is not a proper power.

	\begin{lemm}
	\label{res: cylinder almost invariant by max finite subgroup}
		Let $g$ be a hyperbolic element of $G$ and $H$ a subgroup of $G$ fixing pointwise $\{g^-,g^+\}$.
		Let $F$ be the maximal finite subgroup of $H$.
		The cylinder $Y_g$ of $g$ is contained in the $48\delta$-neighborhood of $C_F$.
	\end{lemm}
	
	\begin{proof}
		According to Lemma~\ref{res: cylinder quasi-convex}, the union $\Gamma$ of all geodesics joining $g^-$ to $g^+$ is $8\delta$-quasi-convex.
		It is also $F$-invariant.
		By Lemma~\ref{res: characteristic subset elliptic group - non empty}, there exists a point $x$ in $C_F\cap Y_g$.
		The subgroup $F$ being $g$ invariant, for every $n \in \Z$, $g^nx$ belongs to $C_F\cap Y_g$.
		Note that $g^nx$ tends to $g^+$ (\resp $g^-$) and $n$ approaches $+\infty$ (\resp $-\infty$).
		In particular, for every $y \in Y_g$ there exist $n,m \in \Z$ such that $\gro {g^nx}{g^mx}y \leq 40\delta$.
		Since $C_F$ is $8\delta$-quasi-convex, $y$ lies in the $48\delta$-neighborhood of $C_F$.
	\end{proof}

	\begin{lemm}
	\label{res: pseudo nu = 1}
		Assume that every elementary subgroup of $G$ is cyclic.
		Let $g, h \in G$.
		If $g$ and $hgh^{-1}$ generate an elementary subgroup then so do $g$ and $h$.
	\end{lemm}
	
	\begin{proof}
		We denote by $H$ the subgroup of $G$ generated by $g$ and $hgh^{-1}$.
		We distinguish two cases. 
		If $g$ is hyperbolic then $g$ and $hgh^{-1}$ are two hyperbolic isometries with the same accumulation points in $\partial X$.
		In particular $h$ stabilizes the set $\{g^-, g^+\}$.
		It follows that $g$ and $h$ belong to the elementary subgroup $E(g)$.
		
		\paragraph{} Assume now that $g$ has finite order.
		The subgroup $H$ is elementary, thus cyclic.
		In particular it has to be finite.
		The isometries $g$ and $hgh^{-1}$ generate two subgroups of $H$ with the same order.
		However two such subgroups in a cyclic group are equal.
		Therefore there exists $m \in \Z$ such that $hgh^{-1} = g^m$.
		Thus any element of the subgroup generated by $g$ and $h$ can be written $g^ph^q$ with $p,q \in \Z$.
		Since $g$ has finite order, this subgroup is also elementary.
	\end{proof}

	\begin{lemm}
	\label{res: extracting subset without inverse}
		We assume that every elementary subgroup of $G$ is cyclic.
		Let $n \in \N$.
		Let $g$ and $h$ be two hyperbolic elements of $G$ which are not proper powers.
		Either $g$ and $h$ generate a non-elementary subgroup of $G$ or $\left<g^n\right> = \left< h^n\right>$.		
	\end{lemm}
	
	\begin{proof}
		Assume that $g$ and $h$ generate an elementary subgroup.
		This subgroup is infinite and cyclic.
		Since $g$ and $h$ are not proper powers, they are either equal or inverse.
		Hence $\left<g^n\right> = \left< h^n\right>$.
	\end{proof}

	\paragraph{Group invariants.}
	We now introduce several invariants associated to the action of $G$ on $X$.
	During the final induction, they will be useful to ensure that the set of relations we are looking at satisfy a small cancellation assumption.
	
	\begin{defi}
	\label{def: injectivity radius}
		Let $P$ be a subset of $G$.
		The \emph{injectivity radius} of $P$ on $X$, denoted by $\rinj PX$ is 
		\begin{displaymath}
			\rinj PX = \inf \set{\len[stable] g}{g \in P, \text{ hyperbolic}}
		\end{displaymath}
	\end{defi}
	
	\begin{prop}[see {\cite[Prop. 3.1]{Del96}}]
	\label{res: positive injectivity radius}
		There exists $a>0$ such that for every hyperbolic element $g\in G$ we have $\len[stable] g \in a\N$.
		In particular $\rinj GX >0$.
	\end{prop}

	\begin{defi}
	\label{def: invariant A}
		We denote by $\mathcal A$ the set of pairs $(g,h)$ generating a non-elementary subgroup of $G$ such that  $\len g \leq 1000\delta$ and $\len h \leq 1000\delta$. 
		The parameter $A(G,X)$ is given by
		\begin{displaymath}
			A(G,X) = \sup_{(g,h) \in \mathcal A} \diam \left( A_g^{+17\delta} \cap A_h^{+17\delta}\right)
		\end{displaymath}
	\end{defi}
	
	\label{rem: dependency delta}
	The invariant $A(G,X)$ depends implicitly on the hyperbolicity constant $\delta$.
	Although the notation does not make this dependency explicit, we should keep in mind that it plays an important role. 
	For instance, we have the following lemma:
	\begin{lemm}
	\label{theo: rescaling A}
		Let $\lambda$ be a positive number.
		We denote by $\lambda X$ the space $X$ endowed with the rescaled metric $\lambda \distV[X]$ and view it as a $\lambda\delta$-hyperbolic space.
		Then $A(G, \lambda X) = \lambda A(G,X)$.
	\end{lemm}
	
	\begin{proof}
		Let $g$ be an element of $G$.
		Its translation length satisfies $\len[espace=\lambda X] g = \lambda \len[espace=X] g$.
		Since $\lambda X$ is a $\lambda\delta$-hyperbolic space, the axis of $g$ in $\lambda X$ is exactly the image in $\lambda X$ of the axis $A_g$ of $g$ in $X$. 
		We will denote it by $\lambda A_g$.
		Let $g$ and $g'$ be two elements of $G$ that do not generate an elementary subgroup and whose translation lengths in $\lambda X$ are at most than $1000\lambda\delta$.
		In particular, we have $\len[espace=X]g, \len[espace=X]{g'} \leq 1000 \delta$.
		By definition of $A(G,X)$, we get
		\begin{eqnarray*}
			\diam \left( \lambda A_g^{+17 \lambda\delta} \cap \lambda A_{g'}^{+17 \lambda\delta}  \right) 
			& = & \lambda \diam \left( A_g^{+17 \delta} \cap A_{g'}^{+17\delta}  \right) \\
			& \leq & \lambda A(G,X).
		\end{eqnarray*}
		After taking the upper bound for all $g$ and $g'$, we obtain $A(G, \lambda X) \leq \lambda A(G,X)$.
		In the same way, $A(G, \lambda X) \geq \lambda A(G,X)$.
		This establishes the desired equality.
	\end{proof}
	
	\begin{prop}
	\label{res: overlap two axes}
		We assume that every elementary subgroup of $G$ is cyclic.
		Let $g$ and $h$ be two elements of $G$ which generate a non-elementary subgroup.
		\begin{enumerate}
			\item \label{enu: overlap axes short}
			If $\len g \leq 1000 \delta$, then $\diam \left( A_g^{+17\delta} \cap A_h^{+ 17\delta}\right)\leq \len h + A(G,X) + 158\delta$.
			\item \label{enu: overlap axes general}
			Without assumption on $g$ we have, 
			\begin{displaymath}
				\diam \left( A_g^{+17\delta} \cap A_h^{+ 17\delta}\right) \leq \len g + \len h + \max\{\len g, \len h\} + A(G,X) + 676\delta.
			\end{displaymath}

		\end{enumerate}
	\end{prop}
	
	\begin{proof}
		 We prove Point~\ref{enu: overlap axes short} by contradiction. 
		 Assume that
		 \begin{displaymath}
		 	\diam \left( A_g^{+17\delta} \cap A_h^{+17\delta}\right) > \len h + A(G,X) + 158\delta.
		 \end{displaymath}
		By definition of $A(G,X)$ we have $\len h > 1000\delta$, otherwise $g$ and $h$ would generate an elementary subgroup.
		We denote by $\gamma : \R \rightarrow X$ a nerve of $h$.
		Its $5\delta$-neighborhood contains $A_g$ (Lemma~\ref{res: axis in cylinder}) therefore by Proposition~\ref{res: intersection of thickened quasi-convex}
		\begin{displaymath}
			\diam\left(A_g^{+17\delta}\cap \gamma^{+11\delta}\right) > \len h + A(G,X) + 110\delta.
		\end{displaymath}
		In particular there exist two points on $\gamma$,  $x=\gamma(s)$ and $x'=\gamma(s')$  which also belong to the $28\delta$-neighborhood of $A_g$ and such that 
		\begin{equation}
		\label{eqn: overlap axes short}
			\dist x{x'} >  \len h + A(G,X) + 88\delta.
		\end{equation}
		By replacing if necessary $h$ by $h^{-1}$ we can assume that $s \leq s'$.
		By stability of quasi-geodesics for all $t \in \intval s{s'}$, $\gro x{x'}{\gamma(t)} \leq 5\delta$.
		Since the $28\delta$-neighborhood of $A_g$ is $2\delta$-quasi-convex (see Lemma~\ref{res: neighborhood of a quasi-convex}), $\gamma(t)$ lies in the $35\delta$-neighborhood of $A_g$.
		Thus $\dist{g\gamma(t)}{\gamma(t)} \leq \len g + 70\delta$.
		According to (\ref{eqn: overlap axes short}) there exists $t \in \intval s{s'}$ such that $\dist {x'}{\gamma(t)} =  \len h$.
		We put $y= \gamma(t)$.
		By construction $hx = \gamma(s+ \len h)$ and $hy = \gamma(t + \len h)$.
		Note that $\dist {s'}t \geq \len h$, thus $s + \len h$ and $t + \len h$ belong to $\intval s{s'}$.
		Hence 
		\begin{displaymath}
			\dist{ghx}{hx}, \dist{ghy}{hy} \leq \len {hgh^{-1}} + 70\delta.
		\end{displaymath}
		It follows from Proposition~\ref{res: axes is quasi-convex}, that $x$ and $y$ belong to the $42\delta$-neighborhood of $hA_g$.
		Consequently $x$ and $y$ are two points of $A_g^{+35\delta} \cap hA_g^{+42\delta}$.
		By Proposition~\ref{res: intersection of thickened quasi-convex},
		\begin{displaymath}
			\diam \left( A_g^{+17\delta} \cap A_{hgh^{-1}}^{+17\delta}\right)\geq \dist xy - 88\delta \geq  \dist{x'}x - \dist{x'}y - 88\delta > A(G,X).
		\end{displaymath}
		Moreover, $\len {hgh^{-1}} = \len g \leq 1000\delta$.
		By definition of $A(G,X)$ the isometries $g$ and $hgh^{-1}$ generate an elementary group.
		It follows from Lemma~\ref{res: pseudo nu = 1} that $g$ and $h$ also generate an elementary group.
		Contradiction.
		
		\paragraph{} 
		We now prove Point~\ref{enu: overlap axes general}.
		According to the previous point we can assume that $\len g > 1000\delta$ and $\len h > 1000 \delta$.
		Without loss of generality we can suppose $\len h \geq \len g$.
		Imagine now that
		\begin{displaymath}
			\diam \left( A_g^{+17\delta} \cap A_h^{+ 17\delta}\right) > \len g + 2\len h + A(G,X) + 676\delta.
		\end{displaymath}
		We denote by $\gamma$ a  nerve of $h$.
		Its $5\delta$-neighborhood contains $A_h$ thus
		\begin{displaymath}
			\diam \left(\gamma^{+ 11\delta} \cap A_g^{+ 17\delta}\right) > \len g + 2\len h + A(G,X) + 628\delta.
		\end{displaymath}
		In particular there exit $x = \gamma(s)$, $x'= \gamma(s')$ lying in the $28\delta$-neighborhood of $A_g$ such that
		\begin{displaymath}
			\dist x{x'}> \len g + 2\len h + A(G,X) + 616\delta.
		\end{displaymath}
		Without loss of generality we can assume that $s \leq s'$.
		As previously, the restriction of $\gamma$ to $\intval s{s'}$ is contained in the $35\delta$-neighborhood of $A_g$.
		We apply Lemma~\ref{res : quasi-geodesic behaving like a nerve}.
		By replacing if necessary $g$ by $g^{-1}$, for every $t \in \intval s{s'}$ if $t \leq s'-\len g$ then $\dist {g\gamma(t)}{\gamma(t+ \len g)} \leq 220\delta$.
		Consequently, for every $t \in \intval s{s'}$ such that $t \leq s' - \len g - \len h$ we have
		\begin{displaymath}
			\dist{gh\gamma(t)}{hg\gamma(t)}
			\leq \dist{g\gamma(t + \len h)}{h\gamma(t+\len g)} + 220\delta
			\leq 440\delta.
		\end{displaymath}	
		It follows that the translation length of the isometry $u = h^{-1}g^{-1}hg$ is at most $1000\delta$ and for all $t \in \intval s{s'}$ if  $t \leq s' - \len g - \len h$ then $\gamma(t)$ is in the $227\delta$-neighborhood of $A_u$ (Proposition~\ref{res: axes is quasi-convex}).
		Let $y = \gamma(t)$ be the point of $\gamma$ such that $\dist {x'}y = \len h + \len g$.
		In particular $x$ and $y$ belong to the $227\delta$-neighborhood of $A_u$ and $A_h$.
		Therefore
		\begin{displaymath}
			\diam \left( A_u^{+17\delta} \cap A_h^{+17\delta}\right)\geq \dist xy - 458\delta \geq \dist {x'}x - \dist{x'}y - 458\delta >  \len h + A(G,X) + 158\delta.
		\end{displaymath}
		It follows from the previous point that $h$ and $u$ generate an elementary subgroup.
		Hence so do $h$ and $g^{-1}hg$.
		However $h$ is a hyperbolic isometry. 
		Consequently $g$ and $h$ generate an elementary group. 
		Contradiction.
	\end{proof}

\section{Rotation families}
	\label{sec: rotation family}
	
	\paragraph{}In this section we follow the presentation of rotation family given by F.~Dahmani, V.~Guirardel and D.~Osin in \cite{Dahmani:2011vu}.
	Four our purpose we only need a weaker statement about rotation family. 
	On the other hand, we work with an eventually non-geodesic space.
	In order to have a self contained text, we recall the main ideas of the proofs.
	
	\paragraph{}Let $G$ be a group acting by isometries on a $\delta$-hyperbolic length space $X$.
	Note that in this this section we do not require the space $X$ to be geodesic or proper.
	Similarly there is no assumption on the action of $G$ on $X$.
	
	\begin{defi}
	\label{def: rotation family}
		Let $\sigma >0$.
		A \emph{$\sigma$-rotation family} is a collection $\mathcal R$ of pairs $(H,v)$ where $H$ is a subgroup of $G$  and $v$ a point of $X$ satisfying the following properties.
		\begin{labelledenu}[R]
			\item \label{enu: rotation family - large angle}
			For every $(H,v) \in \mathcal R$, $H$ is  subgroup of $\stab v$ such that for every $x \in B(v, \sigma/10)$, for every $h \in H\setminus\{1\}$, $\dist{hx}x = 2\dist vx$.
			\item \label{enu: rotation family - apices apart}
			For every $(H,v), (H',v') \in \mathcal R$, if $(H,v) \neq (H',v')$, then $\dist v{v'} \geq \sigma$.
			\item \label{enu: rotation family - G invariant}
			$\mathcal R$ is stable under the action of $G$ defined as follows. For all $g \in G$, for all $(H,v) \in \mathcal R$, $g.(H,v) = (gHg^{-1}, gv)$.
		\end{labelledenu}
	\end{defi}
	
	\rem It follows from~\ref{enu: rotation family - apices apart} and \ref{enu: rotation family - G invariant} that for every $(H,v) \in \mathcal R$, $H$ is actually a normal subgroup of $\stab v$.

	\notas Let $(H,v) \in \mathcal R$. The idea is that each element  $h \in H$ acts on $X$ like a rotation of center $v$ and very large angle - see Axiom~\ref{enu: rotation family - large angle}.
	Therefore $v$ is called an \emph{apex} and $H$ a \emph{rotation group}.
	We denote by $v(\mathcal R)$ the set of all apices.
	Similarly $H(\mathcal R)$ stands for the set of all rotation groups $H$.
	Given a subset $Y$ of $X$, we denote by $K_Y$ the subgroup of $G$ generated by all $H$'s, where $(H,v) \in \mathcal R$ and $v \in Y$.
	
	\subsection{Fundamental theorem}
	\label{sec: rotation family - fundamental theorem}
	
	\begin{theo}
	\label{res: rotation family - fundamental theorem}
		There exists a positive number $\sigma_0=\sigma_0(\delta)$ depending only on $\delta$ with the following property.
		Let $\mathcal R$ be a $\sigma$-rotation family with $\sigma \geq \sigma_0$ and $K$ the (normal) subgroup of $G$ generated by all rotation groups $H \in H(\mathcal R)$.
		Let $x \in X$ and $g \in K$.
		If $g$ does not belong to any of the rotation groups $H\in H(\mathcal R)$ then $\dist {gx}x \geq \sigma -166\delta$.
	\end{theo}

	\paragraph{}
	The rest of this section is dedicated to the proof of the theorem.
	We need first to define $\sigma_0$.
	According to Propositions~\ref{res: stability quasi-geodesic} and \ref{res: stability discrete quasi-geodesic} there exists $\sigma_0 > 0$ depending only on $\delta$ with the following properties.
	\begin{itemize}
		\item For every $l \in\intval 0\delta$ any $\sigma_0/25$-local $(1,l)$-quasi-geodesic is a (global) $(2,l)$-quasi-geodesic.
		\item Assume that $y_0, \dots, y_{m+1}$ is a sequence of points of $X$ such that for every $i \in \intvald 1{m-1}$, $\dist{y_{i+1}}{y_i} \geq \sigma_0$ and for every $i \in \intvald 1m$, $\gro{y_{i+1}}{y_{i-1}}{y_i} \leq 7\delta$ then for every $i \in \intvald 0{m+1}$, $\gro {y_0}{y_{m+1}}{y_i} \leq 12\delta$.
		Moreover, for all $x \in X$ there exists $i \in \intvald 0m$, such that $\gro{y_{i+1}}{y_i}x \leq \gro {y_0}{y_{m+1}}x + 22\delta$.

	\end{itemize}
	Without loss of generality we can require that $\sigma_0 \geq 10^{20}\delta$.
	From now on we assume that $\mathcal R$ is a $\sigma$-rotation family with $\sigma \geq \sigma_0$.
	
	\begin{lemm}
	\label{res: rotation family - small product at the apex}
		Let $(H,v) \in \mathcal R$.
		Let $h \in H\setminus\{1\}$.
		For every $x \in X$, $\gro x{hx}v \leq  2\delta$.
	\end{lemm}
	
	\begin{proof}
		Let $x \in X$.
		If $\dist xv < \sigma/10$ the lemma is just a consequence of \ref{enu: rotation family - large angle}.
		Hence we can assume that $\dist xv \geq \sigma/10$.
		We denote by $y$ a point of $X$ such that $\dist vy = \sigma/20$ and $\gro xvy \leq \delta$.
		By hyperbolicity,
		\begin{displaymath}
			\min \left\{ \gro {hy}{hx}v , \gro{hx}xv, \gro xyv \right\} \leq \gro {hy}yv + 2\delta = 2\delta.
		\end{displaymath}
		However $\gro xyv = \dist vy - \gro xvy >2\delta$.
		Therefore the minimum cannot be achieved by $\gro xyv$ or $\gro {hx}{hy}v$.
		Consequently $\gro {hx}xv \leq 2 \delta$.
	\end{proof}
	
	\begin{defi}
	\label{defi: windmill}
		A non-empty subset $W$ is \emph{windmill} if it satisfies the following conditions.
		\begin{labelledenu}[W]
			\item \label{enu: windmill - quasi convex}
			$W$ is $2\delta$-quasi-convex,
			\item \label{enu: windmill - K_W invariant}
			$W$ is stable under the action of $K_W$,
			\item \label{enu: windmill - no apex in the neighborhood}
			for every $v \in v(\mathcal R)$, if $d(v,W) \leq  \sigma/10$, then $v$ belongs to $W$,
			\item \label{enu: windmill - almost free action}
			for every $g \in K_W$, for every $x \in X$, if $g$ does not belong to any rotation group $H \in H(\mathcal R)$ then $\dist{gx}x \geq  \sigma -166\delta$.
		\end{labelledenu}
	\end{defi}
		
	\begin{prop}
	\label{res: windmill - induction step}
		Let $W$ be a windmill.
		There exists a windmill $W'$ which contains the $\sigma/10$-neighborhood of $W$.
	\end{prop}
	
	\begin{proof}
		Let us denote by $V$ the following set of apices
		\begin{displaymath}
			V = \set{\fantomB v \in v(\mathcal R)\setminus W}{d(v,W) \leq3\sigma/10}.
		\end{displaymath}
		If $V$ is empty then the $\sigma/10$-neighborhood of $W$ is also a windmill.
		Therefore we may assume that $V$ is not empty.
		Note that $V$ is invariant under the action of $K_W$.
		We denote by $S$ (like \emph{sail}) the hull of $W \cup V$ (see Definition~\ref{def: hull}).
		Let $W'$ be the $\sigma/10$-neighborhood of $K_S.S$. 
		In particular $W'$ contains the $\sigma/10$-neighborhood of $W$. 
		The goal is to prove that $W'$ is a windmill.
		
		\begin{lemm}
		\label{res: windmill - apices in the neighborhood of S}
			Let $v \in v(\mathcal R)$.
			If $d(v,S) \leq \sigma/5$ then $v \in S$.
		\end{lemm}

		\begin{proof}
			By definition of hulls, a $\delta$-projection of $v$ on $S$ lies on a $(1,\delta)$-quasi-geodesics between two points $y,y'$ of $W\cup V$.
			In particular $\gro y{y'}v \leq d(v,S) + 2\delta \leq \sigma/5 + 2\delta$.
			We assume that $y,y'\in V$.
			The proof for the other cases works in the same way.
			Let us denote by $z$ and $z'$ respective $\delta$-projections of $y$ and $y'$ on $W$.
			By hyperbolicity
			\begin{equation}
			\label{eqn: windmill - apices in the neighborhood of K.S}
					\min \left\{ \fantomB \dist yv - \dist yz, \gro z{z'}v, \dist{y'}v - \dist {y'}{z'} \right\} \leq \gro y{y'}v +2 \delta \leq \frac \sigma 5 + 4\delta.
			\end{equation}
			Assume first that the minimum is achieved by $\dist yv - \dist yz$ (the proof works similarly for $\dist{y'}v - \dist {y'}{z'}$).
			By construction $\dist yz$ is bounded above by $3\sigma/10 + \delta$, thus $\dist yv < \sigma$.
			Nevertheless the distance between two distinct apices of $\mathcal R$ is at least $\sigma$. 
			Therefore $y=v$.
			Hence $v \in V$.
			Assume know that the minimum in~(\ref{eqn: windmill - apices in the neighborhood of K.S}) is achieved by $\gro z{z'}v$.
			The windmill $W$ being $2\delta$-quasi-convex, we have 
			\begin{displaymath}
				d(v,W) \leq \gro z{z'}v + 2 \delta \leq \frac \sigma 5 + 6\delta \leq \frac {3\sigma}{10}.
			\end{displaymath}
			By definition of $V$, $v$ is necessarily of point of $W \cup V$.
		\end{proof}
		
		\paragraph{}
		Recall that $\mathcal R$ is $G$-invariant.
		It follows from Lemma~\ref{res: windmill - apices in the neighborhood of S} that every apex contained in the $\sigma/5$-neighborhood of $K_S.S$ actually belongs to $S$.
		Since $W'$ is the $\sigma/10$ neighborhood of $K_S.S$, it satisfies \ref{enu: windmill - no apex in the neighborhood}.
		Moreover, all the apices contained in $W'$ lies in $K_S.S$.
		Hence $K_{W'} = K_S$.
		In particular $K_S.S$ and thus $W'$ are $K_S$-invariant.
		This corresponds to \ref{enu: windmill - K_W invariant}.
		
		\begin{lemm}
		\label{res: windmill - gromov product at the apex}
			Let $(H,v) \in \mathcal R$ such that $v \in V$. 
			Let $x,y \in S$ and $h \in H\setminus\{1\}$.
			Then $\gro x{hy}v \leq 7\delta$.
		\end{lemm}
		
		\begin{proof}
			Recall that $S$ is the hull of $W \cup V$. 
			According to Lemma~\ref{res: gromov product and hull}, it is sufficient to prove that for all $x,y \in W\cup V$, $\gro x{hy}v \leq 4\delta$.
			Let $x,y \in W\cup V$.
			Note that if $x=v$ or $y=v$, then $\gro x{hy}v = 0$ ($h$ fixes $v$).
			Therefore we can assume that $x$ and $y$ are distinct from $v$.
			We denote by $r$ a $\delta$-projection of $v$ on $W$.
			We claim that $\gro xvr \leq 7\delta$.
			If $x$ belongs to $W$, then by Proposition~\ref{res: proj quasi-convex}, $\gro xvr \leq 3\delta$.
			Assume now that $x$ is a point of $V\setminus\{v\}$. 
			Fix $p$ a $\delta$-projection of $x$ on $W$.
			By Proposition~\ref{res: proj quasi-convex}, 
			\begin{displaymath}
				\dist pr \leq \max \left\{\fantomB \dist xv - \dist xp - \dist vr + 14\delta, 7\delta \right\}.
			\end{displaymath}
			However $x$ and $v$ are two distinct apices of $V$. 
			It follows that $\dist xv \geq \sigma$ whereas  $\dist xp$ and $\dist vr$ are at most $3\sigma/10 + \delta$.
			By triangle inequality $\dist pr > 7\delta$.
			Consequently we necessarily have 
			\begin{displaymath}
				\dist xp + \dist pr + \dist rv \leq \dist xv + 14 \delta
			\end{displaymath}
			In particular $\gro xvr \leq 7\delta$, which proves our claim.
			Similarly $\gro yrv \leq 7 \delta$.
			Lemma~\ref{res: rotation family - small product at the apex} combined with the hyperbolicity condition leads to
			\begin{equation}
			\label{eqn: windmill - gromov product at the apex}
				\min\left\{ \fantomB \gro rxv,  \gro x{hy}v,  \gro {hy}{hr}v\right\} 
				\leq \gro r{hr}v + 2\delta
				\leq 4 \delta.
			\end{equation}
			Since $W$ is a windmill, all the apices of $\mathcal R$ in the $\sigma/10$-neighborhood of $W$ are actually contained in $W$. 
			If follows that $\dist vr > \sigma/10$.
			Hence $\gro xrv = \dist vr - \gro xvr$ is bounded below by  $\sigma/10-7\delta$.
			The same holds for $\gro {hy}{hr}v = \gro yrv$.
			Consequently the minimum in (\ref{eqn: windmill - gromov product at the apex}) is necessarily achieved by $\gro x{hy}v$.
			Hence $\gro x{hy}v \leq 4 \delta$.
		\end{proof}

		\begin{lemm}
		\label{res: windmill - pre quasi-convex}
			Let $y, y' \in S$. 
			Let $g \in K_S$.
			There exists a sequence of points $y=y_0, \dots, y_{m+1}=gy'$ of $X$ satisfying the following properties
			\begin{enumerate}
				\item \label{enu: windmill - pre quasi-convex - points in S}
				for all $i \in \intvald 1{m+1}$ there exists $g_i \in K_S$ such that $g_i^{-1}y_{i-1}$ and $g_i^{-1} y_i$ belong to $S$,
				\item \label{enu: windmill - pre quasi-convex - distance}
				for all $i \in \intvald 1{m-1}$, $\dist {y_{i+1}}{y_i} \geq \sigma$.
				\item \label{enu: windmill - pre quasi-convex - gromov produt}
				for all $i \in \intvald 1m$, $\gro {y_{i-1}}{y_{i+1}}{y_i} \leq 7 \delta$,
				\item \label{enu: windmill - pre quasi-convex - quasi-geodesic}
				For all $x \in X$ there exists $i \in \intvald 0m$ such that $\gro{y_{i+1}}{y_i}x \leq \gro y{gy'}x + 22\delta$.
			\end{enumerate}
			Moreover, if $m \leq 1$ then there exists $(H,v) \in \mathcal R$, such that $v \in V$ and $g \in H.K_W$.
		\end{lemm}
		
		\begin{proof}
			If $g \in K_W$ then the points $y_0 = y$ and $y_1 = gy'$ lie in $S$ and hence satisfy the conclusion of the lemma.
			Assume now that $g$ does not belong to $K_W$.
			The group $K_S$ is generated by $K_W$ and the rotation groups of the pairs $(H,v) \in \mathcal R$ where $v \in V$.
			It follows that $g$ can be written $g = u_0 h_1 u_1 \dots u_{m-1}h_m u_m$ where $m \in \N^*$ and
			\begin{enumerate}
				\item for all $i \in \intvald 0m$, $u_i$ belongs to $K_W$,
				\item for all $i \in \intvald 1m$,  there exists $(H_i,v_i) \in \mathcal R$ such that $v_i \in V$ and $h_i \in H_i \setminus\{1\}$,.
			\end{enumerate}
			We choose such a decomposition of $g$ which minimizes $m$.
			We claim that for all $i \in \intvald 1{m-1}$, $u_iv_{i+1}\neq v_i$.
			Assume on the contrary that this assertion is false.
			Using the action of $G$ on $\mathcal R$, $u_iH_{i+1}u_i^{-1} =H_i$.
			Therefore we can write $h_i u_i h_{i+1} u_{i+1} = \left( h_i u_ih_{i+1}u_i^{-1}\right)\left( u_i u_{i+1}\right)$, where $h_i u_ih_{i+1}u_i^{-1}\in H_i$ and $u_i u_{i+1} \in K_W$.
			This leads to a shorter decomposition of $g$. 
			Contradiction.
			For all $i \in \intvald 1m$, we put $g_i = u_0 h_1 u_1 \dots h_{i-1}u_{i-1} $ and $y_i = g_iv_i$.
			Moreover, we put $g_{m+1}=g$, $y_0 = y$ and $y_{m+1} = gy'$.
			Note that, if $m= 1$ then $g$  can be written $g = u_0h_1u_1 = \left(u_0h_1u_0^{-1}\right)\left(u_0u_1\right)$, where $u_0,u_1 \in K_W$ and $h_1 \in H_1$.
			However the set $V$ is invariant under the action of $K_W$.
			Therefore $u_0h_1u_0^{-1}$ is an element of the rotation group $u_0H_1u_0^{-1}$, whose apex $u_0v_1$ belongs to $V$.
			This proves the last assertion of the lemma.
			
			\paragraph{}Let $i \in \intvald 2{m+1}$. 
			By construction $g_i^{-1}y_i$ is a point of $S$. 
			By definition of rotation family $h_{i-1}$ fixes the apex $v_{i-1}$, hence $g_i^{-1}y_{i-1} = u_{i-1}^{-1}h_{i-1}^{-1}v_{i-1}  =  u_{i-1}^{-1}v_{i-1}$.
			However $V$ is invariant under the action of $u_{i-1} \in K_W$, thus $g_i^{-1}y_{i-1}$ belongs to $S$.
			On the other hand $y_0 = y$ and $y_1 = u_0v_1$ belong to $S$.
			This completes the proof of Point~\ref{enu: windmill - pre quasi-convex - points in S}.
			
			\paragraph{}Let $i \in \intvald 1{m-1}$.
			The apex $v_i$ is fixed by $h_i$ therefore
			\begin{displaymath}
				\dist {y_{i+1}}{y_i} = \dist{g_ih_iu_iv_{i+1}}{g_iv_i} = \dist {u_iv_{i+1}}{v_i}.
			\end{displaymath}
			However, we explained that $u_i v_{i+1}$ and $v_i$ are necessarily two distinct apices of $\mathcal R$, therefore $\dist{y_{i+1}}{y_i}\geq \sigma$.
			This proves Point~\ref{enu: windmill - pre quasi-convex - distance}.

			\paragraph{}Let $i \in \intvald 1m$. By construction $g_i^{-1} y_i=v_i$ whereas $g_i^{-1}y_{i-1}$ belongs to $S$.
			On the other hand $g_i^{-1}y_{i+1} = g_i^{-1}g_{i+1}g_{i+1}^{-1}y_{i+1} = h_iu_ig_{i+1}^{-1}y_{i+1}$.
			However $u_i$ stabilizes $S$ which contains $g_{i+1}^{-1}y_{i+1}$, hence $g_i^{-1}y_{i+1}$ belongs to $h_iS$.
			By Lemma~\ref{res: windmill - gromov product at the apex}, $\gro {y_{i-1}}{y_{i+1}}{y_i} = \gro{g_i^{-1}y_{i-1}}{g_i^{-1}y_{i+1} }{v_i}$ is bounded above by $7 \delta$, which proves Point~\ref{enu: windmill - pre quasi-convex - gromov produt}.
			We chose the constants $\sigma$ and $\delta$ in such a way that we can apply Proposition~\ref{res: stability discrete quasi-geodesic} to the sequence $y_0,\dots,y_{m+1}$.
			Point~\ref{enu: windmill - pre quasi-convex - quasi-geodesic} follows from the stability of discrete quasi-geodesics.
		\end{proof}

		\begin{lemm}
		\label{res: windmill - W' quasi-convex}
			The set $K_S.S$ is $28\delta$-quasi-convex whereas $W'$ is $2\delta$-quasi-convex.
		\end{lemm}
		
		\begin{proof}
			The set $W'$ was defined as the $\sigma/10$-neighborhood of $K_S.S$.
			According to Lemma~\ref{res: neighborhood of a quasi-convex}, it is sufficient to show that $K_S.S$ is $28\delta$-quasi-convex.
			Let $x \in X$ and $y,y' \in K_S.S$.
			It follows from Lemma~\ref{res: windmill - pre quasi-convex} that there exist $z,z' \in S$ and $g \in K_S$ such that $\gro {gz}{gz'}x \leq \gro y{y'}x + 22\delta$.
			However $S$ being a hull, it is $6\delta$-quasi-convex.
			Therefore 
			\begin{displaymath}
				d(g^{-1}x,K_S.S) \leq d(g^{-1}x,S) \leq \gro {gz}{gz'}x +6\delta\leq \gro y{y'}x + 28\delta.
			\end{displaymath}
			By construction $K_S.S$ is invariant under the  action of $K_S$.
			It follows that  $d(x,K_S.S) \leq \gro y{y'}x + 28\delta$.
		\end{proof}
		
		\begin{lemm}
		\label{res: windmill - pre translation length}
			Let $(H,v) \in \mathcal R$ such that $v \in V$.
			Let $h \in H\setminus\{1\}$ and $u \in K_W\setminus \{1\}$.
			For all $y \in S$, $\dist {huy}y \geq \sigma - 16\delta$.
		\end{lemm}
				
		\begin{proof}
			We need to distinguish two cases.
			
			\paragraph{Case 1.} \emph{There exists $(H',v') \in \mathcal R$ such that $v' \in W$ and $u \in H'\setminus\{1\}$.}
			By Lemma~\ref{res: rotation family - small product at the apex}, $\gro {uy}y{v'}\leq 2\delta$ i.e., $\dist y{v'} = \dist{uy}{v'}$ is bounded above by $\dist{uy}y /2 + 2\delta$.
			The triangle inequality yields
			\begin{eqnarray*}
				\dist v{v'}
				& \leq & \min\left\{ \dist vy +\dist y{v'}, \dist v{uy} +\dist{uy}{v'}  \right\} \\
				& \leq & \min\left\{ \dist vy, \dist v{uy} \right\} +\frac 12 \dist{uy}y + 2\delta \\
				& \leq & \dist {uy}v + \dist vy +2\delta.
			\end{eqnarray*}
			However $v$ ad $v'$ are two distinct apices of $\mathcal R$, thus $\dist v{v'} \geq \sigma$.
			It follows that $ \dist {uy}v + \dist vy$ is at least $\sigma -2\delta$.
			Recall that $y$ and $uy$ are two points of $S$.
			By Lemma~\ref{res: windmill - gromov product at the apex}, we obtain
			\begin{displaymath}
				\dist {huy}y \geq \dist {huy}v + \dist vy - 14\delta = \dist {uy}v + \dist vy - 14\delta \geq \sigma -16 \delta.
			\end{displaymath}

			\paragraph{Case 2.} \emph{Assume that $u$ does not belong to a rotation group.} 
			Let us denote by $p$ and $r$ respective $\delta$-projections of $y$ and $v$ on $W$.
			In particular $up$ is a $\delta$-projection of $uy$ on $W$.
			By projection on a quasi-convex (Lemma~\ref{res: proj quasi-convex}) we have
			\begin{eqnarray*}
				\dist {up}r & \leq  & \max\left\{ \dist{uy}v - \dist vr - \dist yp + 14 \delta ,7 \delta \right\}\\ 
				 \dist pr & \leq &  \max \left\{ \dist yv - \dist vr - \dist yp + 14 \delta ,7 \delta \right\} 
			\end{eqnarray*}
			Since $W$ is a windmill and $u \in K_W$, then $\dist {ur}p + \dist pr \geq \dist {ur}r \geq  \sigma - 166\delta$.
			Therefore the two previous maxima cannot be both achieved by $7\delta$. 
			Assume for instance that the first maximum is not  achieved by $7\delta$ (the other case is symmetric).
			By Lemma~\ref{res: rotation family - small product at the apex}, we get $\dist{huy}y \geq \dist{uy}v + \dist vy - 14\delta$ which leads to 
			\begin{displaymath}
				\dist{huy}y \geq \dist yp + \dist {up}r + \dist rv + \dist vy - 28\delta.
			\end{displaymath}
			However $r$ is a $\delta$-projection of $v$ on $W$, thus
			\begin{displaymath}
				\dist vy+\dist yp \geq \dist vp \geq \dist rp - \gro pvr \geq \dist rp - 3\delta.
			\end{displaymath}
			Consequently the previous inequality becomes
			\begin{displaymath}
				\dist{huy}y \geq  \dist {up}r +\dist rp + \dist rv - 31 \delta \geq \dist{up}p + \dist rv- 31\delta.
			\end{displaymath}
			Nevertheless Axiom~\ref{enu: windmill - almost free action} for $W$ gives $\dist {up}p \geq  \sigma-166\delta$.
			Moreover, by construction $\dist vr \geq d(x,W) \geq \sigma/10$.
			Thus $ \dist{huy}y \geq \sigma$.
		\end{proof}
		
		\begin{lemm}
		\label{res: windmill - translation length}
			Let $g \in K_S$ such that for every $H \in H(\mathcal R)$, $g$ does not belong to $H$.
			For all $x \in X$, $\dist {gx}x \geq \sigma - 166\delta$.
		\end{lemm}
		
		\begin{proof}
			Let $x \in X$.
			Since $W$ is already a windmill we can assume that $g$ does not belong to $K_W$, otherwise the result would follow from \ref{enu: windmill - almost free action}.
			We denote by $y$ a $\delta$-projection of $x$ on $K_S.S$.
			The rotation family $\mathcal R$ is invariant under the action of $G$.
			Without lost of generality we may assume that $y \in S$.
			We claim that $\dist {gy}y \geq \sigma - 48\delta$.
			According to Lemma~\ref{res: windmill - pre quasi-convex} there exists a sequence of points $y = y_0, \dots , y_{m+1} = gy$ satisfying the following conditions
			\begin{enumerate}
				\item for all $i \in \intvald 1{m-1}$, $\dist {y_{i+1}}{y_i} \geq \sigma$.
				\item for all $i \in \intvald 1m$, $\gro {y_{i-1}}{y_{i+1}}{y_i} \leq 7 \delta$,
			\end{enumerate}
			Assume that $m \geq 2$.
			The constant $\sigma$ has been chosen to apply the stability of discrete quasi-geodesic.
			By Proposition~\ref{res: stability discrete quasi-geodesic}, $\gro y{gy}{y_1}$ and $\gro{y_1}{gy}{y_2}$ are at most $12\delta$.
			It follows that 
			\begin{displaymath}
				\dist{gy}y \geq \dist {gy}{y_2} +\dist{y_2}{y_1} + \dist{y_1}y -108\delta \geq \sigma- 48\delta.
			\end{displaymath}
			Assume now that $m \leq 1$.
			According to Lemma~\ref{res: windmill - pre quasi-convex}, there exists $(H,v) \in \mathcal R$ such that $v \in V$ and $g \in H.K_W$.
			In particular $g$ can be written $g = hu$ where $h \in H$ and $u \in K_W$.
			Since $g$ does not belong to $K_W$, $h$ is non-trivial.
			By assumption $g$ does not belong to $H$ thus $u \neq 1$.
			Applying Lemma~\ref{res: windmill - pre translation length}, $\dist{gy}y \geq \sigma - 16\delta$, which completes the proof of our claim.
			By Lemma~\ref{res: windmill - W' quasi-convex}, $K_S.S$ is $28\delta$-quasi-convex. 
			Applying Lemma~\ref{res: proj quasi-convex} we obtain
			\begin{displaymath}
				\dist {gx}x \geq \dist {gx}{gy} + \dist {gy}y +  \dist yx  - 118\delta \geq  \sigma - 166\delta. \qedhere
			\end{displaymath}
		\end{proof}

		We already proved that Axioms \ref{enu: windmill - K_W invariant} and \ref{enu: windmill - no apex in the neighborhood} for $W'$ follow from Lemma~\ref{res: windmill - apices in the neighborhood of S}.
		Axioms \ref{enu: windmill - quasi convex} and \ref{enu: windmill - almost free action} respectively correspond to Lemmas~\ref{res: windmill - W' quasi-convex} and \ref{res: windmill - translation length}.
		Hence $W'$ is a windmill.
	\end{proof}
	
	\begin{proof}[Proof of Theorem~\ref{res: rotation family - fundamental theorem}]	
		Let $g \in K$.
		We choose an apex $v \in v(\mathcal R)$.
		The set $\{v\}$ is a windmill.
		Iterating Proposition~\ref{res: windmill - induction step}, we obtain a windmill $W$ containing sufficiently many apices of $v(\mathcal R)$ so that $g \in K_W$.
		It follows from \ref{enu: windmill - almost free action} that if  $g$ does not belong to any rotation group $H \in H(\mathcal R)$ then for every $x \in X$ $\dist {gx}x \geq \sigma -166\delta$.
	\end{proof}
	
	\begin{coro}
	\label{res: rotation family - free outside the apices}
		Let $l \geq 0$.
		Let $x\in X$ such that for every apex $v \in v(\mathcal R)$, $\dist xv \geq l$.
		Let $g \in K\setminus\{1\}$.
		Then $\dist {gx}x \geq \min\left\{ 2l, \sigma/10 \right\}$.
	\end{coro}
	
	\rem In particular the group $K$ acts freely discontinuously on the space $X \setminus v(\mathcal R)$.
	By freely discontinuously we mean that for every $x \in X \setminus v(\mathcal R)$ there exists $r >0$ such that for every $g \in K$ if $gB(x,r)$ intersects $B(x,r)$ then $g$ is trivial.

	\begin{proof}
		If $g$ does not belong to a rotation group then by Theorem~\ref{res: rotation family - fundamental theorem}, $\dist {gx}x \geq  \sigma -166\delta$.
		Therefore we can assume that there exists $(H,v) \in \mathcal R$ such that $g \in H\setminus\{1\}$.
		If $\dist xv \geq \sigma/10$, then by Lemma~\ref{res: rotation family - small product at the apex}, $\dist {gx}x \geq 2 \dist xv - 4\delta \geq \sigma/10$.
		Otherwise the definition of rotation family yields $\dist {gx}x = 2\dist xv \geq 2l$.
	\end{proof}
	
	\begin{coro}
	\label{res: rotation family - stabilizer in K}
		Let $(H,v) \in \mathcal R$.
		We have $\stab v\cap K = H$.
		Moreover, for every $x \in B(v,\sigma/5)$, for every $g \in K\setminus H$ we have $\dist {gx}x > 3\sigma/5$.
	\end{coro}

	\begin{proof}
		By construction $H$ lies in $\stab v \cap K$.
		Let $g \in \stab v \cap K$ which is not trivial.
		In particular $gv=v$.
		By Theorem~\ref{res: rotation family - fundamental theorem}, there exists $(H',v') \in \mathcal R$ such that $g \in H' \setminus \{1\}$.
		Assume now that $v \neq v'$.
		By Lemma~\ref{res: rotation family - small product at the apex}, $\dist {gv}v \geq 2\dist v{v'} - 4\delta > 0$.
		Contradiction.
		Hence $v = v'$ and $g$ belongs to $H$.
		
		\paragraph{} Let us consider now $x \in B(v,\sigma/5)$ and $g\in K\setminus H$.
		Assume that our second assertion is false.
		Then by triangle inequality $\dist {gv}v \leq \dist{gx}x + 2\dist xv <\sigma$.
		However the distance between two distinct apices is at least $\sigma$.
		Therefore $g$ belongs to $\stab v \cap K = H$.
		Contradiction.
	\end{proof}

	\subsection{Consequences}
	
	We keep here the notations and assumptions of the previous section.
	In particular, $G$ is a group acting by isometries on the $\delta$-hyperbolic length space $X$ and $\mathcal R$ is a $\sigma$-rotation family with $\sigma\geq \sigma_0$ where $\sigma_0$ is the parameter given by Theorem~\ref{res: rotation family - fundamental theorem}.
	The space $\bar X$ is defined to be the quotient of $X$ by $K$.
	If $x$ is a point of $X$ we denote by $\bar x$ its image by the canonical map $\nu : X \rightarrow \bar X$.
	Moreover, we write $\bar v(\mathcal R)$ for the image of $v(\mathcal R)$ in $\bar X$.
	We endow $\bar X$ with the pseudo-metric defined in the following way. 
	\begin{displaymath}
		\forall x,x' \in X, \quad \dist[\bar X]{\bar x}{\bar x'} = \inf_{g \in K} \dist[X]{gx}{x'}.
	\end{displaymath}
	By construction the quotient $\bar G = G/K$ acts by isometries on $\bar X$.
	Given an element $g$ of $G$ we denote by $\bar g$ its image by the projection $\pi : G \twoheadrightarrow \bar G$.
	It follows from Corollary~\ref{res: rotation family - stabilizer in K} that for every $(H,v) \in \mathcal R$ this map induces an embedding $\stab v/H \hookrightarrow \bar G$.
	Furthermore $X\setminus v(\mathcal R)$ is a covering space of $\bar X \setminus \bar v(\mathcal R)$.
	
	\begin{prop}
	\label{res: rotation family - bar X metric space}
		$\bar X$ is a metric length space.
	\end{prop}
	
	\begin{proof}
		One only needs to prove that the pseudo-distance on $\bar X$ is definite positive. 
		The fact that $\bar X$ is a length space follows from the length structure on $X$ \cite[Chap. I.5, Lemma 5.20]{BriHae99}.
		Let $x$ and $x'$ be two points of $X$ such that $\dist{\bar x}{\bar x'} = 0$.
		Assume first that $x$ is not an apex.
		The set of apices $v(\mathcal R)$ is closed and $G$-invariant, thus neither is $x'$.
		Since $K$ acts freely discontinuously on $X\setminus v(\mathcal R)$, $x$ and $x'$ are necessarily in the same $K$-orbit.
		Assume now that $x$ and $x'$ are both apices.
		By definition of the pseudo-metric there exists $g \in K$ such that $\dist {gx'}{x} < \sigma$.
		However the distance between two distinct apices is at least $\sigma$.
		Hence $gx = x'$.
		Consequently, in both cases $\bar x = \bar x'$.
	\end{proof}

	\begin{prop}
	\label{res: rotation family - nu local isometry}
		Let $r \in(0,\sigma/40]$.
		Let $x \in X$ such that for all $v \in v(\mathcal R)$, $\dist xv \geq 2r$.
		The map $\nu: X \rightarrow \bar X$ induces an isometry from $B(x, r)$ onto $B(\bar x, r)$.
	\end{prop}
	
	\rem In particular, the map $\nu : X \rightarrow \bar X$ induces a local isometry from $X\setminus v(\mathcal R)$ onto $\bar X\setminus\bar v(\mathcal R)$.
	
	\begin{proof}
		By construction the ball $B(x,2r)$ is contained in $X\setminus v(\mathcal R)$.
		According to Corollary~\ref{res: rotation family - free outside the apices}, $\nu$ induces a bijection from $B(x,2r)$ onto its image.
		It follows that it induces an isometry from $B(x, r)$ onto its image which is exactly $B(\bar x, r)$.
	\end{proof}	
	
	\begin{lemm}
	\label{res: rotation family - bar X hyperbolic around the apices}
		For every $v \in v(\mathcal R)$ the ball $B(\bar v, \sigma/5)$ of $\bar X$ is $2\delta$-hyperbolic.
	\end{lemm}
	
	\begin{proof}
		First note that for every point $x \in B(v,\sigma/5)$ we have $ \dist{\bar x}{\bar v}=\dist xv$ (this is a consequence of Corollary~\ref{res: rotation family - stabilizer in K}).
		Let $\bar x$, $\bar y$ and $\bar z$ be three points of $B(\bar v, \sigma/5)$.
		We denote by $x$, $y$ and $z$ respective pre-images of $\bar x$, $\bar y$ and $\bar z$ in $B(v,\sigma/5)$.
		By definition of the metric there exists two sequences $(g_n)$ and $(h_n)$ of elements of $K$ such that $\dist{g_nx}{y}$ and $\dist{h_nz}{y}$ respectively converge to $\dist{\bar x}{\bar y}$ and $\dist {\bar z}{\bar y}$ as $n$ approaches infinity.
		In particular if $n$ is sufficiently large we get by triangle inequality $\dist{g_nv}{v}<\sigma$.
		It follows that $g_n$ belongs to the stabilizer of $v$ and thus $g_n x$ also lies in $B(v,\sigma/5)$.
		In the same way for $n$ sufficiently large $h_nz$ belongs to $B(v,\sigma/5)$.
		However $X$ is $\delta$-hyperbolic, therefore $g_nx$, $y$, $h_nz$ and $v$ satisfy the four points inequality~(\ref{eqn: hyperbolicity condition 1}).
		Consequently
		\begin{displaymath}
			\gro{\bar x}{\bar z}{\bar v} \geq \gro{g_nx}{h_nz}v \geq \min\left\{ \gro{g_nx}{y}v, \gro y{h_nz}v\right\} -\delta.
		\end{displaymath}
		Taking the limit as $n$ approaches infinity we obtain
		\begin{displaymath}
			\gro{\bar x}{\bar z}{\bar v} \geq  \min\left\{ \gro{\bar x}{\bar y}{\bar v}, \gro {\bar y}{\bar z}{\bar v}\right\} -\delta.
		\end{displaymath}
		By Lemma~\ref{res: criterion for hyperbolicity}, $B(\bar x,\sigma/5)$ is $2\delta$-hyperbolic.
	\end{proof}
	
	\begin{lemm}
	\label{res: rotation family - bar X delta-simply connected}
		The space $\bar X$ is $50\delta$-simply connected i.e., its fundamental group is normally generated by free homotopies of loops of diameter at most $50\delta$.
	\end{lemm}
	
	\begin{proof}
		Let $\bar \gamma$ be a loop in $\bar X$ based at $\bar x$ and $\gamma_1$ a lift of $\bar \gamma$ in $X$.
		If $x$ is the initial point of $\gamma$ there exists $g \in K$ such that the terminal point of $\gamma_1$ is $gx$.
		Note that $K$ is generated by elliptic isometries of $X$.
		Moreover, the translation length of an elliptic isometry of $X$ is at most $32\delta$.
		Therefore there exists an other path $\gamma_2$ of $X$ joining $x$ to $gx$ such that its image in $\bar X$ can be written as a product of loops of diameter at most $32\delta$.
		The space $X$ being $\delta$-hyperbolic its Rips complex $P_{4\delta}(X)$ is simply-connected \cite[Chap. 5, Prop. 1.1]{CooDelPap90}.
		It implies that $X$ is $50\delta$-simply-connected.
		On the other hand $\gamma_2^{-1}\gamma_1$ is a loop of $X$.
		Consequently it can be written in $X$ and thus in $\bar X$ as a product of loops of diameter at most $50\delta$.
		Hence $\bar X$ is $50\delta$-simply-connected.
	\end{proof}

	\begin{prop}
	\label{res: rotation family - bar X globally hyperbolic}
		$\bar X$ is $600\delta$-hyperbolic
	\end{prop}
	
	\begin{proof}
		According to Propositions~\ref{res: rotation family - nu local isometry} and \ref{res: rotation family - bar X hyperbolic around the apices} every ball $\bar B$ of radius $\sigma/40$ in $\bar X$ is $2\delta$-hyperbolic.
		Indeed if the distance between center of $\bar B$ and every apex is at least $\sigma/10$ then it is isometric to a ball of $X$ which is $\delta$-hyperbolic.
		Otherwise there is $v \in v(\mathcal R)$ such that $\bar B$ lies in $B(\bar v,\sigma/5)$ which is $2\delta$-hyperbolic.
		On the other hand, by Proposition~\ref{res: rotation family - bar X delta-simply connected}, $\bar X$ is $50\delta$-simply-connected.
		Recall that we chose $\sigma \geq \sigma_0$ with $\sigma_0 > 10^{20}\delta$. 
		Therefore we can apply the Cartan-Hadamard Theorem (see Theorem~\ref{res: cartan hadamard}).
		It follows that $\bar X$ is $600 \delta$-hyperbolic.
	\end{proof}
	
	\paragraph{} We now explain how to lift figures of $\bar X$ in $X$. 
	To that end, we first introduce a very general construction, coming from the theory of covering spaces, to lift a path.
	Let $\bar x$ be a point of $\bar X$ and $x$ a preimage of $\bar x$ in $X$.
	Let $\bar \gamma : I \rightarrow \bar X$ be a path starting at $\bar x$.
	We assume that $\bar \gamma$ does not go through $\bar v(\mathcal R)$.
	Since $K$ acts freely discontinuously on $X\setminus v(\mathcal R)$, there exists a unique path $\gamma : I \rightarrow X$ starting at $x$ and lifting $\bar \gamma$.
	We write $x + \bar \gamma$ for the terminal point of $\gamma$.
	It is a preimage of the terminal point of $\bar \gamma$ that we denote $\bar x + \bar \gamma$.
	If $\bar \gamma_1$ and $\bar \gamma_2$ are to paths with the same extremities and homotopic relative to their endpoints in $\bar X\setminus \bar v(\mathcal R)$ then $x + \bar \gamma_1$ and $x+\bar \gamma_2$ define the same point of $X$.
	
	\paragraph{}
	In \cite{DelGro08} the authors make an intensive use of this topological point of view. 
	They need nevertheless to work with orbifolds to deal with the torsion.
	For our purpose we prefer a more geometrical approach that follows from the rotation families.
	The statements in the remainder of this section corresponds to Lemme~5.9.4 and Lemme~5.10.1 of \cite{DelGro08}.
	
	\paragraph{}Since the map $X\setminus v(\mathcal R) \rightarrow \bar X \setminus \bar v(\mathcal R)$ is a local isometry, $\bar \gamma$ and its lift $\gamma$ have the same length.
	If $\bar \gamma$ is a $(1,l)$-quasi-geodesic, then
	\begin{displaymath}
		L(\gamma) = L(\bar \gamma) \leq \dist {(\bar x +\bar \gamma)}{\bar x}+l \leq \dist {(x + \bar \gamma)}x + l.
	\end{displaymath}
	Here $L(\gamma)$ and $L(\bar \gamma)$ stand for the respective lengths of the paths $\gamma$ and $\bar \gamma$.
	In particular $\gamma$ is also a $(1,l)$-quasi-geodesic.
	Moreover, we have $\dist {(x + \bar \gamma)}x \leq \dist {(\bar x +\bar \gamma)}{\bar x}+l$.
	
\begin{figure}[htbp]
	\centering
		\hfill
		\subfigure[Triangle in the quotient $\bar X$.]{
			\includegraphics[width=0.35\textwidth]{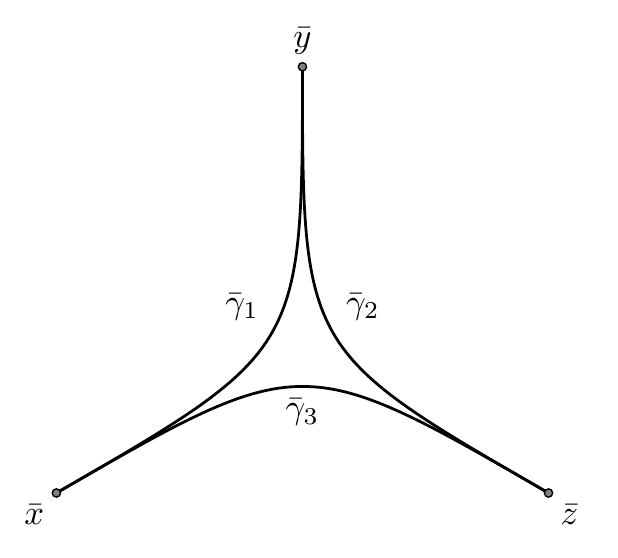}
			\label{fig: lifting triangle - quotient}
		}
		\hfill
		\subfigure[Lift in the space $X$.]{
			\includegraphics[width=0.45\textwidth]{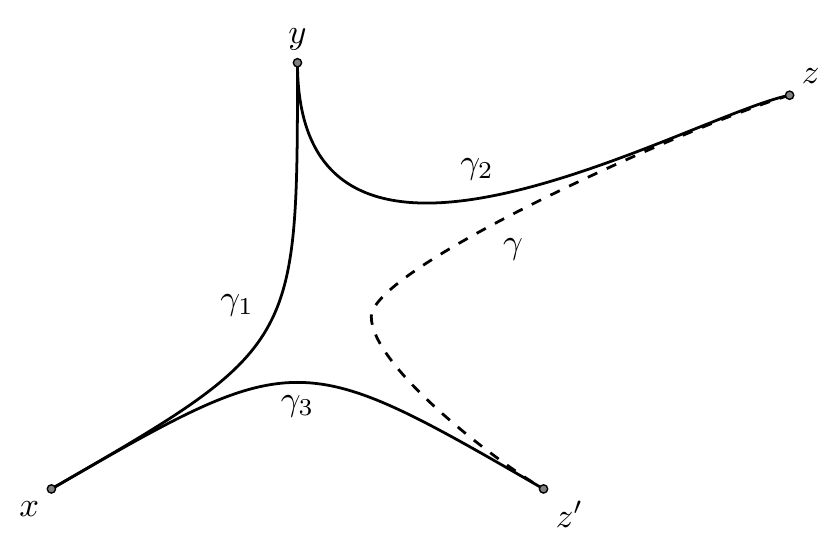}
			\label{fig: lifting triangle - coneoff}
		}
		\hfill
	\caption{Lifting a triangle}
	\label{fig: lifting triangle}
	\end{figure}

	\begin{lemm}
	\label{res: left quasi-geodesic and homotopies}
		Let $\bar x$ be a point of $\bar X$ and $x$ a preimage of $\bar x$ in $X$.
		Let $\bar \gamma_1$, $\bar \gamma_2$ and $\bar \gamma_3$ be $(1,\delta)$-quasi-geodesic paths satisfying the following:
		\begin{enumerate}
			\item $\bar \gamma_1$ and $\bar \gamma_3$ start at $\bar x$,
			\item $\bar \gamma_2$ starts at $\bar x + \bar \gamma_1$ and ends at $\bar x + \bar \gamma_3$; see Figure~\ref{fig: lifting triangle - quotient}.
		\end{enumerate}
		Assume the for every $v \in v(\mathcal R)$, these paths do not enter the ball $B(\bar v, \sigma/10 + 7\delta)$.
		Then $(x+\bar \gamma_1)+\bar \gamma_2 = x +\bar \gamma_3$.
	\end{lemm}
	
	\begin{proof}
		For simplicity of notation, let us respectively denote by $y$, $z$ and $z'$ the points $x+\bar \gamma_1$, $y+\bar \gamma_2$ and $x +\bar \gamma_3$.
		See Figure~\ref{fig: lifting triangle - coneoff}.
		In particular $z$ and $z'$ have the same image in $\bar X$.
		We write $\gamma_1$ (\resp $\gamma_2$, $\gamma_3$) for the lift of $\bar \gamma_1$ (\resp $\bar \gamma_2$, $\bar \gamma_3$) which starts at $x$ (\resp $y$, $x$).
		These paths are $(1,\delta)$-quasi-geodesics.
		Moreover, for all $v \in v(\mathcal R)$ they do not enter $B(v,\sigma/10 + 7\delta)$.
		Let $\eta \in (0,\delta)$.
		We denote by $\gamma$ a $(1,\eta)$-quasi-geodesic of $X$ joining $z$ to $z'$.
		By hyperbolicity, $\gamma$ is contained in the $(\eta+6\delta)$-neighborhood of $\gamma_1 \cup \gamma_2 \cup \gamma_3$.
		Consequently, if $\eta$ is sufficiently small, $\gamma$ does not intersect any of the balls $B(v,\sigma/10)$ where $v \in v(\mathcal R)$.
		According to Proposition~\ref{res: rotation family - nu local isometry} the path $\bar \gamma$, image of $\gamma$ in $\bar X$, is a $\sigma/25$-local $(1,\eta)$-quasi-geodesic whose length is the same as the one of $\gamma$.
		However we chose $\sigma$ in such a way that we can apply the stability of quasi-geodesics.
		In particular $\bar \gamma$ is  a (global) $(2,\eta)$-quasi-geodesic.
		Its endpoints are the same, namely $\bar z = \bar z'$.
		Therefore its length is at most $\eta$.
		Consequently $\dist z{z'} \leq L(\gamma) \leq \eta$.
		This inequality holds for every $\eta >0$, hence $z = z'$ i.e., $(x+\bar \gamma_1)+\bar \gamma_2 = x +\bar \gamma_3$.
	\end{proof}
	
	\begin{prop}
	\label{res: lift quasi-convex isometry}
		Let $\alpha \geq 0$.
		Let $\bar Y$ be a $\alpha$-quasi-convex subset of $\bar X$ such that for every $v \in v(\mathcal R)$,  $\bar Y$ does not intersect $B(\bar v, \sigma/10 + \alpha +8\delta)$.
		Let $\bar y_0$ be a point of $\bar Y$ and $y_0$ a preimage of $\bar y_0$ in $X$.
		There exists a subset $Y$ of $X$ containing $y_0$ such that the map $\nu : X \rightarrow \bar X$ induces an isometry from $Y$ onto $\bar Y$.
	\end{prop}
	
	\begin{proof}
		We construct $Y$ in the following way.
		Let $\bar y$ be a point of $\bar Y$.
		We choose a $(1,\delta)$-quasi-geodesic $\bar \gamma$ joining $\bar y_0$ to $\bar y$.
		Since $\bar Y$ is $\alpha$-quasi-convex, $\bar \gamma$ lies in the $(\alpha + \delta)$-neighborhood of $\bar Y$. 
		In particular, it does not go through a vertex of $\bar X$.
		More precisely $\bar \gamma$ does not run through any ball $B(\bar v, \sigma/10  +7\delta)$ where $v \in v(\mathcal R)$.
		Consequently, $y_0 + \bar \gamma$ is well-defined.
		We define $Y$ to be the set of all points $y_0 +\bar \gamma$ obtained in this way.
		By construction $\nu$ maps $Y$ onto $\bar Y$.
		Hence it is sufficient to prove that the restriction of $\nu$ to $Y$ preserves the distances.
		Let $y_1$ and $y_2$ be two points of $Y$.
		By definition there exist $(1,\delta)$-quasi-geodesics $\bar \gamma_1$ and $\bar \gamma_2$ respectively joining $\bar y_0$ to $\bar y_1$ and $\bar y_2$ such that $y_1 = y_0 + \bar \gamma_1$ and $y_2 = y_0 + \bar \gamma_2$.
		Let $\eta \in (0,\delta)$ and $\bar \gamma$ be a $(1,\eta)$-quasi-geodesic joining $\bar y_1$ to $\bar y_2$.
		According to Lemma~\ref{res: left quasi-geodesic and homotopies}, $(y_0 +\bar \gamma_1) +\bar \gamma = y_0 + \bar \gamma_2$ i.e., $y_1+\bar \gamma = y_2$.
		Hence for every sufficiently small $\eta>0$, $\dist {y_1}{y_2} \leq \dist{\bar y_1}{\bar y_2}+\eta$.
		It follows that $\dist {y_1}{y_2} = \dist{\bar y_1}{\bar y_2}$.
	\end{proof}
	
	\begin{prop}
	\label{res: lift quasi-convex isometry and stabilizer}
		Let $\alpha \geq 0$ and $d \geq \alpha$.
		Let $\bar Y$ be an $\alpha$-quasi-convex subset of $\bar X$ such that for every $v \in v(\mathcal R)$,  $\bar Y$, does not intersect $B(\bar v, \sigma/10 + d +1208\delta)$.
		Let $\bar y_0$ be a point of $\bar Y$ and $y_0$ a preimage of $\bar y_0$ in $X$.
		There exists a subset $Y$ of $X$ containing $y_0$ with the following properties:
		\begin{enumerate}
			\item \label{enu: lift quasi-convex isometry and stabilizer - isometry}
			the map $\nu : X \rightarrow \bar X$ induces an isometry from $Y$ onto $\bar Y$,
			\item \label{enu: lift quasi-convex isometry and stabilizer - stabilizer}
			for every $\bar g \in \bar G$ such that $\bar g \bar Y$ lies in the $d$-neighborhood of $\bar Y$ there exists a preimage $g \in G$ of $\bar g$ with the following property. For every $y,z \in Y$, $\dist{gy}z = \dist{\bar g\bar y }{\bar z}$.
		\end{enumerate}
		In particular the projection $\pi : G \rightarrow \bar G$ induces an isomorphism from $\stab Y$ onto $\stab {\bar Y}$.
	\end{prop}
	
	\begin{proof}
		We denote by $\bar Z$ the $d$-neighborhood of $\bar Y$.
		Recall that $\bar X$ is $600 \delta$-hyperbolic.
		By Proposition \ref{res: neighborhood of a quasi-convex}, $\bar Z$ is $1200\delta$-quasi-convex, thus it satisfies the assumptions of Proposition~\ref{res: lift quasi-convex isometry}.
		As in this proposition we construct a subset $Z$ of $X$ containing $y_0$ such that the map $\nu : X \rightarrow \bar X$ induces an isometry from $Z$ onto $\bar Z$.
		We write $Y$ for the preimage of $\bar Y$ in $Z$.
		In particular $\nu$ maps $Y$ isometrically onto $\bar Y$.
		Let $\bar g \in G$ such that $\bar g \bar Y \subset \bar Z$.
		By construction, there exists $g \in G$ such that $gy_0$ is the unique preimage of $\bar g\bar y_0$ in $Z$.
		Let $y \in Y$.
		By assumption $\bar g \bar y$ is point of $\bar Z$.
		We claim that $gy$ is the (unique) preimage of $\bar g \bar y$ in $Z$.
		For simplicity of notation we put $\bar y_1 = \bar g \bar y_0$ and $\bar y_2 = \bar g \bar y$.
		We denote by $y_1  = g y_0$ and $y_2$ their respective preimages in $Z$.
		There exists a $(1,\delta)$-quasi-geodesic $\bar \gamma$ (\resp $\bar \gamma_1$, $\bar \gamma_2$) joining $\bar y_0$ to $\bar y$ (\resp $\bar y_1$, $\bar y_2$) such that $y = y_0 +\bar \gamma$, (\resp $y_1 = y_0 +\bar \gamma_1$, $y_2 = y_0 + \bar \gamma_2$).
		Note that $\bar g \bar \gamma$ is a $(1,\delta)$-quasi-geodesic joining $\bar y_1 = \bar g \bar y_0$ to $\bar y_2 = \bar g \bar y$.
		Hence by Lemma~\ref{res: left quasi-geodesic and homotopies}, $(y_0 +\bar \gamma_1) + \bar g \bar \gamma = y_0 +\bar \gamma_2$ i.e., $g y_0 +\bar g \bar \gamma = y_2$.
		However $g \gamma$ is exactly the lift of $\bar g\bar \gamma$ starting at $gy_0$ thus $y_2 = gy$ which proves our claim.
		Point~\ref{enu: lift quasi-convex isometry and stabilizer - stabilizer} follows from the fact that $\nu$ induces an isometry from $Z$ onto $\bar Z$.
		
		\paragraph{} Point~\ref{enu: lift quasi-convex isometry and stabilizer - stabilizer} implies that the projection $\pi : G \rightarrow \bar G$ maps $\stab Y$ onto $\stab{\bar Y}$.
		We prove now that this map is also one-to-one.
		Let $g \in \stab Y$ whose image in $\stab{\bar Y}$ is trivial i.e., $g \in K$.
		The point $gy_0$ belongs to $Y$.
		Moreover, this is a lift of $\bar g\bar y_0 = \bar y_0$.
		Using the isometry between $Y$ and $\bar Y$ we get that $gy_0 = y_0$.
		Since $K$ acts freely on $X \setminus v(\mathcal R)$ and $y_0$ is not a vertex of $X$ we get $g=1$.
	\end{proof}

% !TEX root = notes.tex

\section{Cone over a metric space}
\label{sec: cone}
	
	\subsection{Definition and metric.}
	In this section we fix a number $\rho>0$.
	Its value will be made precise later.
	It should be thought of as a very large parameter.
	
	\begin{defi}
		Let $Y$ be a metric space. 
		The \emph{cone of radius $\rho$ over $Y$}, denoted by $Z(Y)$ is the topological quotient of $Y\times \left[0,\rho\right]$ by the equivalence relation that identifies all the points of the form $(y,0)$.
	\end{defi}
	
	\notas The equivalence class of $(y,0)$, denoted by $v$ is called the \emph{apex} of the cone. 
	By abuse of notation, we still write $(y,r)$ for the equivalence class of $(y,r)$.
	
	\begin{prop}[see {\cite[Chap. I.5, Prop. 5.9]{BriHae99}}]
	\label{res: def distance cone}
		The cone $Z(Y)$ is endowed with a metric characterized in the following way. 
		Let $x=(y,r)$ and $x'=(y',r')$ be two points of $Z(Y)$ then
		\begin{displaymath}
			\cosh \dist x{x'} = \cosh r \cosh r' - \sinh r\sinh r' \cos \angle y{y'},
		\end{displaymath}
		where $\angle y{y'}$ is the \emph{angle at the apex} defined by $\angle y{y'} = \min \left\{ \pi , {\dist y{y'}}/\sinh \rho\right\}$.
		Moreover, if $Y$ is a length space, so is $Z(Y)$.
	\end{prop}
	
\begin{figure}[htbp]
		\centering
		\includegraphics{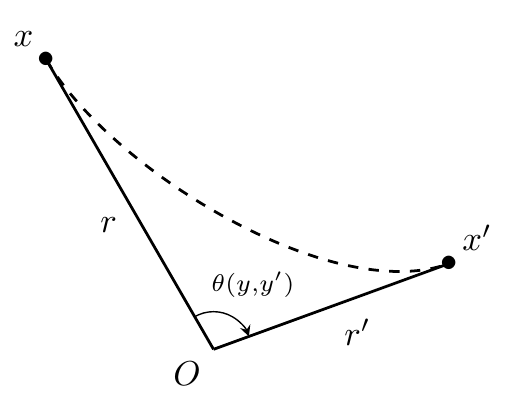}
		\caption{Geometric interpretation of the distance in the cone.}
		\label{fig: distance cone}
	\end{figure}

	In fact, the cone $Z(Y)$ is the ball of radius $\rho$ of the cone $C_{-1}\left(Y/\pi \sinh \rho\right)$ defined in \cite[Chap. I.5]{BriHae99}.
	The distance between two points $x=(y,r)$ and $x'=(y',r')$ of $Z(Y)$ has the following geometric interpretation.
	Consider a geodesic triangle in the hyperbolic plane $\H_2$ such that lengths of two sides are respectively $r$ and $r'$ and the angle between them is $\angle y{y'}$.
	According to the law of cosines, $\dist x{x'}$ is exactly the length of the third side of the triangle (see Figure~\ref{fig: distance cone}).
	
	\begin{prop}[see {\cite[Chap. I.5, Prop. 5.10]{BriHae99}}]
		Let $x=(y,r)$ and $x'=(y',r')$ be two points of $Z(Y)$.
		\begin{itemize}
			\item If $r,r' > 0$ and $\dist y{y'} < \pi \sinh \rho$ then there is a one-to-one correspondence between the geodesics of $Y$ joining $y$ and $y'$ and the geodesics of $Z(Y)$ between $x$ and $x'$.
			\item In all other cases $\dist x{x'} = r + r'$.
			Moreover, there is a unique geodesic between $x$ and $x'$. 
			It goes through the apex of the cone.
		\end{itemize}
	\end{prop}
	
	\exs If $Y$ is a circle of perimeter $2 \pi \sinh \rho$ endowed with the length metric then $Z(Y)$ is the hyperbolic disc of $\H_2$ of radius $\rho$.
	If $Y$ is the real line, then $Z(Y)\setminus\{v\}$ is the universal cover of the punctured hyperbolic disc of radius $\rho$. 
	
	\paragraph{}
	In order to compare the space $Y$ and its cone we introduce the map $\iota : Y \rightarrow Z(Y)$ which sends $y$ to $(y,\rho)$.
It follows from the definition of the metric on $Z(Y)$ that for all $y,y' \in Y$, $\dist[Z(Y)]{\iota(y)}{\iota(y')} = \mu\left( \dist[Y] y{y'}\right)$, where $\mu$ is a map from $\R_+$ into $\R_+$ characterized by
\begin{displaymath}
	\forall t\geq 0, \quad \cosh \mu(t) = \cosh^2 \rho - \sinh^2 \rho \cos \left(\min \left\{ \pi, \frac t{\sinh \rho}\right\}\right).
\end{displaymath}

\begin{prop}
\label{res: map mu}
	The map $\mu$ is continuous, concave, non-decreasing.
	Moreover, we have the followings.
	\begin{enumerate}
		\item \label{enu: mu - comparison map}
		for all $t \geq 0$,
		\begin{math}
		\displaystyle
			t - \frac 1{24}\left(1+\frac 1{\sinh^2\rho}\right) t^3 \leq \mu(t) \leq t.
		\end{math}
		\item \label{emu: mu - lower bound}
		for all $t \in \intval 0{\pi \sinh \rho}$, $t \leq \pi \sinh (\mu(t)/2)$.
	\end{enumerate}
\end{prop}

\begin{figure}[h]
\centering
	\includegraphics{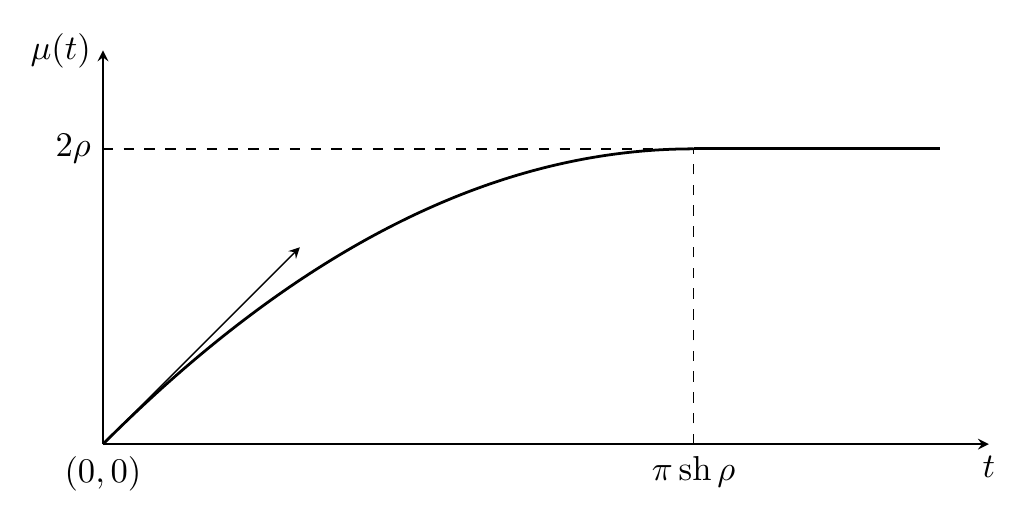}
\caption{Graph of $\mu$.}
\label{fig: graph-mu}
\end{figure}

\begin{proof}
The shape of the graph of $\mu$  is given on Figure~\ref{fig: graph-mu}.
The proof is left to the reader.
\end{proof}

	\subsection{Hyperbolicity of a cone.}
	
	\begin{prop}[Berestovski\u\i\ see {\cite[Chap. II.3, Th. 3.14]{BriHae99}}]
	\label{res: curvature of a cone - limit case}
		Let $Y$ be a metric space.
		If $Y/\pi\sinh \rho$ is $\operatorname{CAT}(1)$ then the cone $Z(Y)$ is $\catold$.
		In particular if $Y$ is a tree then $Z(Y)$ is $\boldsymbol \delta$-hyperbolic.
	\end{prop}
	
	\rem Recall that $\boldsymbol \delta$ is the hyperbolicity constant of the hyperbolic plane $\H_2$.
	
	\begin{prop}
	\label{res: hyperbolicity of a cone}
		Let $Y$ be a metric space.
		The cone $Z(Y)$ is $2\boldsymbol \delta$-hyperbolic.
	\end{prop}
	
	\begin{proof}
	Every triple of points of $Y$ can be isometrically embedded into a tripod.
	Therefore by Proposition~\ref{res: curvature of a cone - limit case} the cone over three points of $Y$ is isometrically embedded into a $\boldsymbol \delta$-hyperbolic space.
	Consequently, for every triple of points in $Z(Y)$, these points together with the apex $v$ of the cone satisfy the four points inequality (\ref{eqn: hyperbolicity condition 1}).
	The proposition follows then from Lemma~\ref{res: criterion for hyperbolicity}.
	\end{proof}
	
	\subsection{Group action on a cone}
	Let $Y$ be a metric space endowed with an action by isometries of a group $H$.
	This action naturally extends to an action by isometries on $Z(Y)$ in the following way.
	For every point $x=(y,r)$ of $Z(Y)$, for every $h\in H$ we put $hx=(hy,r)$.
	Note that $H$ fixes the apex $v$ of the cone.
	Therefore this action is not necessarily proper (even if the action of $H$ on $Y$ is).
	One should think that $H$ acts on $Z(Y)$ as a rotation group with apex $v$.
	
	\begin{lemm}
	\label{res: group acting as a rotation on cone}
		Assume that for every $h \in H$, $\len h \geq \pi \sinh \rho$.
		Then for every point $x \in Z(Y)$, for every $h \in H\setminus\{1\}$, $\dist{hx}x = 2\dist xv$.
	\end{lemm}
	
	\begin{proof}
		We denote by $(y,r)$ the point $x$.
		By assumption $\dist{hy}y \geq \pi \sinh \rho$, hence $\angle y{y'} = \pi$..
		Therefore $\dist{hx}x=2r = 2\dist xv$.
	\end{proof}

	We now assume that the action of $H$ on $Y$ is proper.
	We denote by $\bar Y$ the quotient $Y/H$. 
	For all $y \in Y$, we write $\bar y$ for the image of $y$ in $\bar Y$.
	The space $\bar Y$ is endowed with a metric defined by $\dist{\bar y}{\bar y'} = \inf_{h \in H}\dist y{hy'}$.
	As we previously explained, the action of $H$ on $Z(Y)$ may not be proper.
	Nevertheless the formula $\dist{\bar x}{\bar x'} = \inf_{h \in H}\dist x{hx'}$ still defines a metric on $Z(Y)/H$.
	Moreover, the spaces $Z(Y)/H$ and $Z(Y/H)$ are isometric.

	\begin{lemm}
	\label{res: metric quotient cone}
		Let $l \geq 2\pi \sinh \rho$.
		We assume that for every $h \in H\setminus \{1\}$, $\len h \geq l$.
		Let $x=(y,r)$ and $x'=(y',r')$ be two points of $Z(Y)$.
		If $\dist[Y]y{y'} \leq l - \pi \sinh \rho$ then $\dist{\bar x}{\bar x'} = \dist x{x'}$.
	\end{lemm}

\begin{proof}
	Since $Z(Y/H)$ and $Z(Y)/H$ are isometric, the distance between $\bar x$ and $\bar x'$ in $Z(Y)/H$ is given by 
	\begin{displaymath}
		\cosh \left( \dist{\bar x}{\bar x'}\right) = \cosh r \cosh r' - \sinh r \sinh r' \cos \left( \min \left\{ \pi , \frac {\dist[\bar Y]{\bar y}{\bar y'}}{\sinh \rho}\right\}\right).
	\end{displaymath}
	If $\dist y{y'} <l/2$, then we have $\dist{\bar y}{\bar y'} = \dist y{y'}$.
	It follows that $\dist{\bar x}{\bar x'} = \dist x{x'}$.
	Assume now that $\dist y{y'} \geq l/2$. 
	In particular $\dist y{y'} \geq \pi \sinh \rho$.
	Thus $\dist x{x'} = r + r'$.
	On the other hand, using the triangle inequality, for all $h \in H\setminus\{1\}$, $\dist  y{hy'} \geq l - \dist y{y'}$, thus $\dist {\bar y}{\bar y'} \geq \pi \sinh \rho$.
	Consequently, $\dist{\bar x}{\bar x'} = r+ r' = \dist x{x'}$.
\end{proof}

\section{Cone-off construction}
\label{sec: cone-off construction}

\paragraph{}
The goal of this section is to construct a metric space called \emph{cone-off} obtained by attaching a family $\mathcal Z$ of cones on a base space $X$.
In particular we would like to understand its curvature.
Nevertheless during the exposition we will never use the fact that the spaces we attach are cones.
Therefore we explain the ideas in a more general situation.
In this process the spaces $Z \in \mathcal Z$ are not attached according to an isometry.
Therefore one needs a way to measure the distortion between the glued spaces and the base $X$.
This is the role of the comparison map defined below.

\begin{defi}
\label{def: comparison map}
	Let $a >0$.
	An \emph{$a$-comparison} map is a non-decreasing, concave map $\mu: \R_+ \rightarrow \R_+$ such that for all $t \in  \R_+$, $t-at^3 \leq \mu(t) \leq t$.
\end{defi}

It follows immediately that for all $t \geq 0$, $\mu(t)=0$ if and only if $t=0$.
Moreover, $\mu$ is subadditive, i.e. for all $s,t \in \R_+$, $\mu(s+t) \leq \mu(s)+\mu(t)$.
Hence $\mu$ is 1-Lipschitz.
The map $\mu$ studied in Proposition~\ref{res: map mu} is an $a$-comparison map with $a = (1+1/\sinh^2 \rho)/24$.

\begin{defi}
\label{def: mu X family}
	Let $a >0$ and $\mu$ be an $a$-comparison map.
	Let $X$ be a metric space.
	A \emph{$(\mu,X)$-family} $\mathcal Z$ is a collection of triples $\left(Z,Y,\iota\right)$ where $Z$ is a metric length space and $\iota$ a map from a non-empty subset $Y$ of $X$ into $Z$ such that $\iota(Y)$ is closed in $Z$ and for all $y,y' \in Y$, 
	\begin{equation}
	\label{eqn: def mu X family}
		 \mu\left( \dist[X]y{y'}\right)\leq \dist[Z]{\iota(y)}{\iota(y')}.
	\end{equation}
\end{defi}

\begin{defi}
\label{def:pre-cone off}
	Let $a >0$ and $\mu$ be an $a$-comparison map.
	Let $X$ be a metric space and $\mathcal Z$ a $(\mu, X)$-family.
	The \emph{cone-off over $X$ relatively to $\mathcal Z$} denoted by $\dot X(\mathcal Z)$ (or simply $\dot X$) is obtained by attaching for all $\left(Z,Y,\iota\right) \in \mathcal Z$ the space $Z$ on $X$ along $Y$ according to $\iota$.
\end{defi}

In other words the space $\dot X$ is the quotient of the disjoint union of $X$ and all the $Z$'s by the equivalence relation which identifies every point $y\in Y$ with its image $\iota(y) \in Z$.
To simplify the notations, we use the same letter to design a point of this disjoint union and its equivalence class in $\dot X$.
For the moment $\dot X$ is just a set of points.
Our goal is to define a metric on $\dot X$ and study its properties.

\subsection{Metric on the cone-off}

\paragraph{} We endow the disjoint union of $X$ and all the $Z$'s with the distance induced by $\distV[X]$ and $\distV[Z]$.
This metric is not necessarily finite: the distance between two points in distinct components is infinite.
Let $x$ and $x'$ be two points of $\dot X$.
We define $\dist[SC] x{x'}$ to be the infimum over the distances between two points in the previous disjoint union whose classes in $\dot X$ are respectively $x$ and $x'$.

\rems
\begin{enumerate}
	\item Let $(Z,Y,\iota) \in \mathcal Z$.
	If $x \in Z\setminus Y$ and $x' \notin Z$, then $\dist[SC] x{x'} = + \infty$.
	In particular $\distV[SC]$ is not a distance on $\dot X$ (it does not satisfy the triangle inequality).
	\item Let $x$ and $x'$ be two points of $X$.
	Using the properties of $\mu$ combined with (\ref{eqn: def mu X family}) we get
	\begin{displaymath}
		\mu\left( \dist[X]x{x'}\right) \leq \dist[SC] x{x'} \leq \dist[X]x{x'}.
	\end{displaymath}
\end{enumerate}

\begin{defi}[Chain between two points]
	Let $x$ and $x'$ be two points of $\dot X$.
	A \emph{chain} between $x$ and $x'$ is a finite sequence $C=\left(z_1,\dots,z_m\right)$ such that $z_1=x$ and $z_m=x'$.
	Its \emph{length}, denoted by $l(C)$, is
	\begin{displaymath}
		l(C) = \sum_{j=1}^{m-1} \dist[SC]{z_{j+1}}{z_j}.
	\end{displaymath}
\end{defi}

\begin{prop}[Pseudo-distance, see {\cite[Chap. I.5, Prop. 5.19]{BriHae99}}]
\label{res: pseudo-distance attaching spaces}
	The following map defines a pseudo-distance on $\dot X$,
	\begin{displaymath}
		\begin{array}{ccccl}
			\dot X \times \dot X & \rightarrow & \R_+ && \\
			(x,x') & \rightarrow & \dist[\dot X] x{x'} & = & \inf \left\{ l(C) | C \text{ chain between $x$ and $x'$}\right\}.
		\end{array}
	\end{displaymath}
\end{prop}

\paragraph{} 
By construction, the canonical maps $X \rightarrow \dot X$ and $Z \rightarrow \dot X$ are 1-Lipschitz.
In the remainder of this section we prove that $\distV[\dot X]$ is in fact a distance on $\dot X$.
The next lemma is just a consequence of the triangle inequality.
The proof is left to the reader.

\begin{lemm}
\label{res: chain in X}
	Let $x$ and $x'$ be two points of $\dot X$.
	For all $\eta > 0$ there is a chain $C = \left(z_1,\dots, z_m\right)$ between them such that $l(C) \leq \dist[\dot X] x{x'}+ \eta$ and for all $j \in \intvald 2{m-1}$, $z_j$ belongs to $X$.
\end{lemm}

\begin{comment}
	Let $\eta > 0$.
	By definition there is a chain $C = \left(z_1,\dots, z_m\right)$ between $x$ and $x'$ such that $l(C) \leq \dist[\dot X] x{x'}+ \eta$.
	Assume that there is $j \in \intvald 2{m-1}$ such that $z_j \notin X$.
	There exists $(Z,Y,\iota) \in \mathcal Z$ such that the point $z_j$ belongs to $Z\setminus Y$.
	Since $l(C)$ is finite, $z_{j-1}$ and $z_{j+1}$ also belong to $Z$.
	The triangle inequality gives then
	\begin{eqnarray*}
		\dist[SC] {z_{j+1}}{z_{j-1}} 
		& \leq & \dist[Z]{z_{j+1}}{z_{j-1}} \\
		& \leq & \dist[Z]{z_{j+1}} {z_j} + \dist[Z]{z_j}{z_{j-1}} =  \dist[SC]{z_{j+1}} {z_j} + \dist[SC]{z_j}{z_{j-1}}.
	\end{eqnarray*}
	Consequently, the length of the chain obtained by removing $z_j$ is smaller than $l(C)$.
	After removing every point of $C$ which does not belong to $X$ we obtain the chain of the lemma.
\end{comment}

\begin{lemm}
\label{res: metric on  dot X and Zi coincide}
	Let $(Z,Y,\iota) \in \mathcal Z$.
	Let $x \in Z\setminus Y$.
	Let $d(x,Y)$ be the distance between $x$ and $\iota(Y)$ computed with $\distV[Z]$.
	For all $x' \in \dot X$, if $\dist[\dot X] x{x'} < d(x,Y)$ then $x'$ belongs to $Z$.
	Moreover, $\dist[\dot X] x{x'} = \dist[Z] x{x'}$.
\end{lemm}

\begin{proof}
	Let $\eta > 0$ such that $\dist[\dot X] x{x'} + \eta < d(x,Y)$.
	A chain $C$ between $x$ and $x'$ whose length is bounded above by $\dist[\dot X] x{x'} + \eta$ cannot go outside of $Z$.
	Therefore all its points are contained in $Z\setminus Y$.
	By triangle inequality we obtain $\dist[Z] x{x'} \leq l(C) \leq \dist[\dot X] x{x'} + \eta$.
	This inequality holds for all $\eta > 0$.
	Hence $\dist[Z] x{x'}\leq \dist[\dot X] x{x'}$.
	The other inequality follows from the definition of $\distV[\dot X]$.
\end{proof}

\begin{lemm}
\label{res: metric on X compare to metric on dot X}
	For all $x, x' \in X$, $\mu\left(\dist[X] x{x'}\right) \leq \dist[\dot X] x{x'} \leq \dist[X] x{x'}$.
\end{lemm}

\rem Recall that $\mu$ is 1-Lipschitz and therefore continuous.
The lemma shows in particular that the space $X$ and its image in $\dot X$ have the same topology.

\begin{proof}
	The inequality $\dist[\dot X] x{x'} \leq \dist[X] x{x'}$ follows directly from the definition of $\distV[\dot X]$.
	Let $\eta > 0$.
	 By Lemma~\ref{res: chain in X} there exists a chain $C=\left(z_1,\dots,z_m\right)$ between $x$ and $x'$ whose points belong to $X$ such that $l(C) \leq \dist[\dot X] x{x'} + \eta$.
	 Recall that for all $j \in \intvald 1{m-1}$, $\mu\left(\dist[X]{z_{j+1}}{z_j}\right) \leq \dist[SC]{z_{j+1}}{z_j}$.
	 Using the subadditivity of $\mu$ we get
	 \begin{displaymath}
	 	\mu\left(\dist[X] x{x'}\right)
		\leq \sum_{j=1}^{m-1} \mu\left(\dist[X]{z_{j+1}}{z_j}\right)
		\leq \sum_{j=1}^{m-1} \dist[SC]{z_{j+1}}{z_j}
		= l(C).
	 \end{displaymath}
	Thus for all $\eta > 0$, $\mu\left(\dist[X] x{x'}\right) \leq \dist[\dot X] x{x'} + \eta$, which gives the other inequality.
\end{proof}

\begin{prop}
\label{res: pseudo metric is in fact a metric}
	For all $x,x' \in \dot X$, $\dist[\dot X] x{x'} = 0$ if and only if $x=x'$.
	In particular, $\distV[\dot X]$ is a distance on $\dot X$.
\end{prop}

\begin{proof}
	Suppose that $\dist[\dot X] x{x'}=0$.
	We distinguish two cases.
	\begin{enumerate}
		\item 
		Assume that there exists $(Z,Y,\iota) \in \mathcal Z$ such that $x \in Z\setminus Y$.
		Since $\iota(Y)$ is closed in $Z$ the distance $d(x,Y)$ is positive.
		Lemma~\ref{res: metric on  dot X and Zi coincide} states that $x'$ belongs to $Z$ and $\dist[Z] x{x'} = \dist[\dot X] x{x'}=0$.
		Thus $x=x'$.
		\item
		If $x$ and $x'$ both belong to $X$ then Lemma~\ref{res: metric on X compare to metric on dot X} gives $\mu\left(\dist[X] x{x'}\right) \leq \dist[\dot X] x{x'} = 0$.
		Using the properties of $\mu$, we get $\dist[X] x{x'}=0$. Hence $x=x'$.
	\end{enumerate}
	The  other implication of the proposition is obvious.
\end{proof}

\rem
The proof also tells that the maps $X \rightarrow \dot X$ and $Z \rightarrow \dot X$ are actually embeddings.

\begin{prop}
\label{res: length structure on the cone-off}
	The metric $\distV[\dot X]$ endows $\dot X$ with a length structure.
\end{prop}

\begin{proof}
	By assumption the spaces $X$ and $Z$ are endowed with a length structure.
	Thus $\dot X$ is obtain by attaching together length spaces.
	Since $\distV[\dot X]$ is a metric, $\dot X$ endowed with $\distV[\dot X]$ is a length space (see \cite[Chap. I.5, Lemma 5.20]{BriHae99}).
\end{proof}

\subsection{Uniform approximation of the distance}

To study the curvature of $\dot X$ we should understand how the cone-off construction behaves with respect to ultra-limits.
To that end, we need to approximate the distance between two points of $\dot X$ by a chain such that the number of points involved in this chain only depends on the error and not on the space $X$ or the collection $\mathcal Z$.
More precisely we prove the following result.

\begin{prop}
\label{res: uniform approximation X}
	Let $\epsilon > 0$.
	There exists $M \geq 0$ which only depends on $\epsilon$ and $a$ with the following property.
	Let $x$ and $x'$ be two points of $X$ and $C$ a chain between them whose points all belong to $X$.
	There is a subchain $C'$ of $C$ joining $x$ and $x'$, which does does not contain more than $M(l(C)+1)$ points, such that $l(C') \leq (1+\epsilon) l(C) + \epsilon$.
\end{prop}

\begin{proof}
	We write $C = \left(z_1,\dots,z_m\right)$ for the chain between $x$ and $x'$.
	Let $\eta > 0$.
	We define a subchain of $C$ denoted by $C_\eta = \left(z_{j_1},\dots z_{j_n}\right)$ as follows.
	\begin{itemize}
		\item 
		Put $j_1 = 1$.
		\item 
		Let $k \geq 1$ such that $j_k < m$.
		If $\dist[X]{z_{j_k+1}}{z_{j_k}} > 2\eta$, then we put $j_{k+1} = j_k+1$, otherwise $j_{k+1}$ is the largest index $j \in \intvald{j_k+1}m$ such that $\dist[X]{z_j}{z_{j_k}} \leq 2\eta$.
		\item
		The process stops when $j_k = m$.
	\end{itemize}
	The chain $C_\eta$ joins $x$ and $x'$.
	We denote by $n$ its number of points.
	
	\begin{lemm}
	\label{res: uniform approximation - lemma comparaison length}
		The lengths of the chains $C$ and $C_\eta$ satisfy $l(C_\eta) \leq l(C) + 8an\eta^3$. 
	\end{lemm}
	
	\begin{proof}
		Let $k \in \intvald 1{n-1}$.
		For all $t \in \R_+$, $\mu(t) \geq t - at^3$.
		By Lemma~\ref{res: metric on X compare to metric on dot X}, we get
		\begin{displaymath}
			\sum_{l=j_k}^{j_{k+1}-1}\dist[SC] {z_{l+1}}{z_l}
			\geq \dist[\dot X] {z_{j_{k+1}}}{z_{j_k}}
			\geq \mu\left( \dist[X] {z_{j_{k+1}}}{z_{j_k}} \right)
			\geq \dist[SC] {z_{j_{k+1}}}{z_{j_k}} -  a\dist[X] {z_{j_{k+1}}}{z_{j_k}}^3.
		\end{displaymath}
		If $\dist[X] {z_{j_{k+1}}}{z_{j_k}} \leq 2\eta$, then
		\begin{displaymath}
			\sum_{l=j_k}^{j_{k+1}-1}\dist[SC] {z_{l+1}}{z_l}
			\geq  \dist[SC] {z_{j_{k+1}}}{z_{j_k}} - 8a\eta^3.
		\end{displaymath}		
		If $\dist[X] {z_{j_{k+1}}}{z_{j_k}} > \eta$, then by definition $j_{k+1} = j_k+1$.
		Thus the last inequality holds.
		By summing over $k$ we obtain
		\begin{displaymath}
			l(C)
			= \sum_{l=1}^{m-1}\dist[SC] {z_{l+1}}{z_l}
			\geq \sum_{k=1}^{n-1}\dist[SC] {z_{j_{k+1}}}{z_{j_k}}  - 8an\eta^3
			= l(C_\eta) - 8an\eta^3. \qedhere
		\end{displaymath}
	\end{proof}

	\begin{lemm}
	\label{res: uniform approximation - bounding number of points}
		If $\eta < \sqrt{ 1/10a}$, then the number of points of $C_\eta$ is bounded as follows
		\begin{displaymath}
			n \leq 4\left(\frac{l(C)}\eta +1\right).
		\end{displaymath}
	\end{lemm}
	
	\begin{proof}
		Let $k$ be an integer of $\intvald 1{n-2}$.
		The distances $\dist[X] {z_{j_{k+1}}}{z_{j_k}}$ and  $\dist[X] {z_{j_{k+2}}}{z_{j_{k+1}}} $ cannot be both bounded above by $\eta$.
		Otherwise $j_{k+1}$ would not be the largest index $j \in \intvald{j_k+1}m$ such that $\dist[X]{z_j}{z_{j_k}} \leq 2\eta$.
		Consequently, 
		\begin{displaymath}
			\dist[SC] {z_{j_{k+2}}}{z_{j_{k+1}}} + \dist[SC] {z_{j_{k+1}}}{z_{j_k}} \geq \mu\left( \eta\right) \geq  \eta  - a\eta^3.
		\end{displaymath}
		By summing over $k$ one gets
		\begin{displaymath}
			\left\lfloor \frac{n-1}2 \right\rfloor \left(\eta  - a\eta^3\right) \leq l(C_\eta) \leq l(C) + 8an\eta^3.
		\end{displaymath}
		A small computation leads to 
		\begin{displaymath}
			n \left(1-5a\eta^2\right) \leq 2\frac {l(C)}\eta +3-3 a\eta^2. 
		\end{displaymath}
		Recall that $\eta < \sqrt{ 1/10a}$.
		Consequently, $n \leq 4(l(C)/\eta +1)$.
	\end{proof}
	
	\proofend{End  of the proof of Proposition~\ref{res: uniform approximation X}}
	Combining Lemmas~\ref{res: uniform approximation - lemma comparaison length} and \ref{res: uniform approximation - bounding number of points} yields
	\begin{displaymath}
		l(C_\eta) \leq \left(1+32a\eta^2 \right)l(C) + 32a \eta^3.
	\end{displaymath}
	If one takes $\eta$ small enough then $l(C_\eta) \leq \left(1+\epsilon \right)l(C) + \epsilon$.
	Note that $\eta$ only depends on $a$ and $\epsilon$.
	Moreover, the number of points of $C_\eta$ is bounded above by $4(l(C)/\eta +1)$.
\end{proof}

\begin{coro}
\label{res: uniform approximation - distance}
	Let $\epsilon > 0$. 
	There is a constant $M$ which only depends on $\epsilon$ and $a$ with the following property.
	For all $x,x' \in \dot X$, there exists a chain $C$ between $x$ and $x'$ which does not contain more than $M\left( \dist[\dot X]x{x'} +1\right)$ points and such that $l(C) \leq (1+\epsilon)\dist[\dot X] x{x'} + \epsilon$.
\end{coro}

\begin{proof}
	By Lemma~\ref{res: chain in X}, there exists a chain $C = \left(z_1,\dots, z_m\right)$ between $x$ and $x'$ such that $l(C) \leq \dist[\dot X] x{x'} +\epsilon/2$ and for all $j \in \intvald 2{m-1}$, $z_j \in X$.
	We now apply Proposition~\ref{res: uniform approximation X} to the chain $C_1 = \left(z_2,\dots, z_{m-1}\right)$.
	There exists a constant $M$ which only depends on $\epsilon$ and $a$ and a subchain $C_2$ of $C_1$ satisfying the followings.
	The chains $C_1$ and $C_2$ have the same extremities.
	The number of points of $C_2$ is bounded above by $M(l(C_1)+1)$.
	Moreover, $l(C_2) \leq(1+\epsilon) l(C_1) + \epsilon/2$.
	We extend $C_2$ by adding $z_1$ at the beginning and $z_m$ at the end.
	The number of points of $C'$ is bounded above by $M(\dist[\dot X] x{x'} + \epsilon/2 +1)+2$.
	Its length satisfies 
	\begin{displaymath}
		l(C') \leq (1+\epsilon) l(C) + \epsilon/2 \leq (1+\epsilon)\dist[\dot X] x{x'} + \epsilon. \qedhere
	\end{displaymath}
\end{proof}

\subsection{Ultra-limit and cone-off}
\label{sec: ultra limit attaching spaces}

	\paragraph{} In this section we study the behavior of $\dot X$ under ultra-limits.
	The data that we consider are the followings.
	Let $a > 0$ and $\omega$ be a non-principal ultra-filter.
	For every $n \in \N$ we choose an $a$-comparison map $\mu_n$, a pointed length space $\left(X_n, x_n^0\right)$ and a $(\mu_n,X_n)$-family $\mathcal Z_n$.
	\begin{defi}
		We say that the distortion of the sequence $(\mathcal Z_n)$ is \emph{$\omega$-controlled} if the following holds. 
		For every sequence of triples $\left(Z_n,Y_n,\iota_n\right) \in \prod_{n \in \N}\mathcal Z_n$, for every $(y_n), (y'_n) \in \prod_{n\in \N} Y_n$, 
		\begin{equation}
		\label{eqn: comparison in the limit}
			\limo \mu_n\left(\dist[X_n]{y_n}{y'_n}\right) = \limo \dist[Z_n]{\iota_n(y_n)}{\iota_n(y'_n)}.
		\end{equation}
	\end{defi}

	\paragraph{}
	Since $\mathcal Z_n$ is a $(\mu_n, X_n)$ family we know that the distortion between $X_n$ and the elements of $\mathcal Z_n$ is controlled from below by $\mu_n$ (see Definition~\ref{def: mu X family}).
	This definition says that, at the limit, the distortion is exactly given by $\mu_n$.
	For the remainder of this section, we assume that the distortion of $(\mathcal Z_n)$ is $\omega$-controlled.
	Our goal is to study the space $\limo \dot X_n(\mathcal Z_n)$.
	Before, we define several objects.
	The metric space $\limo \left(X_n, x_n^0\right)$ is denoted by $X$.
	The map $\mu : \R_+ \rightarrow \R_+$ is defined by $\mu(t) = \limo \mu_n(t)$.
	It is also an $a$-comparison-map.	
	
	\paragraph{}
	For every $n \in \N$, for every $\left(Z_n,Y_n,\iota_n\right) \in \mathcal Z_n$ we choose a point $y^0_n \in Y_n$ such that $\dist{x_n^0}{y^0_n}$ is at most $d\left(x_n^0, Y_n\right) +1$ (the distances here are measured with the metric of $X_n$).
	This point exists since $Y_n$ is non-empty.
	Let $\left(Z_n,Y_n,\iota_n\right) \in \prod_{n \in \N}\mathcal Z_n$ be a sequence of triples.
	We define two limit spaces
	\begin{itemize}
		\item $Z = \limo \left( Z_n, \iota_n(y_n^0)\right)$,
		\item $Y = \limo Y_n$ which is a (possibly empty) subset of $X = \limo X_n$.
	\end{itemize}
	It follows from (\ref{eqn: comparison in the limit}) that the map $\iota : Y \rightarrow Z$ given by $\iota(\limo y_n) = \limo \iota_n(y_n)$ is well defined.
	We write then $(Z,Y,\iota) = \limo(Z_n,Y_n,\iota_n)$.
	This triple also satisfies the following property.
	For all $y,y' \in Y$
	\begin{equation}
	\label{eqn: estimation distance limit Z shorter}
		\dist[Z]{\iota(y)}{\iota(y')} = \mu\left(\dist[X] y{y'}\right) \leq \dist[X]y{y'}.
	\end{equation}
	Note that the set $Y = \limo Y_n$ is non-empty if and only if $\left(\dist{x_n^0}{y_n^0}\right)$ is \oeb.
	We write $\prod_\omega\mathcal Z_n$ for the set of sequences of triples satisfying this last condition.
	
	\paragraph{}
	We endow the product $\prod_{n \in \N}\mathcal Z_n$ with the following equivalence relation.
	Given two sequences of triples, $\left(Z_n,Y_n,\iota_n\right)\sim \left(Z'_n,Y'_n,\iota'_n\right)$ if $\left(Z_n,Y_n,\iota_n\right)=\left(Z'_n,Y'_n,\iota'_n\right)$ \oas.
	In particular they define the same limit triple $(Z,Y,\iota) = \limo(Z_n,Y_n,\iota_n)= \limo(Z'_n,Y'_n,\iota'_n)$.
	We can now define $\mathcal Z$ to be the set of triples $(Z,Y,\iota) = \limo (Z_n,Y_n,\iota_n)$ where $ (Z_n,Y_n,\iota_n) \in \prod_\omega\mathcal Z_n/\sim$.
	It follows from (\ref{eqn: estimation distance limit Z shorter}), that $\mathcal Z$ is a $(\mu, X)$-family.
	It allow us to look a the cone-off $\dot X(\mathcal Z)$ over $X$ relatively to $\mathcal Z$.
	
	\paragraph{}
	Our goal is to compare $\limo \dot X_n\left(\mathcal Z_n\right)$ with the metric space $\dot X(\mathcal Z)$.
	To that end we define the following maps (the second kind of maps are defined for every $(Z,Y,\iota) $ in $\mathcal Z$).
	\begin{displaymath}
		\begin{array}{lcccclccc}
			\psi_X:	& X			& \rightarrow	& \limo \dot X_n	& \quad	& \psi_Z:	& Z			& \rightarrow	& \limo \dot X_n \\
					& \limo x_n	& \rightarrow	& \limo x_n		&		& 		& \limo x_n	& \rightarrow	& \limo x_n
		\end{array}
	\end{displaymath}
	Recall that for every $n \in \N$, the embedding $X_n \hookrightarrow \dot X_n$ (\resp $Z_n \hookrightarrow \dot X_n$) is $1$-Lipschitz.
	Thus the maps $\psi$ and $\psi_Z$ are well-defined and 1-Lipschitz.
	Moreover, for all $(Z,Y,\iota) \in \mathcal Z$, for all $y \in Y$, $\psi_Z \circ \iota(y) = \psi_X (y)$.
	Consequently, they induce a map $\dot \psi : \dot X \rightarrow \limo \dot X_n$ whose restriction to $X$ (\resp $Z$) is $\psi_X$ (\resp $\psi_Z$).
	
	\paragraph{}The map $\dot \psi$ cannot be an isometry. 
	Indeed the space $\limo \dot X_n$ is much larger than $\dot X$.
	Imagine for instance that the spaces $Z_n$ that we attach are uniformly bounded. 
	This will be the case later: we will take for $Z_n$ some cones with a fixed radius. 
	One can find a sequence of point $(x_n)$ of $\prod_{n \in \N} X_n$ such that the distance $\dist {x_n^0}{x_n}$ is bounded in $\dot X_n$ but not in $X_n$.
	Therefore this sequence defines a point of $\limo \dot X_n$ which does not correspond to a point of $\dot X$.
	Nevertheless $\dot \psi$ restricted to the neighborhood of $X$ in $\dot X$ induces a local isometry. 
	More precisely we are going to prove the following result.
	
	\begin{prop}
	\label{res: dot psi local isometry}
		Let $x$ be a point of $X$.
		Let $t \in \R_+$.
		The map $\dot \psi$ induces an isometry from $B\left(x, \mu(t)/2\right)$ onto $B(\dot \psi(x), \mu(t)/2)$.
	\end{prop}
	
	The rest of this section is dedicated to the proof of Proposition~\ref{res: dot psi local isometry}.
	We begin by establishing that $\dot \psi$ is 1-lipschitz, and then that it induces a local isometry.
	
	\begin{lemm}
	\label{res: psi dot 1 lipschitz for dist SC}
		If $x$ and $x'$ are two points of $\dot X$ then $\dist[SC] x{x'} \geq \dist {\dot \psi(x)}{\dot \psi(x')}$.
	\end{lemm}
	
	\begin{proof}
		If $\dist[SC] x{x'} = + \infty$ there is nothing to prove.
		Thus we may assume that $x$ and $x'$ both belong either to $X$ or to one of the $Z$ where $(Z,Y \iota) \in \mathcal Z$.
		Suppose that $x,x' \in X$.
		We can write $x = \limo x_n$ and $x'=\limo x'_n$ where $x_n,x'_n \in X_n$.
		By construction of $\dot X_n$ we have for all $n \in \N$, $\dist[X_n] {x_n}{x'_n} \geq \dist[\dot X_n]{x_n}{x'_n}$.
		After taking the $\omega$-limit we get $\dist[X] x{x'} \geq \dist{\psi_X(x)}{\psi_X(x')}$.
		Consider now $(Z,Y,\iota) = \limo (Z_n,Y_n,\iota_n) \in \mathcal Z$.
		If $x$ and $x'$ both belong to $Z$ we prove in the same way that $\dist[Z] x{x'} \geq \dist{\psi_Z(x)}{\psi_Z(x')}$.
		It follows from the definition of $\distV[SC]$ that $\dist[SC] x{x'} \geq \dist {\dot \psi(x)}{\dot \psi(x')}$.
	\end{proof}
	
	\begin{prop}
	\label{res: psi dot 1 lipschitz}
		The map $\dot \psi$ is 1-lipschitz.
	\end{prop}
	
	\begin{proof}
		Let $x$ and $x'$ be two points of $\dot X$.
		Let $C = \left(z_0, \dots ,z_m\right)$ be a chain of points joining $x$ and $x'$.
		By Lemma~\ref{res: psi dot 1 lipschitz for dist SC}, we have
		\begin{displaymath}
			\dist{\dot \psi(x)}{\dot \psi (x')}
			\leq \sum_{j = 0}^{m-1} \dist{\dot\psi(z_{j+1})}{\dot \psi (z_j)}
			\leq \sum_{j = 0}^{m-1} \dist[SC]{z_{j+1}}{z_j}
			= l(C).
		\end{displaymath}
		This inequality holds for all chains joining $x$ and $x'$, thus $\dist{\dot \psi(x)}{\dot \psi (x')} \leq \dist[\dot X] x{x'}$.
	\end{proof}
	
	\begin{lemm}
	\label{res: compare distance ultra-limit tending to zero}
		Let $(x_n)$ and $(x'_n)$ be two sequences of $\prod_\omega \dot X_n$ such that $\limo \dist[\dot X_n]{x_n}{x'_n}=0$.
		\begin{enumerate}
			\item Assume that $x_n, x'_n \in X_n$ \oas, then $\limo \dist[X_n] {x_n}{x'_n} =0$.
			\item Assume that there exists a sequence of triples $(Z_n,Y_n,\iota_n) \in \prod_{n\in \N}\mathcal Z_n$ such that $x_n, x'_n \in Z_n$ \oas, then $\limo \dist[Z_n] {x_n}{x'_n} =0$.
		\end{enumerate}
	\end{lemm}
	
	\rem In this lemma the sequences $(x_n)$ and $(x'_n)$ do not necessarily define points of $X$ or $Z \in \mathcal Z$.
	In particular we do not assume that $(Z_n,Y_n,\iota_n)$ belongs to $\prod_\omega \mathcal Z_n$.
	In other words $Y= \limo Y_n$ may be empty.
	
	\begin{proof}
		Assume that $x_n,x'_n$ belong to $X_n$ \oas.
		By Lemma~\ref{res: metric on X compare to metric on dot X},
		\begin{displaymath}
			\dist[X_n] {x_n}{x'_n} - a \dist[X_n] {x_n}{x'_n}^3 
			\leq \mu_n\left(\dist[X_n] {x_n}{x'_n} \right) 
			\leq \dist[\dot X_n] {x_n}{x'_n} \text{\oas}.
		\end{displaymath}
		In particular, $\limo \dist[X_n]{x_n}{x'_n} =0$.
		
		\paragraph{} Assume now that there is $(Z_n,Y_n,\iota_n) \in \prod_{n\in \N}\mathcal Z_n$ such that $x_n, x'_n \in Z_n$ \oas.
		By Lemma~\ref{res: chain in X}, for every $n \in \N$, there is a chain $C_n$ between $x_n$ and $x'_n$ such that $\limo l(C_n) =0$ and every point of $C_n$ distinct from $x_n$ and $x'_n$ belongs to $X_n$. 
		By doubling if necessary the first and the last point of $C_n$ we can assume the following.
		The second point $y_n$ and the last point but one $y'_n$ of $C_n$ belong to $Y_n$.
		Moreover, $\dist[SC]{x_n}{y_n} = \dist[Z_n]{x_n}{y_n}$ and $\dist[SC]{x'_n}{y'_n} = \dist[Z_n]{x'_n}{y'_n}$.
		It follows that 
		\begin{displaymath}
			l(C_n)  \geq \dist[Z_n] {x_n}{y_n} + \dist[\dot X_n] {y_n}{y'_n} + \dist[Z_n] {y'_n}{x'_n}.
		\end{displaymath}
		In particular the following quantities $\dist[Z_n] {x_n}{y_n}$, $\dist[\dot X_n] {y_n}{y'_n}$ and $\dist[Z_n] {y'_n}{x'_n}$ converge to zero.
		The triangle inequality leads to $\limo \dist[Z_n] {x_n}{x'_n} =   \limo \dist[Z_n] {y_n}{y'_n}$.
		However $y_n,y'_n$ belong to $Y_n$ \oas.
		It follows from the previous point that $\limo \dist[X_n] {y_n}{y'_n} =0$.
		Combined with (\ref{eqn: comparison in the limit}) it gives
		\begin{displaymath}
			\limo \dist[Z_n] {x_n}{x'_n} = \limo \dist[Z_n] {y_n}{y'_n}  = \limo \mu_n\left(\dist[X_n]{y_n}{y'_n} \right) = \mu\left(\limo\dist[X_n]{y_n}{y'_n}\right)=0. \qedhere
		\end{displaymath}
	\end{proof}
	
	\begin{lemm}
	\label{res: compare distSC ultra-limit base}
		Let $x$ and $x'$ be two points of $X$. 
		Let $(z_n)$ and $(z'_n)$ be two sequences which belong to $\prod_\omega \dot X_n$ such that $\psi_X (x) = \limo z_n$ and $\psi_X (x') = \limo z'_n$.
		Then $\limo \dist[SC] {z_n}{z'_n} \geq \dist[SC] x{x'}$.
	\end{lemm}
	
	\begin{proof}
		We can write that $x= \limo x_n$ and $x' = \limo x'_n$ where $x_n, x'_n \in X_n$.
		By definition $\psi_X(x)$ and $\psi_X(x')$ are respectively $\limo x_n$ and $\limo x'_n$ but seen as points of $\limo \dot X_n$, thus $\limo \dist [\dot X_n] {x_n}{z_n} = 0$ and $\limo \dist [\dot X_n] {x'_n}{z'_n} = 0$.
		We distinguish two cases.
		
		\paragraph{First case.}
		Assume that $z_n$ and $z'_n$ both belong to $X_n$ \oas.
		By Lemma~\ref{res: compare distance ultra-limit tending to zero}, $\limo \dist[X_n]{x_n}{z_n} =0$ and $\limo \dist[X_n]{x'_n}{z'_n} =0$.
		Using the triangle inequality we get
		\begin{displaymath}
			\limo \dist[X_n]{z_n}{z'_n} = \limo \dist[X_n]{x_n}{x'_n} = \dist[X] x{x'} \geq \dist[SC] x{x'}.
		\end{displaymath}
		
		\paragraph{Second case.} Assume that three exists $(Z_n,Y_n,\iota_n) \in \prod_{n\in \N}\mathcal Z_n$ such that $z_n, z'_n \in Z_n$ \oas.
		We write $(Z,Y,\iota)$ for $\limo (Z_n,Y_n,\iota_n)$.
		It follows from Lemma~\ref{res: compare distance ultra-limit tending to zero} and the definition of the metric on $\dot X_n$ that  there exists a sequence $(y_n)$ such that $y_n \in Y_n$ \oas and 
		\begin{displaymath}
			\limo\left( \dist[X_n]{x_n}{y_n} + \dist[Z_n] {y_n}{z_n}\right) = 0.
		\end{displaymath}
		Hence $\limo y_n$ defines a point of $X$ which equals $x$.
		In particular $Y$ is non-empty, which means that the sequence $(Z_n,Y_n,\iota_n)$ belongs in fact to $\prod_\omega\mathcal Z_n$.
		We construct an analogue sequence $(y'_n)$ for $x'$.
		By triangle inequality 
		\begin{displaymath}
			 \limo \dist[Z_n]{z_n}{z'_n} =  \limo \dist[Z_n]{y_n}{y'_n} =\dist[Z] {\iota(x)}{\iota(x')} \geq \dist [SC] x{x'}.
		\end{displaymath}
		
		\paragraph{} Note that the $\omega$-limit and the infimum that defines $\distV[SC]$ can be swapped.
		Consequently, $\limo \dist[SC]{z_n}{z'_n} \geq \dist[SC] x{x'}$.
	\end{proof}

	\begin{lemm}
	\label{res: compare distSC ultra-limit cone}
		Let $(Z,Y,\iota)= \limo(Z_n,Y_n,\iota_n) \in \mathcal Z$.
		Let $x$ and $x'$ be two points of $Z$. 
		Let $(z_n)$ and $(z'_n)$ be two sequences which belong to $\prod_\omega \dot X_n$ such that $\psi_Z (x) = \limo z_n$ and $\psi_Z (x') = \limo z'_n$.
		Then $\limo \dist[SC] {z_n}{z'_n} \geq \dist[SC] x{x'}$.
	\end{lemm}

	\begin{proof}
		The proof is essentially the same as for the previous lemma.
		By (\ref{eqn: estimation distance limit Z shorter}) $\dist[SC] x{x'} = \dist[Z] x{x'}$.
		We can write that $x= \limo x_n$ and $x' = \limo x'_n$ where $x_n, x'_n \in Z_n$.
		By definition $\psi_Z (x)$ and $\psi_Z(x')$ are respectively $\limo x_n$ and $\limo x'_n$ but seen as points of $\limo \dot X_n$, thus $\limo \dist [\dot X_n] {x_n}{z_n} = 0$ and $\limo \dist [\dot X_n] {x'_n}{z'_n} = 0$.
		We distinguish three cases.
		
		\paragraph{First case.}
		Assume that $z_n$ and $z'_n$ both belong to $X_n$ \oas.
		It follows from Lemma~\ref{res: compare distance ultra-limit tending to zero} and the definition of the metric on $\dot X_n$ that  there exists a sequence $(y_n)$ such that $y_n \in Y_n$ \oas and 
		\begin{displaymath}
			\limo\left( \dist[Z_n]{x_n}{y_n} + \dist[X_n] {y_n}{z_n}\right) = 0.
		\end{displaymath}
		We construct an analogue sequence $(y'_n)$ for $x'$.
		The triangle inequality gives
		\begin{displaymath}
			\limo \dist[X_n]{z_n}{z'_n} 
			= \limo \dist[X_n]{y_n}{y'_n} 
			\geq \mu\left(\limo \dist[X_n]{y_n}{y'_n}  \right).
		\end{displaymath}
		It follows from (\ref{eqn: comparison in the limit}) that
		\begin{displaymath}
			\limo \dist[X_n]{z_n}{z'_n} 
			\geq \limo \dist[Z_n]{y_n}{y'_n}
			= \limo \dist[Z_n]{x_n}{x'_n} = \dist[SC] x{x'}.
		\end{displaymath}

		\paragraph{Second case.} Assume that $z_n,z'_n \in Z_n$ \oas.
		It follows from Lemma~\ref{res: compare distance ultra-limit tending to zero} that the limits $\limo \dist[Z_n]{x_n}{z_n}$ and $\limo \dist[Z_n]{x'_n}{z'_n}$ equal zero.
		By the triangle inequality, $\limo \dist[Z_n] {z_n}{z'_n} = \limo \dist[Z_n]{x_n}{x'_n} = \dist[SC] x{x'}$.
						
		\paragraph{Third case.} Assume that there exists $(Z'_n,Y'_n,\iota'_n) \in \prod_{n \in \N}\mathcal Z_n$ whose limit $(Z',Y',\iota')$ is distinct from $(Z,Y,\iota)$ and such that $z_n,z'_n \in Z'_n$ \oas.
		It follows from Lemma~\ref{res: compare distance ultra-limit tending to zero} and the definition of the metric on $\dot X_n$ that  there exist two sequences $(y_n)$ and $(t_n)$ such that $y_n \in Y_n$, $t_n \in Y'_n$ \oas and 
		\begin{displaymath}
			\limo\left( \dist[Z_n]{x_n}{y_n} + \dist[X_n] {y_n}{t_n} + \dist[Z'_n]{t_n}{z_n}\right) = 0.
		\end{displaymath}
		In the same way we define two sequences $(y'_n)$ and $(t'_n)$ for $x'$.
		Using the triangle inequality and (\ref{eqn: comparison in the limit}) we have
		\begin{displaymath}
			\begin{array}{lclcl}
				\limo \dist[Z'_n]{z_n}{z'_n} & = & \limo \dist[Z'_n]{t_n}{t'_n} &=& \limo \mu_n \left(\dist[X_n]{t_n}{t'_n}\right), \\
				\limo \dist[Z_n]{x_n}{x'_n} & = & \limo \dist[Z_n]{y_n}{y'_n} &=& \limo \mu_n \left(\dist[X_n]{y_n}{y'_n}\right).
			\end{array}
		\end{displaymath}
		Nevertheless $\mu_n$ is 1-Lipschitz thus
		\begin{eqnarray*}
			\dist{ \mu_n \left(\dist[X_n]{t_n}{t'_n}\right)}{ \mu_n \left(\dist[X_n]{y_n}{y'_n}\right)}
			& \leq & \dist{\dist[X_n]{t_n}{t'_n}}{\dist[X_n]{y_n}{y'_n}} \\
			& \leq & \dist[X_n] {y_n}{t_n} +\dist[X_n] {y'_n}{t'_n}
		\end{eqnarray*}
		Consequently, $\limo \dist[Z'_n] {z_n}{z'_n} = \limo \dist[Z_n]{x_n}{x'_n} = \dist[SC] x{x'}$.
				
		\paragraph{} Note that the $\omega$-limit and the infimum that defines $\distV[SC]$ can be swapped.
		It follows that  $\limo \dist[SC]{z_n}{z'_n} \geq \dist[SC] x{x'}$.
	\end{proof}
	
	\begin{lemm}
	\label{res: ultra limit compare length chains}
		Let $C = \left(z^0, \dots, z^m\right)$ be a chain between two points of $\dot X$.
		For all $j \in \intvald 0m$ we consider a sequence $(z_n^j) \in \prod_\omega \dot X_n$ such that $\dot \psi(z^j) = \limo z_n^j$.
		For all $n \in \N$, we define a chain of $\dot X_n$ as follows: $C_n = (z_n^0, \dots, z_n^m)$.
		Then $l(C) \leq \limo l(C_n)$.
	\end{lemm}

	\begin{proof}
		It follows directly from Lemmas~\ref{res: compare distSC ultra-limit base} and \ref{res: compare distSC ultra-limit cone}.
	\end{proof}
	
	\begin{lemm}
	\label{res: dot psi pre surjectivity}
		Let $x = \limo x_n$ be a point of $X$.
		Let $t \in \R_+$.
		Let $z$ be a point of $B(\dot \psi(x), \mu(t))$.
		There exists $x' \in \dot X$ such that $\dot \psi(x') =z$.
	\end{lemm}
	
	\begin{proof}
		Recall that $B(\dot \psi(x), \mu(t))$ is a subset of $\limo \dot X_n$.
		Thus we can write $z = \limo z_n$, where $z_n \in \dot X_n$.
		Since $z$ belongs to $B(\dot \psi(x), \mu(t))$, $\limo \dist[\dot X_n]{x_n}{z_n} < \mu(t)$.
		We distinguish two cases.
		
		\paragraph{First case.} Assume that $z_n$ belongs to $X_n$ \oas.
		By Lemma~\ref{res: metric on X compare to metric on dot X}
		\begin{displaymath}
			\limo \mu_n \left(\dist[X_n]{x_n}{z_n}\right) \leq \limo \dist[\dot X_n] {x_n}{z_n} <\mu(t) = \limo \mu_n(t).
		\end{displaymath}
		Since $\mu_n$ is non-decreasing, $\dist[X_n]{x_n}{z_n} < t$, \oas.
		Thus $\limo z_n$ defines a point of $X$, whose image by $\dot \psi$ is $z$.
		
		\paragraph{Second case.} Assume that there is $(Z_n,Y_n,\iota_n) \in \prod_{n\in \N}\mathcal Z_n$ such that $z_n \in Z_n$ \oas.
		By the definition of the distance on $\dot X_n$ there is a sequence $(y_n)$ such that $y_n \in Y_n$ \oas and
		\begin{displaymath}
			\limo\left(\dist[\dot X_n]{x_n}{y_n} + \dist[Z_n]{y_n}{z_n}\right) \leq \limo \dist[\dot X_n] {x_n}{z_n} < \mu(t) = \limo \mu_n(t).
		\end{displaymath}
		As above we prove that $y = \limo y_n$ is a well defined point of $X$.
		Moreover, since the sequence $\left(\dist[X_n]{x_n}{y_n}\right)$ is \oeb, $(Z,Y,\iota) = \limo (Z_n,Y_n,\iota_n)$ is in fact an element of $\mathcal Z$.
		It follows that $\limo z_n$ defines a point of $Z$ whose image by $\psi_Z$ is $z$.
	\end{proof}
	
	\begin{lemm}
	\label{res: dot psi pre isometry}
		Let $x = \limo x_n$ be a point of $X$.
		Let $t \in \R_+$.
		Let $y,y'$ be two points of $\dot X$ such that $\dot \psi(y), \dot \psi(y')$ belong to $B(\dot \psi(x),\mu(t)/2)$.
		Then $\dist{\dot \psi(y)}{\dot \psi(y')} = \dist y{y'}$.
	\end{lemm}
	
	\begin{proof}
		By assumption $\dist{\dot \psi(y)}{\dot \psi(y')} \leq \dist{\dot \psi(y)}{\dot \psi(x)} + \dist{\dot \psi(x)}{\dot \psi(y')} < \mu(t)$.
		Let $\eta > 0$ such that $\dist{\dot \psi(y)}{\dot \psi(y')} + \eta  < \mu(t)$.
		There exist two sequences $(z_n)$ and $(z'_n)$ of $\prod_\omega \dot X_n$ such that $\dot \psi(y) = \limo z_n$ and $\dot \psi(x') = \limo z'_n$.
		It follows that $\dist[\dot X_n]{z_n}{z'_n} +\eta < \mu(t)$ \oas.
		According to Proposition~\ref{res: uniform approximation - distance} there exists an integer $m$ and, for all $n \in \N$, a chain $C_n = (z_n^0,\dots z_n^m)$ of $\dot X_n$ between $z_n$ and $z'_n$ such that $l(C_n) \leq \dist[\dot X_n]{z_n}{z'_n} + \eta$.
		It is worth pointing out that $m$, the number of points of $C_n$, does not depend on $n$.
		Note also that for all $j \in \intvald 0m$ the distance between $z_n^j$ and one of the points $z_n$ and $z'_n$ is less than $\mu(t)/2$.
		Thus $z^j = \limo z_n^j$ is a well defined point of $\limo \dot X_n$.
		Moreover, its distance to either $\dot \psi(y)$ or $\dot \psi(y')$ is less than $\mu(t)/2$.
		Consequently, $z^j$ belongs to $B(\dot \psi(x), \mu(t))$.
		By Lemma~\ref{res: dot psi pre surjectivity} there is a point $y^j \in \dot X$ such that $\dot \psi(y^j)=z^j$.
		We choose $y^0=y$ and $y^m=y'$.
		Hence $C = (y^0,\dots y^m)$ is a chain of $\dot X$ between $y$ and $y'$.
		According to Lemma~\ref{res: ultra limit compare length chains} its length satisfies
		\begin{displaymath}
			l(C) \leq \limo l(C_n) \leq \limo \dist[\dot X_n]{z_n}{z'_n} + \eta \leq \dist{\dot \psi(y)}{\dot \psi(y')} +\eta.
		\end{displaymath}
		Therefore $\dist[\dot X] y{y'} \leq  \dist{\dot \psi(y)}{\dot \psi(y')} +\eta$.
		This inequality holds for all $\eta >0$.
		Thus $\dist[\dot X] y{y'} \leq  \dist{\dot \psi(y)}{\dot \psi(y')}$.
		Recall that $\dot \psi$ is 1-lipschitz (Lemma~\ref{res: psi dot 1 lipschitz}).
		This gives $\dist[\dot X] y{y'} =  \dist{\dot \psi(y)}{\dot \psi(y')}$.
	\end{proof}
	
	\begin{proof}[Proof of Proposition~\ref{res: dot psi local isometry}]
		Let $y$ and $y'$ be two points of $B\left(x, \mu(t)/2\right)$.
		Since $\dot \psi$ is 1-lipschitz, $\dot \psi(y)$ and $\dot \psi (y')$ belong to the ball $B(\dot \psi(x), \mu(t)/2)$.
		By Lemma~\ref{res: dot psi pre isometry}, $\dist{\dot \psi(y)}{\dot \psi(y')} = \dist y{y'}$.
		Thus $\dot \psi$ preserves the distances.
		It remains to prove that $\dot \psi$ is onto.
		Let $z$ be a point of $B(\dot \psi(x), \mu(t)/2)$.
		According to Lemma~\ref{res: dot psi pre surjectivity} there is $y \in \dot X$ such that $\dot \psi (y) =z$.
		We should prove show that $y$ belong to $B(x, \mu(t)/2)$.
		By construction $\dot \psi (x)$ and $\dot \psi(y)$ are two points of $B(\dot \psi(x), \mu(t)/2)$.
		It follows from Lemma~\ref{res: dot psi pre isometry} that $\dist xy = \dist{\dot \psi (x)}{\dot \psi(y)} < \mu(t)/2$.
	\end{proof}

\subsection{Hyperbolicity of the cone-off}

	\paragraph{}
	Let $a >0$ and $\mu$ be an $a$-comparison map.
	In section we study the curvature of the space $\dot X(\mathcal Z)$ when the base $X$ is a hyperbolic length space and $\mathcal Z$ a $(\mu,X)$-family such that for every $(Z,Y,\iota) \in \mathcal Z$, $Z$ is hyperbolic.
	We start with the special case of tree graded space.
	
	\begin{defi}
		Let $X$ be a metric length space.
		Let $\mathcal P$ be a collection of closed path-connected subsets (called \emph{pieces}).
		We say that $X$ is \emph{tree graded} with respect to $\mathcal P$ if the following holds
		\begin{enumerate}
			\item Every two different pieces have at most one common point.
			\item Every simple loop is contained in contained in one piece. 
		\end{enumerate}
	\end{defi}
	
	\rem The pieces in this definition have no relation the the pieces of the usual small cancellation theory!
	For more details about tree-graded spaces we refer the reader to {Drutu:2005tz}.

	\begin{prop}
	\label{res: curvature of dot X - limit case}
		Let $\delta >0$.
		Let $a >0$ and $\mu$ be an $a$-comparison map.
		Let $X$ be a tree-graded space with respect to $\mathcal P$.
		Let $\mathcal Z$  a $(\mu,X)$-family.
		We assume that there is a one-to-one correspondence between $\mathcal P$ and the the collection of subsets $Y$ belonging to a triple $(Z,Y,\iota)$ of $\mathcal Z$.
		Then $\dot X(\mathcal Z)$ is $\delta$-hyperbolic as well.
	\end{prop}
	
	\rem In particular if $X$ is an $\R$-tree and the $Y$'s are subtrees sharing at most one point then the proposition applies.
	
	\begin{proof}
		By construction the space $\dot X(\mathcal Z)$ is a tree-graded space with respect to the collection $\set{Z}{(Z,Y,\iota) \in \mathcal Z}$.
		Every piece in this new structure is $\delta$-hyperbolic. 
		However, if one glues together two $\delta$-hyperbolic spaces sharing exactly one point one gets a $\delta$-hyperbolic space.
		Therefore $\dot X(\mathcal Z)$ is $\delta$-hyperbolic.
	\end{proof}
	
	The next result is a small ``perturbation'' of the previous one.
	In particular $X$ is no more a tree but a $\delta_0$-hyperbolic space.
	Moreover, we assume that the $Y$'s are $2 \delta_0$-quasi-convex and allow them to have a small overlap.
	To measure that overlap we introduce the parameter $\Delta$.

	\begin{displaymath}
		\Delta(\mathcal Z) = \sup \set{\diam \left(Y_1^{+ 5\delta_0} \cap Y_2^{+ 5\delta_0}\right) }{(Z_1,Y_1,\iota_1)\neq (Z_2,Y_2,\iota_2) \in \mathcal Z}.
	\end{displaymath}
	
	\begin{prop}
	\label{res: curvature of dot X - general case}
		Let $\delta \geq 0$.
		Let $a, \eta, t>0$.
		There exist positive constants $\delta_0 = \delta_0(\delta,\eta,a,t)$ and $\Delta_0= \Delta_0(\delta,\eta,a,t)$ satisfying the following property.
		Let $\mu$ be an $a$-comparison map.
		Let $X$ be a $\delta_0$-hyperbolic length space and $\mathcal Z$ be a $(\mu,X)$-family such that for every $(Z,Y,\iota) \in \mathcal Z$, $Z$ is $\delta$-hyperbolic length space, $Y$ is a $2\delta_0$-quasi-convex subset of $X$ and for all $y,y'\in Y$
		\begin{displaymath}
			\mu\left(\dist[X] y{y'}\right) \leq \dist[Z]{\iota(y)}{\iota(y')} \leq \mu\left(\dist[X]y{y'}\right) + 8\delta_0.
		\end{displaymath}
		If $\Delta(\mathcal Z) \leq \Delta_0$ then every ball of radius $\mu(t)/8$ of $\dot X(\mathcal Z)$ is $(\delta + \eta)$-hyperbolic.
	\end{prop}
	
	\begin{proof}
		We argue by contradiction.
		Assume that the proposition is false.
		For every $n \in \N$ one can find
		\begin{enumerate}
			\item an $a$-comparison map $\mu_n : \R_+ \rightarrow \R_+$,
			\item a geodesic, $\delta_n$-hyperbolic length space $X_n$ where $\delta_n = o(1)$,
			\item a $(\mu_n,X_n)$-family $\mathcal Z_n$ such that for every $(Z,Y,\iota)$ in  $\mathcal Z_n$, the space $Z$ is a $\delta$-hyperbolic length space, $Y_n$ is a $2\delta_n$-quasi-convex subset of $X_n$ and for all $y,y'\in Y_n$
		\begin{equation}
		\label{eqn: curvature of dot X - general case}
			\mu\left(\dist[X_n] y{y'}\right) \leq \dist[Z_n]{\iota(y)}{\iota(y')} \leq \mu\left(\dist[X_n]y{y'}\right) + 8\delta_n.
		\end{equation}
			Moreover, $\Delta(\mathcal Z_n) = o(1)$,
			\item a point $x_n$ in $\dot X_n(\mathcal Z_n)$ such that the ball $B_n=B\left( x_n,\mu_n(t)/8\right)$ is not $(\delta + \eta)$-hyperbolic.
		\end{enumerate}		
		
		\paragraph{}
		Let us fix a non-principal ultra-filter $\omega$.
		First note that $d(x_n,X_n) < 3\mu_n(t)/8$ \oas.
		Otherwise there is a sequence of triples $(Z_n,Y_n,\iota_n) \in \prod_{n\in \N} \mathcal Z_n$ such that $x_n \in Z_n$ \oas.
		Moreover, by Lemma~\ref{res: metric on  dot X and Zi coincide}, the ball $B_n$ is contained in $Z_n$ and the metrics $\distV[\dot X_n]$ and $\distV[Z_n]$ coincide on $B_n$.
		Since $Z_n$ is $\delta$-hyperbolic, so is $B_n$, a contradiction.
		We denote by $x_n^0$ a point of $X_n$ such that $\dist[\dot X_n]{x_n}{x_n^0} \leq 3\mu_n(t)/8$.
		In particular $B_n$ is contained in the ball $B\left(x_n^0, \mu_n(t)/2\right)$ of $\dot X_n$.
		
		\paragraph{}
		Since $X_n$ is a $\delta_n$-hyperbolic length space, the space $X = \limo \left(X_n,x_n^0\right)$ is an $\R$-tree.
		We denote by $x^0$ the point $x^0 = \limo x_n^0$.
		We define an $a$-comparison map $\mu : \R_+ \rightarrow \R_+$ given by $\mu(s) = \limo \mu_n(s)$ for all $s \in \R_+$.
		It follows from our assumption~(\ref{eqn: curvature of dot X - general case}) that the distortion of $(\mathcal Z_n)$ is $\omega$-controlled.
		As explained in Section~\ref{sec: ultra limit attaching spaces}, we construct a $(\mu,X)$-family $\mathcal Z$ and a map 
		\begin{displaymath}
			\dot \psi : \dot X(\mathcal Z) \rightarrow \limo \left(\dot X_n(\mathcal Z_n), x_n^0\right).	
		\end{displaymath}
		According to Proposition~\ref{res: dot psi local isometry}, $\dot \psi$ induces an isometry from $B\left(x^0, \mu(t)/2\right)$ onto $B(\dot\psi(x^0), \mu(t)/2)$ which is exactly $\limo B(x_n^0,\mu_n(t)/2)$.
		
		\paragraph{}
		Let $(Z,Y,\iota)= \limo(Z_n,Y_n,\iota_n)$ be an element of $\mathcal Z$.
		By assumption for every $n \in \N$, $Z_n$ is $\delta$-hyperbolic, thus so is $Z=\limo Z_n$.
		On the other hand, for every $n \in \N$, $Y_n$ is a $2\delta_n$-quasi-convex subset of $X_n$.
		Thus $Y=\limo Y_n$ is a subtree of $X$.
		Consider now an other triple  $(Z',Y',\iota')= \limo(Z'_n,Y'_n,\iota'_n)$ distinct from $(Z,Y,\iota)$.
		In particular $(Z_n,Y_n,\iota_n) \neq (Z'_n,Y'_n,\iota'_n)$, \oas.
		We assumed that $\Delta(\mathcal Z_n)$ tends to zero as $n$ approaches infinity.
		By Proposition~\ref{res: ultra-limit and diaminter} $\diam(Y\cap Y') = 0$.
		Thus $Y$ and $Y'$ share at most one point.
		Consequently, $\mu$, $X$ and $\mathcal Z$ satisfy the assumptions of Proposition~\ref{res: curvature of dot X - limit case}.
		Hence $\dot X(\mathcal Z)$ is $\delta$-hyperbolic.
		It follows that $\limo B_n$ is also $\delta$-hyperbolic.
		By Proposition~\ref{res: ultra limit of spaces which is hyp}, $B_n$ is $(\delta+\eta)$-hyperbolic \oas, which contradicts our assumptions.
	\end{proof}

% !TEX root = notes.tex

\section{Small cancellation theory}
\label{sec: small cancellation theory}

\subsection{General framework}

\paragraph{}
In this section $X$ is a proper $\delta$-hyperbolic geodesic space endowed with a proper and co-compact action by isometries of a group $G$.
We consider a family $\mathcal Q$ of pairs $(H,Y)$ where $Y$ is a strongly quasi-convex subset of $X$ and $H$ a subgroup of $\stab Y$ acting co-compactly on $Y$.
We assume that $G$ acts on $\mathcal Q$ and that $\mathcal Q /G$ is finite.
The action of $G$ on $\mathcal Q$ is defined as follows.
For every $(H,Y) \in \mathcal Q$, for every $g \in G,$ $g(H,Y) = (gHg^{-1}, gY)$.
In some applications the spaces $Y$'s might not be strongly quasi-convex but simply uniformly quasi-convex (i.e., there exits $\alpha$ such that all the $Y$'s are $\alpha$-quasi-convex).
Then we can substitute $Y$ for an appropriate neighborhood of $Y$ (see Lemma~\ref{res: strong neighborhood of a quasi-convex}) to satisfy the quasi-convexity assumption.
In Section~\ref{sec: gromov monster} we give an example of such an operation.

\paragraph{}
We denote by $K$ the (normal) subgroup generated by the subgroups $H$ where $(H,Y) \in \mathcal Q$.
Our goal is to understand the quotient $\bar G = G/K$.
To that end we consider two parameters which respectively play the role of the length of the largest piece and the length of the smallest relation in the usual small cancellation theory.
\begin{eqnarray*}
	\Delta (\mathcal Q) & = & \sup \set{\diam \left(Y_1^{+5 \delta} \cap Y_2^{+ 5 \delta} \right)}{ (H_1,Y_1) \neq (H_2,Y_2) \in \mathcal Q} \\
	T(\mathcal Q) & = & \inf \set{\len h}{h \in H\setminus\{1\},\; (H,Y) \in \mathcal Q}
\end{eqnarray*}
Let us recall the general strategy sketched in the introduction.
First we construct a space $\dot X$ by attaching on $X$ cones of bases $Y$, where $(H,Y) \in \mathcal Q$.
Under a small cancellation assumption on $\Delta(\mathcal Q)$ and $T(\mathcal Q)$ it turns out that $\dot X$ is hyperbolic.
Moreover the groups $H$ where $(H,Y) \in \mathcal Q$ define a rotation family.
This allow us to apply the results of Section~\ref{sec: rotation family}.
In particular we will see that $\bar G$ acts by isometries on the space $\bar X = \dot X/K$ which is also hyperbolic.
In the next sections we will study the properties of the action of $\bar G$ on $\bar X$.
In particular we provide estimates for the invariants $A(\bar G, \bar X)$ and $\rinj {\bar G}{\bar X}$.

\nota In this section we work with three metric spaces, namely $X$, its cone-off $\dot X$ and the quotient $\bar X$.
Since the map $X \hookrightarrow \dot X$ is an embedding we use the same letter $x$ to designate a point of $X$ and its image in $\dot X$.
We write $\bar x$ for its image in $\bar X$.
Unless stated otherwise, we keep the notation $\distV$ (without mentioning the space) for the distances in $X$ or $\bar X$.
The metric on $\dot X$ will be denoted by $\distV[\dot X]$.

\paragraph{The space $\dot X$.}
Let us now fix $\rho >0$.
Its value will be made precise later (see Theorem~\ref{res: SC - small cancellation theorem}).
It should be thought as a very large parameter.
We are going to build a cone-off over $X$.
To measure the distortion between $X$ and the spaces that we attach we use the map $\mu : \R_+ \rightarrow \R_+$ studied in Proposition~\ref{res: map mu}.
It is an $a$-comparison map with $a = 1/24(1+1/\sinh \rho)$.
Let $(H,Y) \in \mathcal Q$.
We denote by $\distV[Y]$ the length metric on $Y$ induced by the restriction of $\distV[X]$ on $Y$.
Since $Y$ is strongly quasi-convex, $Y$ endowed with $\distV[Y]$ is a length space such that for all $y,y' \in Y$
\begin{equation}
\label{eqn: cone-off construction - strong quasi-convex}
	\dist[X]y{y'} \leq \dist[Y]y{y'} \leq \dist[X]y{y'} +8\delta.
\end{equation}
We consider the cone $Z(Y)$ of radius $\rho$ over $\left(Y, \distV[Y]\right)$.
It comes with a map $\iota : Y \hookrightarrow Z(Y)$ as defined in Section~\ref{sec: cone}.
It follows from~(\ref{eqn: cone-off construction - strong quasi-convex}) and the properties of $\mu$ that for every $y,y' \in Y$.
\begin{equation}
\label{eqn: SC - distortion control}
	\mu\left(\dist[X] y{y'}\right) \leq \dist[Z(Y)]{\iota(y)}{\iota(y')} \leq \mu\left(\dist[X] y{y'}\right) +8\delta
\end{equation}
The set $\mathcal Z$ of all triples $(Z(Y),Y,\iota)$ constructed in this way is a $(\mu,X)$-family.
\begin{defi}
\label{def: SC - definition of cone off}
	The \emph{cone-off} space $\dot X(\mathcal Z)$ or simply $\dot X$ is the space obtained by attaching on $X$ for every $(Z,Y, \iota) \in \mathcal Z$ the cones $Z(Y)$ along $Y$ according to $\iota$.
\end{defi}
By Proposition~\ref{res: length structure on the cone-off} $\dot X$ is a length space.
Note that $X$ is a deformation retract of $\dot X \setminus v(\mathcal Z)$.
Here $v(\mathcal Z)$ stands for the set of all apices of the cones $Z(Y)$.

\begin{lemm}
\label{res: SC - cone off simply-connected}
	The space $\dot X$ is $50\delta$-simply-connected.
\end{lemm}

\begin{proof}
	Since $X$ is a $\delta$-hyperbolic length space, it is $50\delta$-simply connected \cite[Chap. 5, Prop. 1.1]{CooDelPap90}.
	Moreover, the length of a loop contained in $X$ is shorter measured with $\distV[\dot X]$ than with $\distV[X]$.
	Therefore it is sufficient to show that any loop of $\dot X$ is homotopic to a loop in $X$.
	Recall that the subsets $Y$'s over which we built the cones are path connected.
	Therefore any loop of $\dot X$ can be homotoped to a loop avoiding the set of apices $v(\mathcal Z)$.
	The conclusion follows from the fact that $X$ is a deformation retract of $\dot X\setminus v(\mathcal Z)$.
\end{proof}

\paragraph{}
The action of $G$ on $X$ extends by homogeneity in an action on $\dot X$: if $x=(y,r)$ is a point of a cone $Z(Y)$ and $g$ an element of $G$ then $gx$ is the point of the cone $Z(gY)$ defined by $(gy,r)$.
It follows from the definition of $\distV[\dot X]$ that $G$ acts by isometries on $\dot X$.

\begin{lemm}
\label{res : SC - action of G on dot X almost proper}
	Let $x$ be a point of $\dot X$ in the $\alpha$-neighborhood of $X$.
	The set $S$ of elements $g \in G$ such that $\dist[\dot X]{gx}x < 2(\rho - \alpha)$ is finite.
\end{lemm}

\begin{proof}
	We denote by $p$ a projection of $x$ on $X$.
	Let $g \in S$.
	The point $gp$ is a projection of $gx$ on $X$.
	By the triangle inequality 
	\begin{displaymath}
		\mu\left(\dist {gp}p\right) \leq \dist[\dot X] {gp}p \leq \dist[\dot X]{gx}x +2 \alpha < 2\rho.
	\end{displaymath}
	Thus $\dist{gp}p < \pi \sinh \rho$.
	However the action of $G$ on $X$ is proper.
	Therefore the set of elements $g \in G$ such that $\dist{gp}p < \pi \sinh \rho$ is finite, hence so is $S$.
\end{proof}

\begin{prop}
\label{res: SC - cone off globally hyperbolic}
	There exist universal positive numbers $\rho_0 > 10^{20} \boldsymbol \delta$, $\delta_0$ and $\Delta_0$ (i.e. which do not depend on $X$, $G$ or $\mathcal Q$) with the following property.
	If $\rho \geq \rho_0$, $\delta \leq \delta_0$ and $\Delta(\mathcal Q) \leq \Delta_0$ then  $\dot X$ is $900\boldsymbol \delta$-hyperbolic.
\end{prop}

\rem Recall that $\boldsymbol \delta$ is the hyperbolicity constant of the hyperbolic plane $\H_2$.

\begin{proof}
	The proof falls in two steps.
	First we show that $\dot X$ is locally hyperbolic and then we apply the Cartan-Hadamard theorem.
	We first fix $\rho_0 > 0$ such that $\rho_0 > 10^{20} \boldsymbol \delta$ and $(1+1/\sinh \rho_0)/24 \leq 1/50$.
	Since $\rho \geq \rho_0$ the map $\mu$ is a 1/50-comparison map.
	On the other hand we know that for every $t \in \intval 0{\pi \sinh \rho}$, $\mu(t) \geq 2\arcsinh(t/\pi)$ (see Proposition~\ref{res: map mu}-\ref{emu: mu - lower bound}).
	In particular there exists $t>0$ which only depends on $\rho_0$ and $\boldsymbol \delta$ such that $\mu(t) > 8.10^19 \boldsymbol \delta$.
	We denote by $\delta_0(2\boldsymbol \delta, \boldsymbol \delta,1/50,t)$ and $\Delta_0(2\boldsymbol \delta, \boldsymbol \delta,1/50,t)$ the constants provided by Theorem~\ref{res: curvature of dot X - general case}.
	All the cones that we attach are $2\boldsymbol \delta$-hyperbolic length spaces.
	Moreover, the distortion is control by the inequality~(\ref{eqn: SC - distortion control}).
	Assume now that 
	\begin{eqnarray*}
		\delta & \leq & \min \left\{\fantomB \boldsymbol \delta, \delta_0(2\boldsymbol \delta,\boldsymbol \delta, 1/50,t)\right\}, \\
		\Delta(\mathcal Q) & \leq & \Delta_0(2\boldsymbol \delta,\boldsymbol \delta, 1/50,t).
	\end{eqnarray*}
	Then by Theorem~\ref{res: curvature of dot X - general case}, every ball of radius $\mu(t)/8$ of $\dot X$ is $3\boldsymbol \delta$-hyperbolic.
	By Lemma~\ref{res: SC - cone off simply-connected}, $\dot X$ is also $50\boldsymbol \delta$-simply-connected.
	However $\mu(t)/8 > 10^{19}\boldsymbol \delta$.
	It follows from the Cartan-Hadamard Theorem (Theorem~\ref{res: cartan hadamard}) that $\dot X$ is (globally) $900\boldsymbol \delta$-hyperbolic.
\end{proof}

\paragraph{A rotation family.}
To every $(H,Y) \in \mathcal Q$ we associate the pair $(H,v)$ where $v$ is the apex of the cone $Z(Y)$ viewed as a point of $\dot X$.
We denote by $\mathcal R$ the set of all pairs obtained in this way.
Using the notation of Section~\ref{sec: rotation family} the set of apices $v(\mathcal R)$ is the same as $v(\mathcal Z)$.

\begin{lemm}
\label{res: SC rotation family}
	If $T(\mathcal Q) \geq \pi \sinh \rho$ then the collection $\mathcal R$ is a $\sigma$-rotation family where $\sigma =2\rho$.
\end{lemm}

\begin{proof}
	Let $(H,v)$ be a pair of $\mathcal R$.
	Is corresponds to a pair $(H,Y)$ of $\mathcal Q$ where $v$ is the apex of the cone $Z(Y)$.
	We assumed that $H$ stabilizes $Y$, thus it fixes $v$.
	By Lemma~\ref{res: metric on  dot X and Zi coincide} the metric of $Z(Y)$ and $\dot X$ coincide on the ball $B(v,\rho/5)$.
	Since $T(\mathcal Q) \geq \pi \sinh \rho$ Lemma~\ref{res: group acting as a rotation on cone} implies that for every $x \in B(v,\rho/5)$, for every $h \in H\setminus\{1\}$
	$\dist[\dot X]{hx}x = 2\dist[\dot X] xv$.
	This proves Axiom~\ref{enu: rotation family - large angle}  of the definition of rotation family (see Definition~\ref{def: rotation family})
	Axiom~\ref{enu: rotation family - apices apart} follows from the fact that the distance between two apices in $\dot X$ is at least $2\rho$.
	Finally the family $\mathcal Q$ being $G$-invariant, so is $\mathcal R$.
\end{proof}

\paragraph{The space $\bar X$.}

Recall that $K$ is the (normal) subgroup of $G$ generated by the subgroups $H$ with $(H,Y) \in \mathcal Q$.

\begin{defi}
\label{def: SC - definition of the space bar X}
	The space $\bar X$ is the quotient of the cone off $\dot X$ by the subgroup $K$.
\end{defi}

\begin{prop}
\label{res: SC - bar X globally hyperbolic}
	There exist universal positive numbers $\rho_0>10^{20}\boldsymbol \delta$, $\delta_0$ and $\Delta_0$ (i.e. which do not depend on $X$, $G$ or $\mathcal Q$) with the following property.
	Assume that $\rho \geq \rho_0$ and $\delta \leq \delta_0$.
	If in addition $\Delta(\mathcal Q) \leq \Delta_0$ and $T(\mathcal Q) \geq \pi \sinh \rho$ then  $\bar X$ is a $\bar \delta$-hyperbolic space where $\bar \delta \leq 64.10^4\boldsymbol \delta$.
	The group $\bar G$ acts by isometries on it and for all $(H,v) \in \mathcal R$ the projection $G \twoheadrightarrow \bar G$ induces an isomorphism from $\stab v /H$ onto $\stab {\bar v}$. 
\end{prop}

\begin{proof}
	The constant $\sigma_0 = \sigma_0(900 \boldsymbol \delta)$ is the one given by the fundamental theorem of rotation families (Theorem~\ref{res: rotation family - fundamental theorem}).
	We denote by $\delta_0$, $\Delta_0$ and $\rho_0$ the parameters given by Proposition~\ref{res: SC - cone off globally hyperbolic}.
	By increasing if necessary  $\rho_0$ we may assume that $2 \rho_0 \geq \sigma_0$.
	It follows from Proposition~\ref{res: SC - cone off globally hyperbolic} that $\dot X$ is $900 \boldsymbol \delta$-hyperbolic.
	According to Lemma~\ref{res: SC rotation family} the collection $\mathcal R$ that we previously built is a $2\rho$-rotation family.
	However by assumption $2\rho \geq \sigma_0$, thus we can apply all the results from Section~\ref{sec: rotation family} about rotation families.
	Among others the space $\bar X$ is $\bar \delta$-hyperbolic with $\bar \delta \leq 64.10^4\boldsymbol \delta$.
	Moreover, $\bar G$ acts by isometries on it.
	The last statement is a consequence of Corollary~\ref{res: rotation family - stabilizer in K}.
\end{proof}

	\paragraph{}
	For the remainder of this paragraph we assume that $X$, $G$ and $\mathcal Q$ satisfy the assumptions of Proposition~\ref{res: SC - bar X globally hyperbolic}.

	\begin{lemm}
	\label{res: SC - isometries not stabilizing a cone}
		Let $v \in  v(\mathcal R)$. 
		Let $\bar g \in \bar G \setminus \stab{\bar v}$.
		For every $x \in\bar X$, $\dist{\bar g\bar x}{\bar x} \geq 2(\rho - \dist{\bar x}{\bar v})$.
	\end{lemm}
	
	\begin{proof}
		Since $\bar g$ does not fix $\bar v$ the distance between $\bar v$ and $\bar g\bar v$ is at least $2\rho$.
		It follows the triangle inequality that 
		\begin{displaymath}
			2\rho \leq \dist {\bar g \bar v}{\bar g \bar x} + \dist{\bar g\bar x}{\bar x} + \dist{\bar x}{\bar v} = \dist{\bar g\bar x}{\bar x} + 2\dist{\bar x}{\bar v}. \qedhere
		\end{displaymath}
	\end{proof}

	\begin{prop}
	\label{res : SC - action of bar G co-compact}
		The group $\bar G$ acts properly co-compactly on $\bar X$.
	\end{prop}

	\begin{proof}
		Let $\bar x$ be a point of $\bar X$. 
		We claim that the set of elements $\bar g \in \bar G$ such that $\bar g B(\bar x, \bar \delta) \cap B(\bar x, \bar \delta) \neq \emptyset$ is finite.
		This will prove that the action is proper.
		We distinguish two cases.
		
		\paragraph{Case 1.}\emph{There exists $(H,v)\in\mathcal R$ such that $\dist {\bar v}{\bar x} < \rho -\bar \delta$}. 
		Let $\bar g \in \bar G$.
		Assume that $B(\bar x, \bar \delta)$ intersects $\bar g B(\bar x, \bar \delta)$.
		In particular $\dist{\bar g \bar x}{\bar x} \leq 2 \bar \delta$.
		According to Lemma~\ref{res: SC - isometries not stabilizing a cone}, $\bar g$ fixes $\bar v$.
		However $\stab{\bar v}$ is isomorphic to $\stab v /H = \stab Y/H$ which is finite ($H$ acts co-compactly on $Y$).
		Consequently, $\bar g$ belongs to $\stab{\bar v}$, which is finite.
				
		\paragraph{Case 2.}\emph{The point $\bar x$ is $\bar \delta$-close to $\nu(X)$.}
		We denote by $x$ a pre-image of $\bar x$ in $\dot X$ and by $S$ the set of elements $g \in G$ such that $\dist[\dot X] {gx}x \leq 2\bar \delta$.
		Recall that $\rho \geq \rho_0 > 2\bar \delta$.
		According to Lemma~\ref{res : SC - action of G on dot X almost proper}, $S$ is finite.
		Let $\bar g \in \bar G$ such that $\bar gB(\bar x, \bar \delta)$ intersects $B(\bar x, \bar \delta)$.
		In particular $\dist{\bar g \bar x}{\bar x} < 2\bar \delta$.
		By Proposition~\ref{res: rotation family - nu local isometry} there exists $g \in G$ such that $\dist[\dot X]{gx}x = \dist {\bar g\bar x}{\bar x}$.
		Thus $\bar g$ belongs to the image in $\bar G$ of $S$, which is finite.
		
		\paragraph{} The space $\bar X/\bar G$ can be obtained by attaching on $X/G$ finitely many cones of radius $\rho$ over $Y/H$ where $(H,Y) \in \mathcal Q/G$.
		As $G$ (\resp $H$) acts co-compactly on $X$ (\resp $Y$) the space $X/G$ (\resp $Y/H$) is compact.
		Consequently, so is $\bar X/\bar G$.
		Thus the action of $\bar G$ on $\bar X$ is co-compact. 
	\end{proof}

	\begin{prop}
	\label{res : SC - bar X proper}
		The space $\bar X$ is proper and geodesic
	\end{prop}

	\begin{proof}
		The space $\bar X$ is a metric space endowed with an action of $\bar G$ which is proper and co-compact.
		It follows that $\bar X$ is complete and locally compact \cite[Chap. I.8, Ex. 8.4(1)]{BriHae99}.
		On the other hand $\dot X$ and thus $\bar X$ is a length space. 
		The Hopf-Rinow Theorem implies that $\bar X$ is geodesic \cite[Chap.1, Prop. 3.7]{BriHae99}.
	\end{proof}

\paragraph{Small cancellation theorem}
The previous results can be summarized in the following theorem. 
It is an analog of the well-known fact saying that a group whose presentation satisfies the usual $C''(\lambda)$ small cancellation assumption with $\lambda < 1/6$ is hyperbolic (see \cite[Appendix, Th.36]{GhyHar90}).

\begin{theo}[Small cancellation theorem]
\label{res: SC - small cancellation theorem}
	There exist positive constants $\rho_0$, $\delta_0$ and $\Delta_0$ satisfying the following property.
	Let $G$ be a group acting properly co-compactly on a geodesic proper $\delta$-hyperbolic space.
	Let $\mathcal Q$ be a family of pairs $(H,Y)$ such that $Y$ is a strongly quasi-convex subset of $X$ and $H$ a subgroup of $G$ stabilizing $Y$.
	We assume that
\begin{enumerate}
	\item for every $(H,Y) \in \mathcal Q$, $H$ acts co-compactly on $Y$,
	\item $\mathcal Q$ is $G$-invariant i.e., for every $(H,Y) \in \mathcal Q$ for every $g \in G$ $(gHg^{-1}, gY)$ is still an element of $\mathcal Q$,
	\item the quotient $\mathcal Q/G$ is finite. 
\end{enumerate}
	Let $K$ be the (normal) subgroup of $G$ generated by the subgroups $H$'s where $(H,Y) \in \mathcal Q$.
	Let $\rho \geq \rho_0$.
	Let $\dot X$ be the cone-off space obtained by attaching for every $(H,Y) \in \mathcal Q$ a cone of radius $\rho$ and base $Y$ on $X$.
	Let $\bar X$ be the quotient of $\dot X$ by $K$.
	
	\paragraph{}If $\delta \leq \delta_0$, $\Delta(\mathcal Q) \leq \Delta_0$ and $T(\mathcal Q) \geq \pi \sinh \rho$ then $\bar X$ is a geodesic proper $\bar \delta$-hyperbolic space with $\bar \delta \leq 64.10^4 \boldsymbol \delta$.
	The group $\bar G = G/K$ acts properly co-compactly by isometries on $\bar X$.
	For every $(H,Y) \in \mathcal Q$, the projection $G \twoheadrightarrow \bar G$ induces an embedding $\stab Y/H \hookrightarrow \bar G$.
\end{theo}

\rems Note that in the theorem the constants $\delta_0$ and $\Delta_0$ (\resp $\rho_0$) can be chosen arbitrary small (\resp large).
From now on, we will always assume that $\rho_0 > 10^{20} \boldsymbol \delta$ whereas $\delta_0, \Delta_0 < 10^{-10}\boldsymbol\delta$.
These estimates are absolutely not optimal.
We chose them very generously to be sure that all the inequalities that we might need later will be satisfied.
What really matters is their orders of magnitude recalled below.
\begin{displaymath}
	\max\left\{\delta_0, \Delta_0\right\} \ll \boldsymbol \delta  \ll \rho_0 \ll \pi \sinh \rho_0.
\end{displaymath}
An other important point to remember is the following.
The constants $\delta_0$, $\Delta_0$ and $\pi \sinh \rho_0$ are used to describe the geometry of $X$ whereas $\boldsymbol \delta$ and $\rho_0$ refers to the one of $\dot X$ or $\bar X$.
This theorem looks slightly different from Theorem~\ref{res: rotation small cancellation theorem rescaled} given in the introduction.
The later follows from a rescaling argument.
However we prefer to retain that formulation which will make easier to keep track of the parameters.
In Section~\ref{sec: gromov monster} we give an example where we apply this rescaling argument.

\paragraph{Comparison with the original proof.} At this stage, we can explain precisely the main differences between our proof and the one by T.~Delzant and M.~Gromov.
In \cite{DelGro08}, the authors construct an orbifold $\bar Q$ by attaching on $X/G$ appropriate cones. 
In this context the space $\dot X$ can be used as a chart for the neighborhood of the points of $\bar Q$ close to $X/G$.
They interpret $\bar G$ as the fundamental group (in the sense of orbifolds) of $\bar Q$.
In this description the groups $\stab Y/H$ correspond to the stabilizers of the singularities of the orbifold.
They use then a version of the Cartan-Hadamard theorem for orbifolds (\cite[Th. 4.3.1]{DelGro08}).
They prove that if an orbifold is locally hyperbolic (with the correct quantifiers) then it is developable and its fundamental cover  $\bar X$ (in the sense of orbifolds) is globally hyperbolic.
One consequence of the developability is that the stabilizers of the singularities, namely $\stab Y/H$, embeds into $\bar G$.
An other is that the developing map $\nu : \dot X \rightarrow \bar X$ induces a local isometry.
In our approach we substitute the topological point of view for a study of rotation families acting on $\dot X$.
It allows us to construct directly $\bar X$ as a quotient of $\dot X$ without going through $\bar Q$. 
The two consequences of developability mentioned before correspond to our Corollary~\ref{res: rotation family - stabilizer in K} and Proposition~\ref{res: rotation family - nu local isometry}.
T.~Delzant and M.~Gromov also used the description in terms of orbifold to lift elements from $\bar G$ to $G$.
In particular if two paths of $\bar X$ with the same extremities are homotopic relative to their endpoints in $\bar X\setminus \bar v(\mathcal R)$, then their lifts in $\dot X$ will define the same endpoint in $\dot X$.
We prove instead Propositions~\ref{res: lift quasi-convex isometry} and \ref{res: lift quasi-convex isometry and stabilizer}.
An other difference concerns the framework used to describe the geometry of $\dot X$ and $\bar X$.
We only proved that these spaces were hyperbolic.
They are actually endowed with a finer metric structure.
As explained in \cite{DelGro08} they are locally $\operatorname{CAT}(-1,\epsilon)$, which roughly means that every ball of large radius satisfies the $\operatorname{CAT}(-1)$ condition up to an error $\epsilon$.
However this refinement is not needed to prove the infiniteness of Burnside groups. 

\paragraph{}
For the remainder of Section~\ref{sec: small cancellation theory} we assume that $X$, $G$ and $\mathcal Q$ are as in Theorem~\ref{res: SC - small cancellation theorem}. 

\subsection{Isometries of the quotient}

\paragraph{} Since $\bar G$ acts properly co-compactly by isometries on a hyperbolic space, its elements are either elliptic or hyperbolic.
In this section we study how the type of an isometry is related to the one of its preimage in $G$ (seen as an isometry of $X$).

\begin{prop}
\label{res: SC  - lifting elliptic subgroups}
	Let $\bar F$ be a finite subgroup of $\bar G$.
	Either there exists $v \in v(\mathcal R)$ such that $\bar F$ is contained in $\stab{\bar v}$ or $\bar F$ is isomorphic to a finite subgroup of $G$.
\end{prop}

\begin{proof}
	Recall that $C_{\bar F}$ is the set of points $\bar x \in \bar X$ such that for every $\bar g \in \bar F$, $\dist{\bar g\bar x}{\bar x} \leq 10\bar \delta$ (see Definition~\ref{def: characteristic set of a finite subgroup}). 
	It is $\bar F$-invariant and $8\bar \delta$-quasi-convex (see Corollary~\ref{res: characteristic subset elliptic group - quasi convex}).
	If $\bar F$ is not contained in some $\stab{\bar v}$, then according to Lemma~\ref{res: SC - isometries not stabilizing a cone} for every $v \in v(\mathcal R)$, $C_{\bar F}$ does not intersect $B(\bar v, \rho-5\bar \delta)$.
	 By Proposition~\ref{res: lift quasi-convex isometry and stabilizer}, there exists a subset $C$ of $\dot X$ such that the map $\nu : \dot X \rightarrow \bar X$ induces an isometry from $C$ onto $C_{\bar F}$ and the projection $\pi : G \rightarrow \bar G$ induces an isomorphism from $\stab C$ onto $\stab {C_{\bar F}}$.
	 In particular $\bar F$ is isomorphic to a subgroup $F$ of $\stab C$.
\end{proof}

\begin{prop}
\label{res: SC - lower bound injectivity radius}
	We denote by $l$ the greatest lower bound on the stable translation length (in $X$) of hyperbolic elements of $G$ which do not belong to some $\stab v$, $v \in v(\mathcal R)$.
	Then $\rinj{\bar G}{\bar X} \geq \min\left\{ \kappa l/8, \bar \delta\right\}$ where $\kappa = 2\rho / \pi\sinh \rho$.
\end{prop}

\begin{proof}
	Let $\bar g$ be a hyperbolic element of $\bar G$.
	Recall that for every $m \in \N$ we have $m\len[stable]{\bar g} \geq \len{\bar g^m} - 32\bar \delta$.
	Therefore it suffices to find an integer $m$ such that $\len{\bar g^m} \geq  m\min\left\{ \kappa l/8, \bar \delta\right\} +32\bar \delta$.
	We denote by $m$ the largest integer satisfying $m \min\left\{ \kappa l/8, \bar \delta\right\} \leq  40\bar \delta$.
	Assume that $\len{\bar g^m}$ is smaller than $m\min\left\{ \kappa l/8, \bar \delta\right\} +32\bar \delta$.
	In particular $\len{\bar g^m} \leq 72\bar \delta$.
	However $\bar g$ is hyperbolic thus it cannot belong to some $\stab{\bar v}$.
	According to Lemma~\ref{res: SC - isometries not stabilizing a cone} for every $v \in v(\mathcal R)$, $A_{\bar g^m}$ does not intersect $B(\bar v, \rho-36\bar \delta)$.
	By Proposition ~\ref{res: lift quasi-convex isometry and stabilizer}, there exists a subset $A$ of $\dot X$ such that the map $\nu : \dot X \rightarrow \bar X$ induces an isometry from $A$ onto $A_{\bar g^m}$ and the projection $\pi : G \rightarrow \bar G$ induces an isomorphism from $\stab A$ onto $\stab {A_{\bar g^m}}$.
	We denote by $g$ the preimage of $\bar g$ in $\stab A$.
	It is a hyperbolic element of $G$.
	In particular $\len[stable]g \geq l$.
	Let $\bar x$ be a point of $A_{\bar g^m}$, $x$ the preimage of $\bar x$ in $A$ and $y$ a projection of $x$ on $X$.
	In particular $g^my$ is a projection of $g^mx$ on $X$ and $\dist[\dot X] xy \leq 36\bar \delta$.
	Moreover, we have
	\begin{displaymath}
		\mu\left(\dist{g^my}y\right) \leq \dist[\dot X]{g^mx}x + 72\bar \delta = \dist{\bar g^m \bar x}{\bar x} +72 \bar \delta \leq \max\left\{\len{\bar g^m} , 8 \bar \delta \right\} + 72 \bar \delta \leq 144 \bar \delta <2\rho.
	\end{displaymath}
	It follows that $\dist{g^m y}y < \pi \sinh \rho$.
	The function $\mu$ being concave we get $ml \leq \dist{g^my}y \leq 144\kappa^{-1}\bar \delta$, which contradicts the maximality of $m$.
\end{proof}

\subsection{Groups without even-torsion}
\label{sec: groups without even-torsion}

\paragraph{}
This section is specific to the study of Burnside groups of odd exponents.
As explained in the introduction, we will construct by induction a sequence of groups by adjoining at each steps relations of the form $r^n$.
If $n$ is large enough, these relations satisfy a small cancellation assumption.
The difficulty is to prove that the same exponent $n$ can be used at every step.
This can be achieved by estimating the invariant $A$.
More precisely we provide in this section an upper bound for $A(\bar G, \bar X)$ in terms of  $A(G,X)$ and $\bar \delta$.
To that end we make the following hypotheses.

\begin{enumerate}
	\item Every elementary subgroup of $G$ is cyclic - infinite or of finite of odd order.
	\item For every $(H,Y) \in \mathcal Q$, there exists $r \in G$ which is hyperbolic and not a proper power such that $Y$ is the cylinder of $r$ (see definition in Section~\ref{sec: isometries hyperbolic space}) and $H$ the cyclic subgroup of $G$ generated by an odd power of $r$.
\end{enumerate}
We still require that $X$, $G$ and $\mathcal Q$ satisfy the assumptions of Theorem~\ref{res: SC - small cancellation theorem}.
\begin{prop}
	\label{res: SC - lifting overlap of axes}
	Let $\bar g$ and $\bar h$ be two elements of $\bar G$ such that $\len {\bar g}$ and $\len {\bar h}$ are at most $1000 \bar \delta$.
	One of the following holds.
	\begin{enumerate}
		\item There exists $v \in v(\mathcal R)$ such that $\bar g$ and $\bar h$ belong to $\stab{\bar v}$.
		\item There exist respective preimages $g$ and $h$ in $G$  of $\bar g$ and $\bar h$ such that $\len g$ and $\len h$ are at most $\pi \sinh(1034 \bar \delta)$ and
		\begin{displaymath}
			\diam \left( A_{\bar g}^{+ 17\bar \delta} \cap A_{\bar h}^{+ 17\bar \delta} \right) 
			\leq \diam \left( A_g^{+ 17\delta} \cap A_h^{+ 17 \delta} \right) + \pi \sinh (2087\bar \delta).
		\end{displaymath}
	\end{enumerate}
\end{prop}

\rems In the statement of the proposition all the metric objects are measured either with the distance of $X$ or $\bar X$, but not with the one of $\dot X$.
Note that this result actually holds without our additional assumptions on the torsion.

\begin{proof}
	Without loss of generality we can assume that the intersection of the respective $17 \bar \delta$-neighborhoods of $A_{\bar g}$ and  $A_{\bar h}$ is not empty.
	Let us call $\bar Z$ this intersection.
	Assume that there exists $v\in v(\mathcal R)$ and a point $\bar z \in \bar Z$ such that $\dist {\bar v}{\bar z} < \rho -  517\bar \delta$.
	By definition $\bar g$ and $\bar h$ move $\bar z$ by a distance at most $1034 \bar \delta$.
	According to Lemma~\ref{res: SC - isometries not stabilizing a cone}, they belong to $\stab {\bar v}$, which provides the first case.
	
	\paragraph{}
	We now assume that for every $v \in v(\mathcal R)$, $\bar Z$ does not intersect $B(\bar v, \rho -517\bar \delta)$.
	By Lemma~\ref{res: intersection of quasi-convex}, $\bar Z$ is $7\bar \delta$-quasi-convex.
	Moreover, $\bar g$ and $\bar h$ move any point of $\bar Z$ by at most $1034 \bar \delta$.
	According to Proposition~\ref{res: lift quasi-convex isometry and stabilizer}, there exists a subset $Z$ of $\dot X$ and respective preimages $g$ and $h$ of $\bar g$ and $\bar h$ satisfying the following properties.
	\begin{enumerate}
		\item The map $\nu : \dot X \rightarrow \bar X$ induces an isometry from $Z$ onto $\bar Z$.
		\item For every $z \in Z$  we have $\dist[\dot X]{gz}z = \dist{\bar g\bar z}{\bar z}$ and $\dist[\dot X]{hz}z = \dist{\bar h\bar z}{\bar z}$.
	\end{enumerate}
	We now denote by $\bar z$ and $\bar z'$ two points of $\bar Z$ such that
	\begin{displaymath}
		\dist{\bar z}{\bar z'} \geq \diam\left(A_{\bar g}^{+ 17\bar \delta} \cap A_{\bar h}^{+ 17\bar \delta} \right)  - \bar \delta. 
	\end{displaymath}
	The points $z$ and $z'$ stands for their preimages in $Z$.
	We write $x$ and $x'$ for respective projections of $z$ and $z'$ on $X$.
	By assumption, $\bar Z$ lies in the $517\bar \delta$-neighborhood of $\nu(X)$.
	Thus $\dist[\dot X] xz, \dist[\dot X] {x'}{z'} \leq 517\bar \delta$.
	In particular,
	\begin{displaymath}
		\mu \left(\dist{gx}x\right) \leq \dist[\dot X]{gx}{x} \leq \dist{\bar g \bar z}{\bar z} + 1034\bar \delta \leq 2068\bar \delta < \pi \sinh \rho.
	\end{displaymath}
	It follows that $\dist {gx}x \leq \pi \sinh (1034 \bar \delta)$.
	The same holds for $x'$.
	Consequently, $\len {g} \leq \pi \sinh (1034\bar \delta)$.
	Moreover, $x$ and $x'$ belong to the $A$-neighborhood of $A_g$ with $A = \pi \sinh (1034\bar \delta)/2 +7\bar \delta$.
	Similarly they belong to the $A$-neighborhood of $A_h$.
	By Lemma~\ref{res: intersection of thickened quasi-convex}
	\begin{displaymath}
		\dist x{x'} \leq \diam\left(A_g^{+ 17\delta} \cap A_h^{+17 \delta}\right) +  \pi \sinh (1034\bar \delta) +14\bar \delta + 4\delta
	\end{displaymath}
	On the other hand, the map $X \rightarrow \dot X$ shorten the distances.
	Therefore
	\begin{displaymath}
		\dist x{x'} \geq \dist[\dot X] x{x'} \geq \dist[\dot X] z{z'} - 1034\bar \delta \geq \dist {\bar z}{\bar z'} - 1034\bar \delta \geq  \diam\left(A_{\bar g}^{+ 17\bar \delta} \cap A_{\bar h}^{+ 17\bar \delta} \right) - 1035\bar \delta. 
	\end{displaymath}
	The conclusion of the second case follows from the last two inequalities.
\end{proof}

	\begin{coro}
	\label{res: SC - estimate A for the quotient}
		The invariant $A(\bar G, \bar X)$ is bounded above by $A(G,X) + \pi \sinh (10^4\bar \delta)$.
	\end{coro}
	
	\begin{proof}
		We denote by $\mathcal A$ the set of pairs $(\bar g, \bar h)$ of $\bar G$  generating a non-elementary subgroup such that $\len {\bar g}$ and $\len {\bar h}$ are at most$1000 \bar \delta$.
		Let $(\bar g, \bar h) \in \mathcal A$.
		By definition of $\mathcal A$, $\bar g$ and $\bar H$ cannot be both in some $\stab {\bar v}$.
		According to Proposition~\ref{res: SC - lifting overlap of axes}, there exist respective preimages $g$ and $h$ in $G$  of $\bar g$ and $\bar h$ such that $\len g$ and $\len h$ are at most $\pi \sinh(1034 \bar \delta)$ and
		\begin{displaymath}
			\diam \left( A_{\bar g}^{+ 17\bar \delta} \cap A_{\bar h}^{+ 17\bar \delta} \right) 
			\leq \diam \left( A_g^{+ 17\delta} \cap A_h^{+ 17 \delta} \right) + \pi \sinh (2087\bar \delta).
		\end{displaymath}
		However we assumed that every elementary subgroup of $G$ is cyclic.
		Thus by Proposition~\ref{res: overlap two axes},
		\begin{displaymath}
			\diam \left( A_{\bar g}^{+ 17\bar \delta} \cap A_{\bar h}^{+ 17\bar \delta} \right) \leq 2 ( \len g + \len h) + A (G,X) + 667\delta \leq A(G,X) + \pi \sinh (10^4\bar \delta).
		\end{displaymath}
		The last inequality holds for every $(\bar g, \bar h) \in \mathcal A$, which leads to the result.
	\end{proof}
	
	\begin{lemm}
	\label{res - SC no more than one fixed apex}
		If $T(\mathcal Q) \geq 4 \pi \sinh \rho$ then an element $\bar g \in \bar G \setminus \{1\}$ cannot fix more than one apex of $\bar v( \mathcal R)$.
	\end{lemm}
	
	\begin{proof}
		Assume that $\bar g$ fixes an apex of $\bar X$.
		There exists $(H, Y) \in \mathcal Q$ such that the apex $\bar v$ fixed by $\bar g$ is the image in $\bar X$ of the apex of the cone $Z(Y)$.
		We proved that $\stab {\bar v}$ is isomorphic to $\stab v /H = \stab Y/H$.
		However $\stab Y$ is a cyclic group generated by $r\in G$.
		Therefore there exists $m \in \Z$ such that $\bar g$ is the image of $g = r^m$.
		Since $\bar g$ is not trivial $m \neq 0 \operatorname{mod} n$.
		Thus there exists an integer $p \in m\Z + n\Z$ between $n/3$ and $2n/3$.
		In particular $\bar r^p$ is a power of $\bar g$.
		Let $x = (y,r)$ be a point of the cone $Z(Y)$.
		By construction of $p$ we have
		\begin{displaymath}
			\dist {r^py}y \geq p\len[stable] r \geq n\len[stable] r /3 \geq T(\mathcal Q)/3 -32 \delta \geq \pi \sinh \rho.
		\end{displaymath}
		Consequently, $\dist[\dot X] {r^px}x = 2r$.
		On the other hand $y$ is a point of the cylinder of $r$ and thus is contained in the $36\delta$-neighborhood of the axis of $r^p$ (see Lemma~\ref{res: cylinder contained in invariant subset}).
		Hence
		\begin{displaymath}
			\dist {r^py}y \leq \len {r^p} + 72\delta \leq p\len[stable] r + 104\delta \leq \len{r^n} + 104\delta - (n-p)\len [stable] r \leq \len {r^n} - \pi \sinh \rho.
		\end{displaymath}
		According to Lemma~\ref{res: metric quotient cone}, $\dist{\bar r^p \bar x}{\bar x} = \dist[\dot X]{r^p x}x = 2r$.
		It follows that the points of $Z(Y)/H$ contained in the axis of $\bar r^p$ are $4\bar \delta$-close to $\bar v$.
		This axis is nevertheless $14\bar \delta$-quasi-convex.
		Thus it is contained in $B(\bar v, 5\bar \delta)$.
		In particular, $\bar r^p$ cannot fix an other apex than $\bar v$.
		Hence $\bar g$ cannot either. 
	\end{proof}

	\begin{lemm}
	\label{res: SC - no even torsion}
		The group $\bar G$ has no element of order 2.
	\end{lemm}
	
	\begin{proof}
		According to Proposition~\ref{res: SC  - lifting elliptic subgroups}, every element of $\bar G$ with finite order is contained either in the image of a finite subgroup of $G$ or in some $\stab {\bar v}$, with $v \in v(\mathcal R)$.
		It follows from our assumptions at the beginning of the section that none of these groups contain an element of order 2.
		Thus $\bar G$ cannot have even-torsion.
	\end{proof}
	
	\begin{prop}
	\label{res : SC - every elementary subgroup is cyclic}
		If $T(\mathcal Q) \geq 4 \pi \sinh \rho$ then every elementary subgroup of $\bar G$ is cyclic.
	\end{prop}
	
	\begin{proof}
		Let $\bar E$ be an elementary subgroup of $\bar G$.
		Assume first that $\bar E$ is finite.
		By Proposition~\ref{res: SC  - lifting elliptic subgroups} $\bar E$ is either isomorphic to a subgroup of $G$ or contained in some $\stab {\bar v} = \stab v/H$ with $(H,v) \in \mathcal R$. 
		However we assumed that every elementary subgroup of $G$ is cyclic.
		On the other hand if $v$ is the apex of the cone built on $Y$ with $(H,Y) \in \mathcal Q$ then $\stab {\bar v}$ is isomorphic to $\stab Y / H$.
		According to our assumption $\stab Y / H$ is also cyclic.
		Hence so is $\bar E$.
		
		\paragraph{}
		Suppose now that $\bar E$ is infinite. 
		By Lemma~\ref{res: SC - no even torsion} it does not contain an element of order 2. 
		Therefore $\bar E$ is isomorphic to the semi-direct product $\sdp {\bar F} \Z$ where $\bar F$ is the maximal finite subgroup of $\bar E$, $\Z$ is generated by a hyperbolic element $\bar g$ acting by conjugation on $\bar F$.
		We claim that the cylinder $Y_{\bar g}$ of $\bar g$ lies in the $53\bar \delta$-neighborhood of $\nu(X)$.
		Assume on the contrary that there exists $v \in v(\mathcal R)$ and a point $\bar x \in Y_{\bar g}$ such that $\dist {\bar x}{\bar v} < \rho - 53\bar \delta$.
		By Lemma~\ref{res: cylinder almost invariant by max finite subgroup} $Y_{\bar g}$ lies in the $48\bar \delta$-neighborhood of $C_{\bar F}$.
		Thus for every $\bar u \in \bar F$, $\dist {\bar u \bar x}{\bar x} \leq 106\bar \delta$.
		According to Lemma~\ref{res: SC - isometries not stabilizing a cone}  $\bar F$ is contained in $\stab {\bar v}$.
		However $\bar g$ normalizes $\bar F$, therefore $\bar F$ also fixes $\bar g \bar v$.
		It follows from Lemma~\ref{res - SC no more than one fixed apex} that $\bar g \bar v = \bar v$.
		In particular $\bar g$ cannot be hyperbolic.
		Contradiction.
		The cylinder $Y_{\bar g}$ is a $2\bar \delta$-quasi-convex subset of $\bar X$ contained in the $53\bar \delta$-neighborhood of $\nu(X)$.
		Applying Proposition~\ref{res: lift quasi-convex isometry and stabilizer}, there exists a subset $Y$ of $\dot X$ such that 
		\begin{enumerate}
			\item The map $\nu : \dot X \rightarrow \bar X$ induces an isometry from $Y$ onto $Y_{\bar g}$
			\item The projection $G  \rightarrow \bar G$ induces an automorphism from $\stab Y$ onto $\stab {Y_{\bar g}}$.
		\end{enumerate}
		However $\bar E$ stabilizes $Y_{\bar g}$ thus it is isomorphic to a subgroup of $G$, which is necessarily elementary.
		All the elementary subgroups of $G$ being cyclic, so is $\bar E$.
	\end{proof}
	
	\begin{prop}
	\label{res : SC - non elementary quotient}
		Assume that there exists two isometries $u, v \in G$ and a point $x \in X$ such that $0 <\dist {uvu^{-1}v^{-1}x}x < \rho/20$.
		Then the same holds in $\bar X$.
		Moreover $\bar G$ is not elementary.
	\end{prop}

	\begin{proof}
	In order to simplify the notations we put $g = uvu^{-1}v^{-1}$.
	The natural map $X \rightarrow \dot X$ shortens the distances.
	More precisely we have
	\begin{displaymath}
		0 < \mu \left( \dist {gx}x\right) \leq \dist[\dot X] {gx}x \leq \dist {gx}x < \rho/20.
	\end{displaymath}
	However the map $\nu : \dot X \rightarrow \bar X$ induces an isometry from the ball of center $x$ and radius $\rho/20$ onto its image (see Proposition~\ref{res: rotation family - nu local isometry}).
	Consequently, $0 < \dist {\bar g \bar x}{\bar x} < \rho/20$.
	Assume now that $\bar G$ is elementary. 
	By Proposition~\ref{res : SC - every elementary subgroup is cyclic}, it is cyclic and thus abelian.
	This contradicts the fact that $\dist {\bar g \bar x}{\bar x} > 0$.
	\end{proof}

\section{Applications}
\label{sec: applications}

\subsection{Periodic quotients of hyperbolic groups}
\label{sec: application to periodic groups}

The next proposition will play the role of the induction step in the proof of the main theorem.
\begin{prop}
\label{res: SC - induction lemma}
	There exist positive constants $\rho_0$,  $\delta_1$ and an integer $n_0$ with the following properties.
	Let $G$ be a group acting properly co-compactly by isometries on a geodesic proper $\delta_1$-hyperbolic space $X$.
	Let $n_1 \geq n_0$  and $n \geq n_1$ be an odd integer.
	We denote by $R$ the set of hyperbolic elements $r$ of $G$ which are not a proper power and such that $\len r \leq 1000\delta_1$.
	Let $K$ be the normal subgroup of $G$ generated by $\{r^n, r \in R\}$ and $\bar G$ the quotient of $G$ by $K$.
	We make the following assumptions.
	\begin{enumerate}
		\item \label{enu: SC - induction lemma - elementary subgroups}
		$G$ is non elementary. Every elementary subgroup of $G$ is cyclic either infinite of finite with order dividing $n$.
		\item \label{enu: SC - induction lemma - A}
		$A(G,X) \leq 2\pi \sinh (10^4\delta_1)$.
		\item \label{enu: SC - induction lemma - rinj}
		$\rinj GX \geq 20\sqrt {\rho_0\delta_1/n_1}$.
		\item \label{enu: SC - induction lemma - non-commutating elements}
		There exist $u,v \in G$ and a point $x \in X$ such that 
		\begin{displaymath}
			0 <\dist {uvu^{-1}v^{-1}x}x \leq \rho_0/20.
		\end{displaymath}
	\end{enumerate}
	Then there exists a geodesic proper $\delta_1$-hyperbolic space $\bar X$ on which $\bar G$ acts properly, co-compactly by isometries.
	The action of $\bar G$ on $\bar X$ satisfies also Points \ref{enu: SC - induction lemma - elementary subgroups}-\ref{enu: SC - induction lemma - non-commutating elements}.
	Moreover, for every $g \in G$, if $\bar g$ stands for its image in $\bar G$ we have 
	\begin{displaymath}
			\len[stable, espace= \bar X]{\bar g} \leq \frac 1{\sqrt {n_1}} \left(\frac {\pi \sinh \rho_0}{5\sqrt{\rho_0\delta_1}}\right)\len[stable, espace= X]g. 
		\end{displaymath}
	\end{prop}

\voc Let $G$ be a group acting by isometries on a space $X$ and $n$ an integer.
Once $n_1$ and $n$ have been fixed, if $G$ and $X$ satisfy the assumption of the proposition including Points \ref{enu: SC - induction lemma - elementary subgroups}-\ref{enu: SC - induction lemma - non-commutating elements}, we will write that $(G,X)$ \emph{satisfies the induction hypotheses for exponent $n$}.
The proposition says in particular that if $(G,X)$ satisfies the induction hypotheses for exponent $n$ then so does $(\bar G, \bar X)$.

\begin{proof}
	The parameters $\rho_0$, $\delta_0$ and $\Delta_0$ are the one given by the small cancellation theorem (Theorem~\ref{res: SC - small cancellation theorem}).
	We set $\delta_1 = 64.10^4\boldsymbol \delta$.
	We now define the critical exponent $n_0$.
	To that end we consider a rescaling parameter $\lambda_n$ depending on an integer $n$
	\begin{displaymath}
		\lambda_n =\frac {\pi \sinh \rho_0}{ 5\sqrt{n\rho_0\delta_1}}
	\end{displaymath}
	The sequence $(\lambda_n)$ converges to 0 as $n$ approaches infinity. 
	Therefore there exists an integer $n_0$ such that for every $n \geq n_0$
	\begin{eqnarray}
		\label{eqn: induction - delta}
		\lambda_n \delta_1 &  \leq & \delta_0 \\
		\label{eqn: induction - Delta}
		\lambda_n\left(2\pi \sinh(10^4\delta_1) +  86\delta_1\right) & \leq & \min \left\{\Delta_0, \pi \sinh (10^4\delta_1) \right\} \\
		\label{eqn: induction - rinj}
		\frac{100\lambda_n\rho_0\delta_1}{\pi \sinh \rho_0} & \leq & \delta_1 \\
		\label{eqn: induction - non elem}
		\lambda_n \rho_0 & \leq & \rho_0
	\end{eqnarray}
	Let $n_1 \geq n_0$ and $n \geq n_1$ be an odd integer.
	For simplicity of notation we denote by $\lambda$ the rescaling parameter $\lambda = \lambda_{n_1}$.
	Let $G$ be a group acting by isometries on a metric space $X$ such that $(G,X)$ satisfies the induction hypotheses for exponent $n$.
	We denote by $R$ the set of hyperbolic elements $r$ of $G$ which are not a proper power and such that $\len r \leq 1000\delta_1$.
	Let $K$ be the normal subgroup of $G$ generated by $\{r^n, r \in R\}$ and $\bar G$ the quotient of $G$ by $K$.
	We are going to prove that $\bar G$ is a small cancellation quotient of $G$.
	To that end we consider the action of $G$ on the rescaled space $\lambda X$.
	In particular it is a $\delta$-hyperbolic space, with $\delta = \lambda \delta_1 \leq \delta_0$.
	Unless stated otherwise, we will always work with the rescaled space $\lambda X$.
	We define the family $\mathcal Q$ by 
	\begin{displaymath}
		\mathcal Q = \set{\fantomB\left(\left< r^n\right>,Y_r\right)}{r \in R}.
	\end{displaymath}

	\begin{lemm}
	\label{res: SC - induction lemma - small cancellation}
		The family $\mathcal Q$ satisfies the following assumptions: $\Delta \left( \mathcal Q \right) \leq \Delta_0$ and $T\left(\mathcal Q\right) \geq 4\pi \sinh \rho_0$.	
	\end{lemm}
	
	\begin{proof}
		We start with the upper bound of $\Delta (\mathcal Q)$.
		Let $r_1$ and $r_2$ be two elements of $R$ such that $(\langle r_1^n\rangle,Y_{r_1}) \neq (\langle r_2^n\rangle,Y_{r_2})$.
		According to Lemma~\ref{res: cylinder contained in invariant subset}, $Y_{r_1}$ and $Y_{r_2}$ are respectively contained in the $36\delta$-neighborhood of $A_{r_1}$ and $A_{r_2}$, thus by Lemma~\ref{res: intersection of thickened quasi-convex}
		\begin{displaymath}
			\diam\left( Y_{r_1}^{+ 5\delta} \cap Y_{r_2}^{+ 5\delta} \right) \leq \diam\left( A_{r_1}^{+17 \delta} \cap A_{r_2}^{+ 17\delta} \right) + 86\delta.
		\end{displaymath}
		According to Lemma~\ref{res: extracting subset without inverse}, $r_1$ and $r_2$ generate a non-elementary subgroup.
		On the other hand, their translation lengths in $\lambda X$ are at most $1000 \delta$, thus
		\begin{eqnarray*}
			\diam\left( Y_{r_1}^{+ 5\delta} \cap Y_{r_2}^{+ 5\delta} \right) \leq A(G, \lambda X) + 86\delta 
			& \leq & \lambda A(G,X) + 86\lambda\delta_1 \\
			& \leq & \lambda(2\pi \sinh(10^4\delta_1) + 86\delta_1).
		\end{eqnarray*}
		Thus by (\ref{eqn: induction - Delta}), $\Delta(\mathcal Q) \leq \Delta_0$.
		Let us focus now on $T(\mathcal Q)$.
		The injectivity radius of $G$ on $\lambda X$ is bounded below as follows
		\begin{displaymath}
			\rinj G{\lambda X} \geq \lambda\frac{20\sqrt{\rho_0\delta_1}}{\sqrt {n_1}} \geq  \frac{4\pi \sinh \rho_0}{n_1} \geq   \frac{4\pi \sinh \rho_0}n 
		\end{displaymath}
		In particular for every $r \in R$ we have $\len[stable]{r^n} = n \len[stable] r \geq 4\pi \sinh \rho_0$. 
		Hence $T(\mathcal Q) \geq 4\pi \sinh \rho_0$.
	\end{proof}
	
	\paragraph{}On account of the previous lemma, we can now apply the small cancellation theorem (Theorem~\ref{res: SC - small cancellation theorem}) to the action of $G$ on the rescaled space $\lambda  X$ and the family $\mathcal Q$.
	We denote by $\dot X$ the space obtained by attaching on $\lambda X$ for every $(\left<r^n\right>,Y) \in \mathcal Q$, a cone of radius $\rho_0$ over the set $Y_r$.
	The quotient of $\dot X$ by $K$ is the space $\bar X$.
	According to Theorem~\ref{res: SC - small cancellation theorem}, $\bar X$ is geodesic proper and $\delta_1$-hyperbolic.
	Moreover, $\bar G$ acts properly co-compactly  by isometries on $\bar X$.
	We now prove that the action of $\bar G$ on $\bar X$ also satisfies Points \ref{enu: SC - induction lemma - elementary subgroups}-\ref{enu: SC - induction lemma - non-commutating elements}.
	Note that the family $\mathcal Q$ satisfies the additional assumptions of Section~\ref{sec: groups without even-torsion}.
	
	\begin{lemm}
	\label{res: SC - induction - elementary subgroups}
		Every elementary subgroup of $\bar G$ is cyclic either infinite or with order dividing $n$.
	\end{lemm}
	
	\begin{proof}
		By Proposition~\ref{res : SC - every elementary subgroup is cyclic}, every elementary subgroup of $\bar G$ is cyclic.
		Let $\bar F$ be a finite subgroup of $\bar G$.
		According to Proposition~\ref{res: SC  - lifting elliptic subgroups}, $\bar F$ is either isomorphic to a finite subgroup of $G$ or contained in some $\stab{Y_r}/\left<r^n\right>$ where $r \in R$.
		By Assumption~\ref{enu: SC - induction lemma - elementary subgroups} every finite subgroup of $G$ has order dividing $n$.
		On the other hand $\stab{Y_r}/\left<r^n\right>$ is isomorphic to $\Z/n\Z$.
		Therefore in both cases the order of $\bar F$ divides $n$.
	\end{proof}
	
	\begin{lemm}
	\label{res: SC - induction - estimating parameters}
		The constant $A(\bar G, \bar X)$ is bounded above by $2\pi \sinh(10^4\delta_1)$ whereas $\rinj {\bar G}{\bar X}$ is bounded below by $20\sqrt {\rho_0\delta_1/n_1}$.
	\end{lemm}
	
	\begin{proof}
		We start with the upper bound of $A(\bar G, \bar X)$.
		According to Proposition~\ref{res: SC - estimate A for the quotient}, $A(\bar G, \bar X) \leq A(G, \lambda  X) + \pi \sinh(10^4 \delta_1)$.
		However the inequality~(\ref{eqn: induction - Delta}) gives
		\begin{displaymath}
			 A(G, \lambda  X) =  \lambda A(G,  X) \leq 2\lambda  \pi \sinh(10^4\delta_1) \leq \pi \sinh(10^4\delta_1).
		\end{displaymath}
		Thus $A(\bar G, \bar X)$ is bounded above by $2\pi \sinh(10^4\delta_1)$.
		We now focus on the injectivity radius of $\bar G$.
		Let $g$ be a hyperbolic isometry of $G$ which does not belong to the stabilizer of $Y_r$ where $r \in R$.
		Its asymptotic translation length in $\lambda  X$ is larger than $400\lambda \delta_1$.
		Proposition~\ref{res: SC - lower bound injectivity radius} combined with (\ref{eqn: induction - rinj}) gives
		\begin{displaymath}
			\rinj {\bar G}{\bar X} \geq \min \left\{ \frac{100\lambda \rho_0\delta_1}{\pi \sinh \rho_0},  \delta_1 \right\} =  \frac{100\lambda \rho_0\delta_1}{\pi \sinh \rho_0} = 20\sqrt {\frac{\rho_0\delta_1}{n_1}}. \qedhere
		\end{displaymath}
	\end{proof}

	\begin{lemm}
	\label{res: SC - induction - non elementary}
		The group $\bar G$ is non-elementary.
		Moreover, there exist $\bar u, \bar v \in \bar G$ and $\bar x \in \bar X$ such that 
		\begin{displaymath}
			0 <\dist {\bar u\bar v\bar u^{-1}\bar v^{-1}\bar x}{\bar x} \leq \rho_0/20.
		\end{displaymath}
	\end{lemm}
	
	\begin{proof}
		We denote by $\bar u$ and $\bar v$ the images of $u$ and $v$ (see Point~\ref{enu: SC - induction lemma - non-commutating elements}) in $\bar G$ and by $\bar x$ the image of $x$ in $\bar X$.
		In the rescaled space $\lambda  X$ we have
		\begin{displaymath}
			0 <\dist {uvu^{-1}v^{-1}x}x \leq \lambda  \rho_0/20 < \rho_0/20.
		\end{displaymath}
		The conclusion follows from Proposition~\ref{res : SC - non elementary quotient}.
	\end{proof}

	\paragraph{}Lemmas~\ref{res: SC - induction - elementary subgroups}, \ref{res: SC - induction - estimating parameters}, \ref{res: SC - induction - non elementary} prove that $(\bar G, \bar X)$ satisfies the induction hypotheses for exponent $n$.
	The only remaining fact to prove concerns the comparison of the translation lengths.
	
	\begin{lemm}
		For every $g \in G$, we have 
		\begin{displaymath}
			\len[stable, espace= \bar X]{\bar g} \leq \frac 1{\sqrt {n_1}} \left(\frac {\pi \sinh \rho_0}{5\sqrt{\rho_0\delta_1}}\right)\len[stable, espace= X]g.
		\end{displaymath}
	\end{lemm}

	\begin{proof}
		Let $g \in G$.
		The asymptotic translation length of $g$ in the rescaled space $\lambda X$ is $\len[stable, espace= \lambda X]g = \lambda \len[stable, espace= X]g$.
		On the other hand the map $\lambda  X \rightarrow \bar X$ shortens the distances, thus $\len[stable, espace= \bar X]{\bar g}\leq  \lambda \len[stable, espace= X]g$.
	\end{proof}

	This last lemma completes the proof of Proposition~\ref{res: SC - induction lemma}.
\end{proof}

\begin{theo}
\label{res : SC - periodic quotient}
	Let $G$ be a non-elementary torsion-free hyperbolic group.
	There exists a critical exponent $n_0$ such that for every odd integer $n \geq n_0$, the quotient $G/G^n$ is infinite.
\end{theo}

\begin{proof}
	We are actually going to prove that $G/G^n$ is not finitely presented and therefore infinite.
	The main ideas are the followings.
	Using Proposition~\ref{res: SC - induction lemma} we construct by induction a sequence of groups $G_0 \rightarrow G_1\rightarrow G_2 \rightarrow \dots$ where $G_{k+1}$ is obtained from $G_k$ by adding new relations of the form $r^n$.
	Then we prove that the direct limit of these groups is exactly $G/G^n$.
	The group $G/G^n$ cannot be finitely presented otherwise the sequence $(G_k)$ should be ultimately constant. 
	The parameters $\rho_0$, $\delta_1$ and $n_0$ are the one given by Proposition~\ref{res: SC - induction lemma}.
		
	\paragraph{Critical exponent.}
	The group $G$ is torsion-free non-elementary and hyperbolic. 
	Let $X$ be the Cayley graph of $G$ with respect to some finite generating set $S$ of $G$.
	It is a proper geodesic $\delta$-hyperbolic space for some $\delta$ depending on $S$.
	Moreover, $A(G,X)$ is finite and there exist two elements $u,v \in G$ which do not commute. 
	In particular if $x$ is the vertex of $X$ corresponding to the identity then $\dist {uvu^{-1}v^{-1}x}x >0$.
	By rescaling if necessary the space $X$ we can assume the followings
	\begin{enumerate}
		\item $\delta \leq \delta_1$,
		\item $A(G,X) \leq 2\pi \sinh (10^4\delta_1)$,
		\item $0 < \dist {uvu^{-1}v^{-1}x}x \leq \rho_0/20$.
	\end{enumerate}
	According to Proposition~\ref{res: positive injectivity radius}, $\rinj GX >0$.
	Therefore, there exists an integer $n_1 \geq n_0$ such that $\rinj GX \geq 20\sqrt {\rho_0\delta_1/n_1}$ and the constant $c$ defined below is less than $1$.
	\begin{displaymath}
		c = \frac 1{\sqrt{n_1}} \left(\frac {\pi \sinh \rho_0}{5\sqrt{\rho_0\delta_1}}\right).
	\end{displaymath}
	From now on we fix an odd integer $n \geq n_1$.
	
	\paragraph{Initialization.} We put $G_0 = G$ and $X_0 = X$.
	In particular $(G_0,X_0)$ satisfies the induction hypotheses for exponent $n$.
	
	\paragraph{Induction.} We assume that we already constructed the group $G_k$ and the space $X_k$ such that $(G_k,X_k)$ satisfies the induction hypotheses for exponent $n$.
	We denote by $R_k$ the set of hyperbolic elements $r$ of $G_k$ which are not a proper power and such that $\len[espace = X_k] r \leq 1000\delta_1$.
	Let $K_k$ be the normal subgroup of $G_k$ generated by $\{r^n, r \in R_k\}$ and $G_{k+1}$ the quotient of $G_k$ by $K_k$.
	By Proposition~\ref{res: SC - induction lemma}, there exists a metric space $X_{k+1}$ such that $(G_{k+1},X_{k+1})$ satisfies the induction hypotheses for exponent $n$.
	Moreover, for every $g \in G_k$, if we still denote by $g$ its image in $G_{k+1}$ we have $\len[stable, espace = X_{k+1}] g \leq c \len[stable, espace = X_k] g$.
	
	\paragraph{Direct limit.} We denote by $G_\infty$ the direct limit of the sequence $(G_k)$.
	We claim that $G_\infty$ is isomorphic to $G/G^n$.
	\paragraph{}
	Let $g$ be an element of $G$.
	To shorten the notation we will still denote by $g$ its image in $G/G^n$, $G_k$ or $G_\infty$.
	It follows from the construction of the sequence $(G_k)$ that for every $k \in \N$ we have $\len[stable, espace= X_k] g \leq c^k\len[stable, espace=X] g$.
	Recall that $c <1$.
	There exists an integer $k$ such that $\len[stable, espace= X_k] g < 20\sqrt {\rho_0\delta_1/n_1}$.
	The injectivity radius of $G_k$ on $X_k$ is at least $20\sqrt {\rho_0\delta_1/n_1}$.
	As an element of $G_k$ the isometry $g$ has finite order.
	However the order of every finite subgroup of $G_k$ divides $n$ hence $g^n$ is trivial in $G_k$ and thus in $G_\infty$.
	Hence $G^n$ is contained in the kernel of the map $G \twoheadrightarrow G_\infty$.
	
	\paragraph{} On the other hand, for every $k \in \N$, the kernel of the map $G_k \rightarrow G_{k+1}$ is generated by the $n$-th power of some elements of $G_k$.
	It follows that the kernel of the map $G \twoheadrightarrow G_\infty$ is exactly $G^n$.
	Consequently, $G_\infty$ and $G/G^n$ are isomorphic.
	
	\paragraph{Conclusion.} Assume now that $G/G^n$ is finitely presented.
	Let $\left<S|R\right>$ be a finite presentation of $G/G^n$.
	We still denote by $S$ a preimage of $S$ in $G$.
	The set $S$ does not necessarily generate $G$, however since $G$ is also finitely generated there exists $k_0 \in \N$ such that for all $k \geq k_0$, $S$ generates $G_k$.
	On the other hand $R$ is finite, therefore there exists $k \geq k_0$ such that all the relations in $R$ are satisfied in the group $G_k$.
	Hence $G_k$ and $G/G^n$ are isomorphic.
	In particular every element of $G_k$ has finite order.
	This contradicts the fact that $G_k$ is a non-elementary hyperbolic group.
	Consequently, $G/G^n$ is not finitely presented and thus not finite.
\end{proof}

\subsection{A few words about Gromov's monster group}
\label{sec: gromov monster}

\paragraph{}
Iterated small cancellation theory can be used to produce many examples of groups with pathological properties,
The Gromov monster is such an example \cite{Gro03}.
It is built in such a way that its Cayley graph coarsely contains an expander. 
This is an obstruction for the group to coarsely embed in a Hilbert space.
In particular, it cannot satisfy the Baum-Connes conjecture with coefficients.

\paragraph{}
The global strategy is the following.
One starts with a group $G_0$ satisfying Kazdhan's property (T) and an infinite family of graphs $(\Theta_k)$.
Then, one constructs a sequence of groups
\begin{displaymath}
	G_0 \twoheadrightarrow G_1 \twoheadrightarrow \dots \twoheadrightarrow G_k \twoheadrightarrow G_{k+1} \twoheadrightarrow \dots,
\end{displaymath}
where $G_{k+1}$ is a small cancellation quotient of $G_k$ in which $\Theta_k$ is quasi-isometrically embedded.
Moreover the previously embedded graphs $\Theta_\ell$ with $\ell<k$ are still quasi-isometrically embedded in $G_{k+1}$ but with different parameters.
It uses an extension of graphical small cancellation as sketched in the introduction (Example~\ref{exa: small cancellation with graphs}).
The Gromov monster group is the direct limit of the groups $G_k$.

\paragraph{}
This construction involves many different tools such as small cancellation theory, probability and harmonic analysis, arithmetics,...
In this section we only give a few keys to understand the part relying on small cancellation.
In particular we will not explained how the graphs $(\Theta_k)$ need to be chosen and labelled so that they can be used to iterate small cancellation.
For the rest of the proof we refer the reader to Gromov's original paper \cite{Gro03} or \cite{ArzDel08}.
In this section we follow with little adaptations the presentation given by G.~Arzhantseva and T.~Delzant in \cite{ArzDel08}.

\paragraph{}We begin by recalling some additional facts about hyperbolic geometry.
\begin{defi}
	\label{def: GM - quasi-isometry}
	Let $L \geq 0$, $k \geq 1$ and $l\geq 0$.
	A map $f : X_1 \rightarrow X_2$ between two metric spaces is an \emph{$L$-locally $(k,l)$-quasi-isometric embedding} if for every $x,x' \in X_1$ such that $\dist x{x'}\leq L$,
	\begin{displaymath}
		k^{-1}\dist x{x'} - l \leq \dist{f(x)}{f(x')} \leq k \dist x{x'} +l.
	\end{displaymath}
	If $L = + \infty$ we say that $f$ is a (globally) \emph{$(k,l)$-quasi-isometric embedding}.
\end{defi}

In particular a (locally) quasi-isometric embedding of an interval of $\R$ is a (local) quasi-geodesic.
The following proposition is a consequence of the stability of quasi-geodesics.

\begin{prop}
\label{res: GM - stability quasi-isometry}
	Let $\delta \geq 0$.
	Let $k \geq 1$ and $l \geq 0$.
	There exist positive constants $L$ and $\alpha$ with the following property.
	Let $X_1$ and $X_2$ be two geodesic metric spaces such that $X_2$ is $\delta$-hyperbolic.
	If $f : X_1 \rightarrow  X_2$ is an $L$-locally $(k,l)$-quasi-isometric embedding then it is a (globally) $(2k,l)$-quasi-isometric embedding and its image $f(X_1)$ in $X_2$ is $\alpha$-quasi-convex.
\end{prop}

\rem In this proposition the constants $L$ and $\alpha$ only depend on $k$, $\delta$ and $l$.
Using a rescaling argument we obtain that the best possible values for $L$ and $\alpha$ satisfy the following property.
For every $\lambda >1$, $L(k,\lambda l, \lambda \delta) = \lambda L(k,l,\delta)$ and $\alpha(k,\lambda l, \lambda \delta) = \lambda \alpha(k,l,\delta)$.
In particular  $\alpha$ tends to 0 as $\delta$ and $l$ approach 0.

\paragraph{}
We present now one step of the construction, namely how the quotient $G_{k+1}$ is built from $G_k$.
This operation should be thought as a partial analogue for the Gromov monster group of Proposition~\ref{res: SC - induction lemma}.
From now on, $G$ is a torsion-free non-elementary hyperbolic group and $X$ its Cayley graph with respect to some finite generating set $S$ of $G$.
We assume here that $S \cap S^{-1} = \emptyset$.
Let $\Theta$ be a finite connected graph with no vertex of degree 1.
A labeling of $\Theta$ by $S$ is a map $m$ which assigns to each edge of $\Theta$ an orientation and a letter of $S$.
If we fix a base point in $\Theta$ then $m$ induces a homomorphism $m_*$ from $\pi_1(\Theta)$ to $G$ defined as follows.
Given a simplicial loop $\gamma$ in $\Theta$, its image under $m_*$ is the element of $G$ represented by the word over the alphabet $S \cup S^{-1}$ that labels $\gamma$.
Let $T$ be the universal cover of $\Theta$.
The labeling $m$ also induces a $\pi_1(\Theta)$-equivariant simplicial map from $T$ to $X$ that we denote by $f$.
We write $Y$ for the image under $f$ of $T$ in $X$ and $H$ for the image under $m_*$ of $\pi_1(\Theta)$ in $G$.

\paragraph{}
The goal is to understand the quotient $\bar G = G / \lnormal H \rnormal$ where $\lnormal H \rnormal$ is the normal subgroup generated by $H$.
In particular we want to prove that it is a hyperbolic group.
To that end we introduce the following family.
\begin{displaymath}
	\mathcal Q = \set{\left(gHg^{-1},gY\right)\fantomB}{g \in G}.
\end{displaymath}
Note that if $g$ belongs to $H$ then $(gY,gHg^{-1})$ and $(Y,H)$ defines the same pair. 
However since we consider $\mathcal Q$ as a \emph{set} of pairs, it will count only once.
We define then two small cancellation parameters which respectively play the role of the length of largest piece and the length of the smallest relation in the usual small cancellation theory.
\begin{eqnarray*}
	\Delta'(\mathcal Q)  & = & \sup_{g \in G\setminus H} \diam\left( Y^{+12\delta} \cap gY^{+12\delta}\right) \\
	T(\mathcal Q) & = & \inf \set{\len h}{h \in H\setminus\{1\}}
\end{eqnarray*}
Note that the parameter $\Delta'(\mathcal Q)$ is not exactly the same as the one we used before.
Indeed the regular parameter $\Delta(\mathcal Q)$ would be
\begin{displaymath}
	\Delta(\mathcal Q)   =  \sup_{g \in G\setminus(\stab Y\cap\operatorname{Norm}(H))} \diam\left( Y^{+12\delta} \cap gY^{+12\delta}\right),
\end{displaymath}
where $\operatorname{Norm}(H)$ is the normalizer of $H$ in $G$.
In particular $\Delta(\mathcal Q) \leq \Delta'(\mathcal Q)$.
The small cancellation assumptions that will follow will be thus stronger than the one we presented.
As an other consequence we can see that if $\Delta'(\mathcal Q)$ is finite then $H$ is the whole stabilizer of $Y$.
The small cancellation theorem used by G.~Arzhantseva and T.~Delzant in \cite{ArzDel08} can be stated as follows.

\begin{theo}[Compare {\cite[Th. 3.10]{ArzDel08}}]
\label{res: GM - small cancellation}
There exists $\rho_0 \geq 10^{20}\boldsymbol \delta$ such that for all $k \geq 1$, for all $\rho \geq \rho_0$, there exist positive numbers $\delta_2$ and $\Delta_2$ which do not depend on $G$, $\Theta$ or $m$ with the following property.
Let $l \geq 0$.
Let $L = L(k,l,\delta)$ and $\alpha = \alpha(k,l,\delta)$ be the constants given by Proposition~\ref{res: GM - stability quasi-isometry}.
Assume that $f : T \rightarrow X$ is an $L$-locally $(k,l)$-quasi-isometric embedding.
If 
\begin{displaymath}
	\frac \delta{T(\mathcal Q)} \leq \delta_2, \quad
	\frac \alpha{T(\mathcal Q)} \leq 10 \delta_2 \quad \text{and} \quad
	\frac {\Delta'(\mathcal Q)}{T(\mathcal Q)} \leq \Delta_2,
\end{displaymath}
then the followings hold.
\begin{enumerate}
	\item \label{enu: GM - hyperbolic quotient}
	The quotient $\bar G$ is a non-elementary torsion-free hyperbolic group. Moreover the hyperbolicity constant of its Cayley graph with respect to the image of $S$ in $\bar G$ can be bounded by a number only depending on $\rho$, $\delta$, $T(\mathcal Q)$ and $\diam \Theta$.
	\item \label{enu: GM - local embedding}
	The projection $G \twoheadrightarrow \bar G$ induces an embedding from $B(1,R)$ in $\bar G$ where $B(1,R)$ is the ball of $G$ of center $1$ and radius
	\begin{displaymath}
		R =  \frac 1{20}. \frac \rho{\pi\sinh \rho} T(\mathcal Q).
	\end{displaymath}
	\item \label{enu: GM - quasi-isometric embedding of the graph}
	The map $f$ induces a map (that we still denote $f$) from $\Theta$ into the Cayley graph of $\bar G$ satisfying for all $x,x' \in \Theta$
	\begin{displaymath}
		\dist{f(x)}{f(x')} \geq \frac 1{500}.\frac{T(\mathcal Q)}{\diam \Theta}. \frac {\rho}{\pi \sinh \rho} \left(\frac 1{2k}\dist x{x'} - l\right).
	\end{displaymath}
\end{enumerate}
\end{theo}

\begin{proof}[Sketch of proof]
We denote by $\rho_0$, $\delta_0$ and $\Delta_0$ the constants given by Theorem~\ref{res: SC - small cancellation theorem}.
Without loss of generality we can assume that $\rho_0 \geq 10^{20}\boldsymbol \delta$, $\delta_0 \leq 10^{-10}\boldsymbol \delta$ and $\Delta_0 \geq 100 \delta_0$.
Let $k \geq 1$ and $\rho \geq \rho_0$.
We put
\begin{displaymath}
	\delta_2 = \frac{\delta_0}{\pi \sinh \rho} \quad \text{and} \quad \Delta_2 = \frac{\Delta_0 - 28 \delta_0}{\pi \sinh \rho}
\end{displaymath}
Since $f: T \rightarrow$ is a $L$-locally $(k,l)$-quasi-isometric embedding in $X$ it follows from Proposition~\ref{res: GM - stability quasi-isometry} that $f$ is actually a $(2k,l)$-quasi-isometric embedding in $X$ and $Y = f(T)$ is $\alpha$-quasi-convex.

\paragraph{}
We now consider the action of $G$ on the rescaled space $\lambda X$ where $\lambda = \pi \sinh \rho/T(\mathcal Q)$.
In particular $\lambda X$ is $\delta_0$-hyperbolic and for every $h \in H\setminus\{1\}$, $\len[espace=\lambda X] h \geq \pi\sinh\rho$.
Moreover $Y$, viewed as a subset of $\lambda X$ is $10 \delta_0$-quasi-convex.
Note that $Y$ is not necessarily strongly-quasi-convex.
We denote by $Z$ the 12$\delta_0$-neighborhood of $Y$ in $\lambda X$.
By Proposition~\ref{res: strong neighborhood of a quasi-convex}, $Z$ is strongly quasi-convex.
Instead of working with $\mathcal Q$ we consider the family $\mathcal S$ defined as 
\begin{displaymath}
	\mathcal S = \set{\left(gHg^{-1},gZ\right)\fantomB}{g \in G}.
\end{displaymath}
According to Lemma~\ref{res: intersection of thickened quasi-convex}, $\Delta'(\mathcal S)$ is bounded above by $\lambda \Delta'(\mathcal Q) + 28\delta_0$.
Hence $\Delta(\mathcal S) \leq \Delta'(\mathcal S) \leq \Delta_0$.
On the other hand $T(\mathcal S) = \lambda T(\mathcal Q) \geq \pi \sinh \rho$.
Thus $\mathcal S$ satisfies the assumptions of Theorem~\ref{res: SC - small cancellation theorem}.
Let $\dot X$ be the cone-off obtain by attaching for every $g \in G/H$ a cone of radius $\rho$ and base $Z$ on $X$.
Let $\bar X$ be the quotient of $\dot X$ by $\lnormal H\rnormal$.
By Theorem~\ref{res: SC - small cancellation theorem}, $\bar X$ is a proper geodesic hyperbolic space and $\bar G$ acts properly co-compactly by isometries on it.
In particular $\bar G$ is hyperbolic.

\paragraph{}Recall that the embedding $\lambda X \hookrightarrow \dot X$ is 1-Lipschitz.
On the other hand the map $\nu : \dot X \rightarrow \bar X$ induces an isometry from $B(1, \rho/20)$ onto its image (see Proposition~\ref{res: rotation family - nu local isometry}).
It follows that the projection $G \rightarrow \bar G$ induces an embedding from $B(1,R)$ in $\bar G$ where
	\begin{displaymath}
		R = \lambda^{-1} \frac \rho{20} = \frac 1{20} \frac \rho{\pi\sinh \rho} T(\mathcal Q).
	\end{displaymath}
This proves Point~\ref{enu: GM - local embedding}.
Similar techniques provide an estimate of the hyperbolicity constant of $\bar G$ and Point~\ref{enu: GM - quasi-isometric embedding of the graph}.
See \cite[Section 3.4]{ArzDel08} for the details.

\rem By assumption $\rho \geq 10^{20}\boldsymbol \delta$ and
\begin{displaymath}
	\frac \delta{T(\mathcal Q)} \leq \delta_2 = \frac {\delta_0}{\pi\sinh \rho} \leq \frac{10^{-10}\boldsymbol \delta}{\pi \sinh \rho}.
\end{displaymath}
If follows that 
\begin{displaymath}
	R = \frac 1{20} \frac \rho{\pi\sinh \rho} T(\mathcal Q) \geq \frac{10^{200}\delta}{20}.
\end{displaymath}
In particular $R$ is very large compare to the hyperbolicity constant $\delta$ of $X$, the Cayley graph of $G$.
Thus Point~\ref{enu: GM - local embedding} is far from being a vacuous observation.

\paragraph{}We now prove that $\bar G$ is torsion-free and non-elementary.
By assumption $G$ is torsion-free.
It follows from Proposition~\ref{res: SC  - lifting elliptic subgroups} that $\bar G$ is torsion free if and only if $\stab Y/H$ is.
The small cancellation assumption involving $\Delta'(\mathcal Q)$ is much stronger than the regular assumption of Theorem~\ref{res: SC - small cancellation theorem}.
As we already noticed, since $\Delta'(\mathcal Q)$ is finite the subgroup $H$ is exactly $\stab Y$.
Hence $\bar G$ is torsion-free.
In particular every elementary subgroup of $\bar G$ is cyclic.
Assume now that $\bar G$ is elementary.
Then $\bar G$ would be cyclic thus abelian. 
In particular, all the generators of $S$ would commute in $\bar G$.
This is not compatible with the embedding of $B(1,R)$ into $\bar G$.
Hence $\bar G$ is not elementary.
\end{proof}

This theorem can be now used to iterate the small cancellation process.
One has in some sense more freedom than for Burnside groups.
In order to uniformly bound the exponent during the induction we had indeed to provide precise estimates of $A(\bar G,\bar X)$, $\rinj {\bar G}{\bar X}$, etc.
This is not needed here.
In particular if the girth of the next graph one wants to embed is very large compare to the previous parameters then it will be easy to satisfy again the assumptions $\delta/T(\mathcal Q)\leq \delta_2$ and  $\alpha/T(\mathcal Q) \leq 10 \delta_2$.
In this way the group that we obtain at the end is lacunary hyperbolic.
The difficulty here is to exhibit the right family of graphs and an appropriate labeling so that one can repeat the process.

% !TEX root = notes.tex

\appendix

\section{Appendix: Cartan-Hadamard Theorem}
\label{sec: cartan hadamard}

\paragraph{}
In \cite{DelGro08} the authors state and exploit a Cartan-Hadamard theorem for orbifolds.
They prove that if an orbifold is locally hyperbolic (with the appropriate quantifiers) then it is developable and its universal cover (in the sense of orbifold) is globally hyperbolic.
As we explained, rotation families provide a means for avoiding orbifolds.
However we still need a ``regular'' Cartan-Hadamard theorem.
This result was given first in Gromov's original paper about hyperbolic groups \cite{Gro87}.
Other proofs using isoperimetric  inequalities can be found in \cite[Chap. 8, Th. 8.1.2]{Bowditch:1991wl} and \cite{Papasoglu:1996uv}.
In this appendix we present an other proof based on topology.
It relies on the ideas used by T.~Delzant and M.~Gromov for the case of orbifolds.

\paragraph{}Let $X$ be a length space.
We say that $X$ is \emph{$\sigma$-locally $\delta$-hyperbolic} if every ball of radius $\sigma$ is $\delta$-hyperbolic i.e., for every four points $x$, $y$, $z$ and $t$ contained in such a ball 
\begin{displaymath}
\gro xzt \geq \min \left\{ \gro xyt, \gro yzt\right\} - \delta.
\end{displaymath}
If its fundamental group of is normally generated by free homotopy classes of loops of diameter less than $r$ then $X$ is said to be \emph{$r$-simply connected}.
In particular, if $X$ is simply connected, it is $r$-simply connected for every $r>0$.

\begin{theo}[Cartan-Hadamard Theorem]
\label{res: cartan hadamard}
Let $\delta \geq 0$ and  $\sigma > 10^7\delta$.
Let $X$ be a length space.
If $X$ is $\sigma$-locally $\delta$-hyperbolic and $10^{-5}\sigma$-simply connected then $X$ is (globally) $300\delta$-hyperbolic.
\end{theo}

\paragraph{General strategy.}
Let us first give the main ideas of the proof.
We are going to demonstrate a version of stability of local quasi-geodesics in $X$ (Proposition~\ref{res: CH - stability of quasi-geodesics}).
The global hyperbolicity will follow from it (Proposition~\ref{res: CH - thin triangles}).
To that end we endow the set of paths of $X$ with a binary relation (Definition~\ref{def: CH - following}). 
Two related paths notably fellow travel and have the same endpoints.
Restricted to the set of local quasi-geodesics, this relation turns out to be an equivalence relation (Proposition~\ref{res: CH - transitivity}).
After fixing a base point $x_0$ in $X$, we look at the set $\tilde X$ of equivalence classes of local quasi-geodesics starting at $x_0$.
It comes with a natural map $p : \tilde X \rightarrow X$ which sends each path to its terminal point.
We show that, endowed with the appropriate topology, $\tilde X$ is a path-connected cover of $X$ (Lemma~\ref{res: CH - path connected} and Proposition~\ref{res: CH - cover}).
If $X$ was simply connected this would force $p$ to be a bijection. 
Actually it remains true if $X$ is just $10^{-5}\sigma$-simply connected.
It implies that two paths starting at $x_0$ with the same terminal points are equivalent and therefore at bounded Hausdorff distance.

\paragraph{} 
Before starting the proof we recall a well-known fact about hyperbolic space.
If $X$ is a geodesic hyperbolic space the distance function is quasi-convex \cite[Chap. 10, Cor. 5.3]{CooDelPap90}.
The next proposition gives a precise statement of the analog property for two quasi-geodesics in a hyperbolic length space.
	
	\begin{prop}
	\label{res: quasi-convexity distance}
		Let $X$ be a $\delta$-hyperbolic length space.
		Let $\gamma : [a,b]\rightarrow X$ and $\gamma : [a',b']\rightarrow X$ be respectively $(1,l)$- and  $(1,l')$-quasi-geodesics.
		Let  $f:[0,1] \rightarrow \R$ be the function defined by 
		\begin{displaymath}
			f(\theta) = \dist{\gamma\left(\fantomB (1-\theta)a + \theta b\right)}{\gamma'\left(\fantomB  (1-\theta)a' +\theta b'\right)}
		\end{displaymath}
		For every $\theta \in [0,1]$,
		\begin{displaymath}
			f(\theta) \leq \theta f(0)+ (1-\theta)f(1)+ 2l + 2l' + 8\delta
		\end{displaymath}
	\end{prop}

	\begin{proof}
		Let us first examine the case where $\gamma$ and $\gamma'$ start at the same point.
		We put $x = \gamma(a) = \gamma'(a')$.
		Let $y$ and $y'$) be the respective endpoints of $\gamma$ and $\gamma'$.
		Let $\theta \in [0,1]$.
		The points $s$ and $s'$ respectively stand for $\gamma\left((1-\theta) a + \theta b\right)$ and $\gamma'\left((1-\theta) a' + \theta b'\right)$.
		Since $\gamma$ is a $(1,l)$-quasi-geodesic $s$ satisfies $\gro xys \leq l/2$. 
		Moreover,
		\begin{displaymath}
			\theta\dist xy -l \leq \dist xs \leq \theta\dist xy +l.
		\end{displaymath}
		The same kind of properties holds for $s'$.
		By Lemma~\ref{res: metric inequalities}-\ref{enu: metric inequalities - comparison tripod} we have
		\begin{displaymath}
			\dist s{s'}
			\leq \max \left\{ \dist{\fantomB\dist xs}{\dist x{s'}} + \max\left\{l,l'\right\}, \dist  xs + \dist x{s'} - 2 \gro y{y'}x\right\} + 4 \delta.
		\end{displaymath}
		However using the triangle inequality we have
		\begin{displaymath}
			\dist{\fantomB\dist xs}{\dist x{s'}} \leq \theta \dist {\fantomB\dist xy}{\dist x{y'}}  +l + l' \leq \theta\dist y{y'} +l +l'.
		\end{displaymath}
		On the other hand,
		\begin{eqnarray*}
			\dist  xs + \dist x{s'} - 2 \gro y{y'}x & \leq &  \theta\left(\fantomB\dist xy +\dist x{y'}\right)  - 2 \gro y{y'}x + l + l' \\
			& \leq &  \theta\dist y{y'} - 2(1-\theta) \gro y{y'}x + l + l'.
		\end{eqnarray*}
		Combining these inequalities we obtain that $f(\theta) = \dist s{s'}$ is bounded above by
		\begin{displaymath}
			\theta \dist y{y'} + l+l'+  \max\left\{l,l'\right\} +4\delta
			= \theta f(1) + l+l'+  \max\left\{l,l'\right\} +4\delta.
		\end{displaymath}

		\paragraph{}
		Let us consider now the case where $\gamma$ and $\gamma'$ have different extremities.
		We denote by $x$ and $y$ (\resp $x'$ and $y'$) the endpoints of $\gamma$ (\resp $\gamma'$) such that $f(0) = \dist x{x'}$ and $f(1) = \dist y{y'}$.
		Let $\eta >0$.
		There exists a $(1,\eta)$-quasi-geodesic $\mu : [c,d] \rightarrow X$ such that $\mu(c) = x$ and $\mu(d) = y'$.
		Let $\theta \in [0,1]$.
		The points $s$, $t$ and $s'$ respectively stand for $\gamma\left((1-\theta) a +\theta b\right)$, $\mu\left( (1-\theta)c +\theta d \right)$ and $\gamma'\left( (1-\theta)a' + \theta b'\right)$.
		Hence $f(\theta) = \dist s{s'} \leq \dist st + \dist t{s'}$.
		It follows from the previous case that 
		\begin{eqnarray*}
			\dist st & \leq & \theta\dist y{y'} + l+\eta + \max\left\{l,\eta\right\} +4\delta, \\
			\dist t{s'} & \leq & (1-\theta)\dist x{x'} +\eta+l'+  \max\left\{\eta,l'\right\} +4\delta.
		\end{eqnarray*}
		Thus
		\begin{displaymath}
			f(\theta) \leq (1-\theta) f(0)+ \theta f(1)+ l + l' +2\eta +\max\left\{l,\eta\right\}+ \max\left\{\eta,l'\right\}+ 8\delta. 
		\end{displaymath}
		This last inequality holds for every $\eta>0$ which completes the proof.
	\end{proof}
	
\paragraph{} For the rest of the appendix we fix $\delta \geq 0$ and $\sigma > 10^7\delta$.
We denote by $X$ a $\sigma$-locally $\delta$-hyperbolic and $10^{-5}\sigma$-simply connected length space.

\subsection{Following paths. Definition and first properties.}

We recall that the paths that we consider are rectifiable and parametrized by arc length.
The length of a path $\gamma$ is denoted by $L(\gamma)$.
\begin{defi}
\label{def: CH - following}
	Let $D\geq 0$.
	Let $\gamma : \intval ab \rightarrow X$ and $\gamma' : \intval{a'}{b'} \rightarrow X$ be two paths of $X$.
	We say that $\gamma$ \emph{$D$-follows} $\gamma'$ if there exists a non decreasing map $\theta : \intval ab \rightarrow \intval{a'}{b'}$ satisfying the following properties.
	\begin{enumerate}
		\item $\theta(a)=a'$ and $\theta(b) = b'$,
		\item for every $s \in \intval ab$, $\dist {\gamma(s)}{\gamma'\circ \theta(s)} < D$,
		\item for every $s,t \in \intval ab$, if $\dist st \leq \sigma/4$ then $\dist {\theta(s)}{\theta(t)} \leq \sigma/3$.	
	\end{enumerate}
	If $D = 10^{-2}\sigma$ we simply say that $\gamma$ \emph{follows} $\gamma'$.
\end{defi}

\rems
We keep the notations of the definition and assume that $\gamma$ $D$-follows $\gamma'$.
Let $s,t \in \intval ab$ such that $s\leq t$. The restriction of $\gamma$ to $\intval st$ also $D$-follows the restriction of $\gamma'$ to $\intval {\theta(s)}{\theta(t)}$. 
If $\gamma'$ is a $\sigma/3$-locally $(1,l')$-quasi-geodesic path then for all $s,t \in \intval ab$ such that $\dist st \leq \sigma/4$,
\begin{eqnarray*}
	\dist{\theta(s)}{\theta(t)} 
	& \leq & \dist {\gamma'\circ\theta(s)}{\gamma'\circ\theta(t)} +l' \\
	& \leq & \dist{\gamma'\circ\theta(s)}{\gamma(s)} + \dist{\gamma(s)}{\gamma(t)} + \dist {\gamma(t)}{\gamma'\circ\theta(t)} +l' \\
	& \leq & \dist {s}{t} +2D + l'
\end{eqnarray*}
Applying this fact twice we obtain that for every $s,t \in \intval ab$ such that $\dist st \leq \sigma/2$ we have $\dist{\theta(s)}{\theta(t)} \leq  \dist {s}{t} +4D + 2l'$.
We will use this argument very often.

\paragraph{}Every path $D$-follows itself.
More generally let $\gamma : \intval ab \rightarrow X$ be a path and $x$, $y$ be two points of $X$ such that $\max\left\{\dist x{\gamma(a)}, \dist y{\gamma(b)}\right\} < D$.
There exists a path joining $x$ to $y$ which $D$-follows $\gamma$.
This path can be obtained as a concatenation of $\gamma_x$, $\gamma$ and $\gamma_y$ where $\gamma_x$ (\resp $\gamma_y$) is a $(1,l)$-quasi-geodesic joining $x$ to $\gamma(a)$ (\resp $\gamma(b)$ to $y$) for a sufficiently small $l$. 

\paragraph{}The next lemmas explain how this notion behaves with respect to the concatenation and extension of paths.

\begin{lemm}
\label{res: CH - concatenation of paths}
	Let $l' \in (0,10^{-5}\sigma)$ and $D \in (0,10^{-2}\sigma)$.
	Let $\gamma : \intval ac \rightarrow X$ be a path and $\gamma' : \intval{a'}{c'} \rightarrow X$ be a $\sigma/3$-local $(1,l')$-quasi-geodesic.
	Let $b \in \intval ac$ and $b' \in \intval {a'}{c'}$.
	We assume that the restriction of $\gamma$ to $\intval ab$ (\resp $\intval bc$)  $D$-follows the restriction of $\gamma'$ to $\intval {a'}{b'}$ (\resp $\intval {b'}{c'}$).
	Then $\gamma$ $D$-follows $\gamma'$ as well.
\end{lemm}

\begin{proof}
	The restriction of $\gamma$ to $\intval ab$ $D$-follows the one of $\gamma'$ to $\intval {a'}{b'}$.
	By definition, there exists a map $\theta_- : \intval ab \rightarrow \intval {a'}{b'}$ satisfying the axioms of Definition~\ref{def: CH - following}.
	In the same way we have a map $\theta_+ : \intval bc \rightarrow \intval{b'}{c'}$.
	We define $\theta : \intval ac \rightarrow \intval {a'}{c'}$ such that its restriction to $\intval ab$ (\resp $\intval bc$) is $\theta_-$ (\resp $\theta_+$).
	Let $s,t \in \intval ac$ such that $\dist st \leq \sigma/4$.
	If $s$ ant $t$ both belong to $\intval ab$ (\resp $\intval bc$), then it follows from the properties of $\theta_-$ and $\theta_+$ that $\dist{\theta(s)}{\theta(t)}\leq \sigma/3$.
	Assume now that $s \in \intval ab$ and $t \in \intval {a'}{b'}$.
	In particular, $\dist sb$ and $\dist bt$ are at most $\sigma/4$.
	Since $\gamma'$ is a $\sigma/3$-local $(1,l')$-quasi-geodesic, we obtain
	\begin{displaymath}
		\dist {\theta(s)}{\theta(t)} =   \dist {\theta_-(s)}{\theta_-(b)} +  \dist {\theta_+(b)}{\theta_+(t)} \leq  \dist st + 4D + 2l' \leq \frac \sigma3.
	\end{displaymath}
	Hence for all $s,t \in \intval ac$ if $\dist st \leq \sigma/4$ then $\dist{\theta(s)}{\theta(t)}\leq \sigma/3$.
	The other axioms of Definition~\ref{def: CH - following} can be easily checked, hence $\gamma$ $D$-follows $\gamma'$.
\end{proof}

\begin{lemm}
\label{res: CH - extended following path}
	Let $l' \in (0,10^{-5}\sigma)$ and $D \in (0,10^{-2}\sigma)$.
	Let $\gamma : \intval ab \rightarrow X$ be a path and ${\gamma': \intval{a'}{c'} \rightarrow X}$ be a $\sigma/3$-locally $(1,l')$-quasi-geodesic.
	Let $b' \in \intval {a'}{c'}$ such that $\dist {b'}{c'} \leq \sigma/12-2D-l'$.
	Moreover, we assume that $\dist {\gamma(b)}{\gamma'(c')} < D$.
	If $\gamma$ $D$-follows the restriction of $\gamma'$ to $\intval {a'}{b'}$ then it $D$-follows $\gamma'$ as well.
\end{lemm}

\begin{proof}
	The path $\gamma$ $D$-follows the restriction of $\gamma'$ to $\intval {a'}{c'}$.
	This provides a map $\theta : \intval ab \rightarrow \intval {a'}{c'}$ satisfying the axioms of Definition~\ref{def: CH - following}.
	We define a map $\mu : \intval ab \rightarrow \intval {a'}{b'}$ by $\mu(s) = \theta(s)$ for every $s \in [a,b)$ and $\mu(b) = b'$.
	Let $s,t \in \intval ac$, $s\leq t$ such that $\dist st \leq \sigma/4$.
	If $t \neq b$, the fact follows from the properties of $\theta$ that $\dist{\mu(s)}{\mu(t)}\leq \sigma/3$.
	Assume now that $t=b$.
	Since $\gamma'$ is a $\sigma/3$-local $(1,l')$-quasi-geodesic, we obtain
	\begin{displaymath}
		\dist{\mu(s)}{\mu(t)}
		=\dist{\theta(s)}{\theta(b)}+\dist{c'}{b'}
		\leq \dist sb +2D +l' +\dist {c'}{b'}
		\leq \frac \sigma3.
	\end{displaymath}
	The other axioms of Definition~\ref{def: CH - following} can be easily checked, hence $\gamma$ $D$-follows $\gamma'$.
\end{proof}

\rem The two previous lemmas lead to the following fact.
Let $l \in (0,10^{-5}\sigma)$ and $D \in (0,10^{-2}\sigma)$.
Let $\gamma : \intval ac \rightarrow X$ be a $\sigma/3$-local $(1,l)$-quasi-geodesic. 
Let $b \in \intval ac$ such that $\dist bc < D$.
Then $\gamma$ $D$-follows its restriction to $\intval ab$ and conversely.

\subsection{Transitivity of the relation}

One important property of the notion of \emph{$D$-following} lies in this fact: assume that $\gamma$ follows $\gamma'$;
if $\gamma$ and $\gamma'$ are quasi-geodesics and have close extremities then $\gamma$ actually $D$-follows $\gamma'$ with  $D$ much smaller than $10^{-2}\sigma$.
This assertion is carefully proved in the next two lemmas. 
It is the crucial point needed to define later a binary relation which is transitive.

\begin{lemm}
\label{res: CH - quasi-geodesics follow each other}
	Let $l,l' \in (0, 10^{-5}\sigma)$.
	Let $\gamma : [a,b]\rightarrow X$ (\resp $\gamma : [a',b']\rightarrow X$) be a $(1,l)$-quasi-geodesic (\resp $(1,l')$-quasi-geodesic).
	We denote by $x$ and $y$ (\resp $x'$ and $y'$) the initial and terminal points of $\gamma$ (\resp $\gamma'$).
	We assume that 
	\begin{displaymath}
		\max \left\{ \dist x{x'}, \dist y{y'}\right\} <\sigma/200.
	\end{displaymath}
	If $\gamma$ and $\gamma'$ are contained in a ball of radius $\sigma$ of $X$ then the path $\gamma$ $D$-follows $\gamma'$ and conversely where 
	\begin{displaymath}
		D = \max \left\{ \dist x{x'}, \dist y{y'}\right\} + 3l +3l' +8\delta.
	\end{displaymath}
\end{lemm}

\begin{proof}
	The paths $\gamma$ and $\gamma'$ are contained in a ball of radius $\sigma$.
	Thus everything works as if we were in a hyperbolic space.
	We denote by $\theta$ the increasing affine function which maps $\intval ab$ onto $\intval {a'}{b'}$.
	In particular, $\theta(a) = a'$ and $\theta(b)=b'$.
	By Proposition~\ref{res: quasi-convexity distance}, for every $t \in \intval ab$ we have
	\begin{displaymath}
		\dist {\gamma(t)}{\gamma'\circ \theta(t)} < \max\left\{\fantomB \dist x{x'}, \dist y{y'}\right\} + 3l + 3l' + 8\delta.
	\end{displaymath}
	On the other hand $\gamma'$ is a $(1,l')$-quasi-geodesic.
	Hence for every $s,t \in \intval ab$ such that $\dist st \leq \sigma/4$, we have $\dist{\theta(s)}{\theta(t)} \leq \dist {s}{t} + 2 D + l' \leq \sigma/3$.
	Therefore $\gamma$ $D$-follows $\gamma'$.
	The statement is symmetric, thus $\gamma'$ also $D$-follows $\gamma$.
\end{proof}

\begin{lemm}
\label{res: CH - local quasi-geodesics follow each other}
	Let $l,l' \in (0,10^{-5}\sigma)$.
	Let $\gamma : [a,b]\rightarrow X$ (\resp $\gamma : [a',b']\rightarrow X$) be a $\sigma/3$-local $(1,l)$-quasi-geodesic (\resp $\sigma/3$-local $(1,l')$-quasi-geodesic).
	We denote by $x$ and $y$ (\resp $x'$ and $y'$) the initial and terminal points of $\gamma$ (\resp $\gamma'$).
	We assume that 
	\begin{displaymath}
		\max \left\{ \dist x{x'}, \dist y{y'}\right\} <  \sigma/200.
	\end{displaymath}
	If $\gamma$ follows $\gamma'$, then $\gamma$ actually $D$-follows $\gamma'$ and conversely where 
	\begin{displaymath}
		D = \max \left\{ \dist x{x'}, \dist y{y'}, l + l' +5\delta \right\} + 3l +3l' +8\delta.
	\end{displaymath}
\end{lemm}

\begin{proof}
	For simplicity of notation we put $\tau = \sigma/25$.
	The path $\gamma$ follows $\gamma'$.
	Hence there exists a non-decreasing function $\theta : \intval ab \rightarrow \intval {a'}{b'}$ satisfying the properties of Definition~\ref{def: CH - following}.
	In particular, for every $t \in \intval ab$ we have $\dist{\gamma(t)}{\gamma'\circ\theta(t)}< 10^{-2}\sigma$.
	The point $\gamma(t)$ is actually much closer to $\gamma'$: we claim that for every $t \in \intval {a+\tau}{b-\tau}$, there exists $t' \in \intval{\theta(t-\tau)}{\theta(t+\tau)}$ such that $\dist {\gamma(t)}{\gamma'(t')} \leq  l + l' +5\delta$.
	Note that $2\tau\leq \sigma/4$.
	Therefore $\gamma$ and $\gamma'$ restricted to $\intval{t-\tau}{t+\tau}$ and $\intval{\theta(t-\tau)}{\theta(t+\tau)}$ are respectively $(1,l)$- and $(1,l')$-quasi-geodesics contained in a ball $B$ of radius $\sigma$.
	In particular, the second one is $(l'+3\delta)$-quasi-convex in the $\delta$-hyperbolic ball $B$.
	Moreover $\dist {\gamma(t)}{\gamma(t\pm\tau)}  \geq \sigma/25 - l$.
	It follows from the hyperbolicity of $B$ that 
	\begin{displaymath}
		\gro{\gamma'\circ\theta(t-\tau)}{\gamma'\circ\theta(t+\tau)}{\gamma(t)} \leq \gro{\gamma(t-\tau)}{\gamma(t+\tau)}{\gamma(t)} + 2\delta \leq l + 2 \delta.
	\end{displaymath}
	By quasi-convexity, the distance between $\gamma(t)$ and the path $\gamma'$ restricted to $\intval{\theta(t-\tau)}{\theta(t+\tau)}$ is less than or equal to $l+l'+5\delta$, which establishes our claim.
	
	\paragraph{}
	Let us now consider the subdivision $a=t_0 \leq \dots \leq t_m = b$ of $\intval ab$ such that for all $i\in \intvald 0{m-2}$, $\dist {t_{i+1}}{t_i} = 2\tau$ and $2\tau \leq \dist {t_m}{t_{m-1}}\leq 4\tau$.
	First we put $t'_0=a'$ and $t'_m = b'$.
	According to our previous claim, for every $i \in \intvald 1{m-1}$ there exists $t'_i \in \intval{\theta(t_i-\tau)}{\theta(t_i+\tau)}$ such that $\dist {\gamma(t_i)}{\gamma'(t'_i)} \leq l + l' +5\delta$.
	It follows that for every $i \in \intvald 0m$;
	\begin{displaymath}
		\dist {\gamma(t_i)}{\gamma'(t'_i)} \leq  \max \left\{\fantomB \dist x{x'}, \dist y{y'}, l + l' +5\delta \right\} <  \sigma/200,
	\end{displaymath}
	Let $i \in \intvald 1{m-2}$.
	By construction $t_i+\tau = t_{i+1}-\tau$.
	Since $\theta$ is non-decreasing we have 
	\begin{displaymath}
		\theta(t_i-\tau)\leq  t'_i \leq \theta(t_i+\tau) = \theta(t_{i+1}-\tau) \leq t'_{i+1} \leq \theta(t_{i+1}+\tau).
	\end{displaymath}
	Moreover, $\dist{(t_{i+1}+\tau)}{(t_i-\tau)}=4\tau \leq \sigma/4$.
	By definition of $\theta$, 
	\begin{displaymath}
		\dist{t'_{i+1}}{t'_i} 
		\leq \dist{\theta(t_{i+1}+\tau)}{\theta(t_i-\tau)} 
		\leq \frac \sigma3.
	\end{displaymath}
	The same facts hold for $i=0$ and $i=m-1$.
	Therefore $a'=t'_0 \leq \dots \leq t'_m = b'$ is a subdivision of $\intval {a'}{b'}$ such that for all $i \in \intvald 0{m-1}$, $\dist {t'_{i+1}}{t'_i} \leq \sigma/3$.
	In particular, $\gamma$ restricted to $\intval{t_i}{t_{i+1}}$ and $\gamma'$ restricted to $\intval{t'_i}{t'_{i+1}}$ are respectively $(1,l)$- and $(1,l')$-quasi-geodesics contained in a ball of radius $\sigma$.
	According to Lemma~\ref{res: CH - quasi-geodesics follow each other}, they $D$-follows each other with
	\begin{displaymath}
		D = \max \left\{ \dist x{x'}, \dist y{y'}, l + l' +5\delta \right\} + 3l +3l' +8\delta.
	\end{displaymath}
	By concatenation (Lemma~\ref{res: CH - concatenation of paths}), $\gamma$ $D$-follows $\gamma'$ and conversely.
\end{proof}

\begin{prop}
\label{res: CH - transitivity}
	Let $l \in (0,10^{-5}\sigma)$.
	For every $i \in  \{1,2,3\}$, let $\gamma_i : \intval{a_i}{b_i} \rightarrow X$ be a $\sigma/3$-local $(1,l)$-quasi-geodesic joining $x_i$ to $y_i$.
	We assume that 
	\begin{enumerate}
		\item $\max\left\{ \dist {x_1}{x_2}, \dist {x_2}{x_3},\dist {y_1}{y_2},\dist{y_2}{y_3}\right\}  < \sigma/400$,
		\item $\gamma_1$ follows $\gamma_2$ and $\gamma_2$ follows $\gamma_3$
	\end{enumerate}
	Then $\gamma_1$ follows $\gamma_3$.
\end{prop}

\begin{proof}
	According to Lemma~\ref{res: CH - local quasi-geodesics follow each other} and Hypothesis (i), $\gamma_1$ actually $\sigma/200$-follows $\gamma_2$.
	It gives us a map $\theta_1 : \intval{a_1}{b_1} \rightarrow \intval {a_2}{b_2}$ satisfying the axioms of Definition~\ref{def: CH - following}.
	In the same way, $\gamma_2$ $\sigma/200$-follows $\gamma_3$ providing a map $\theta_2 : \intval{a_2}{b_2} \rightarrow \intval {a_3}{b_3}$.
	We denote by $\theta : \intval{a_1}{b_1}\rightarrow \intval{a_3}{b_3}$ the composition $\theta = \theta_2\circ \theta_1$.
	By triangle inequality, for all $s \in \intval{a_1}{b_1}$,
	\begin{displaymath}
		 \dist {\gamma_1(s)}{\gamma_3\circ\theta(s)} 
		\leq \dist {\gamma_1(s)}{\gamma_2\circ\theta_1(s)} + \dist {\gamma_2\circ\theta_1(s)} {\gamma_3\circ\theta_2\circ\theta_1(s)} 
		< 10^{-2}\sigma.
	\end{displaymath}
	Let $s,t \in \intval {a_1}{b_1}$ such that $\dist st \leq \sigma/4$.
	Since $\gamma_2$ is a $\sigma/3$-local $(1,l)$-quasi-geodesic, we have	
	\begin{displaymath}
		\dist{\theta_1(s)}{\theta_1(t)} \leq \dist st + \frac\sigma{50}  +l \leq \frac \sigma2.
	\end{displaymath}
	However $\gamma_3$ is also a $\sigma/3$-local $(1,l)$-quasi-geodesic, hence
	\begin{displaymath}
		\dist{\theta(s)}{\theta(t)}
		\leq \dist{\theta_1(s)}{\theta_1(t)}  +\frac \sigma{25} +2l 
		\leq \dist st  + \frac{3\sigma}{50} +3l \leq \frac \sigma3.
	\end{displaymath}
	The other properties of Definition~\ref{def: CH - following} can be easily checked, hence $\gamma_1$ follows $\gamma_3$.
\end{proof}

\begin{defi}
\label{def: CH - equivalence relation}
	Let $l \in (0,10^{-5}\sigma)$.
	Let $\gamma : [a,b]\rightarrow X$ and  $\gamma : [a',b']\rightarrow X$ be two $\sigma/3$-local $(1,l)$-quasi-geodesics. 
	We say that $\gamma$ and $\gamma'$ are \emph{related} if they have the same endpoints and follow each other.
\end{defi}
	
This binary relation is reflexive and symmetric.
By Proposition~\ref{res: CH - transitivity}, it is transitive as well.
Therefore it is an equivalence relation.

\rem The equivalence relation used by T.~Delzant and M.~Gromov in \cite{DelGro08} is the following: two paths are equivalent if they have the same extremities and their Hausdorff distance is at most $200 \delta$.
However this formulation is not sufficient to complete the proofs.
Indeed, to be able to use the local hyperbolicity of $X$ one needs to be sure that every small portion of one path is close to a small portion of the other path, which may not be the case with a global assumption on the Hausdorff distance. 
That is the reason why we introduced the notion of $D$-following paths.

\subsection{Existence of following quasi-geodesics}

\begin{prop}
\label{res: CH - following path}
	Let $l' \in (0,10^{-5}\sigma)$.
	Let $x$, $x'$, $y$ and $y'$ be four points of $X$ such that 
	\begin{displaymath}
		\max \left\{ \dist x{x'}, \dist y{y'} \right\} < \sigma/200.
	\end{displaymath}
	Let $\gamma' : \intval {a'}{b'} \rightarrow X$ be a $\sigma/3$-locally $(1,l')$-quasi-geodesic joining $x'$ to $y'$.
	For every $l >0$, there exists a $\sigma/3$-local $(1,l)$-quasi-geodesic $\gamma$  joining $x$ to $y$ which follows $\gamma'$.
\end{prop}

\begin{proof}
	Let $l \in (0,10^{-10}\sigma)$.
	We denote by $\mathcal C$ the set of paths joining $x$ to $y$ which follow $\gamma'$.
	This set is non-empty.
	We fix a path $\gamma : \intval ab \rightarrow X$ of $\mathcal C$ such that $L(\gamma) \leq \inf_{c \in \mathcal C} L(c) + l$.
	By construction there is a map $\theta : \intval ab \rightarrow \intval{a'}{b'}$ satisfying the axioms of Definition~\ref{def: CH - following}.
	For simplicity of notation, we put $\tau = 3\sigma/100$
	
	\begin{lemm}
	\label{res: CH - following path - lemma 1}
		For every $t \in \intval {a+\tau}{b-\tau}$, the distance $d$ between $\gamma(t)$ and $\gamma'$ restricted to $\intval{\theta(t-\tau)}{\theta(t+\tau)}$ is bounded above by $5l'+l +14\delta$.
	\end{lemm}

	\begin{proof}
		Let $t \in \intval {a+\tau}{b-\tau}$.
		Note that $2\tau \leq \sigma/4$.
		Hence by definition of $\theta$, $\dist{\theta(t+\tau)}{\theta(t-\tau)}$ is at most $\sigma/3$.
		Consequently, the restriction of $\gamma'$ to  the interval $\intval{\theta(t-\tau)}{\theta(t+\tau)}$ is a $(1,l')$-quasi-geodesic lying in a ball $B$ of radius $\sigma$.
		In particular, this path is $(l'+3\delta)$-quasi-convex in $B$.
		One can choose $B$ in such a way that it contains $\gamma(t-\tau)$, $\gamma(t)$ and $\gamma(t+\tau)$.
		This allow us to use projections on a quasi-convex subset (see Lemma~\ref{res: proj quasi-convex}).
		We define four points of $B$ (see Fig.~\ref{fig: following paths}).
		\begin{figure}[htbp]
			\centering
			\includegraphics{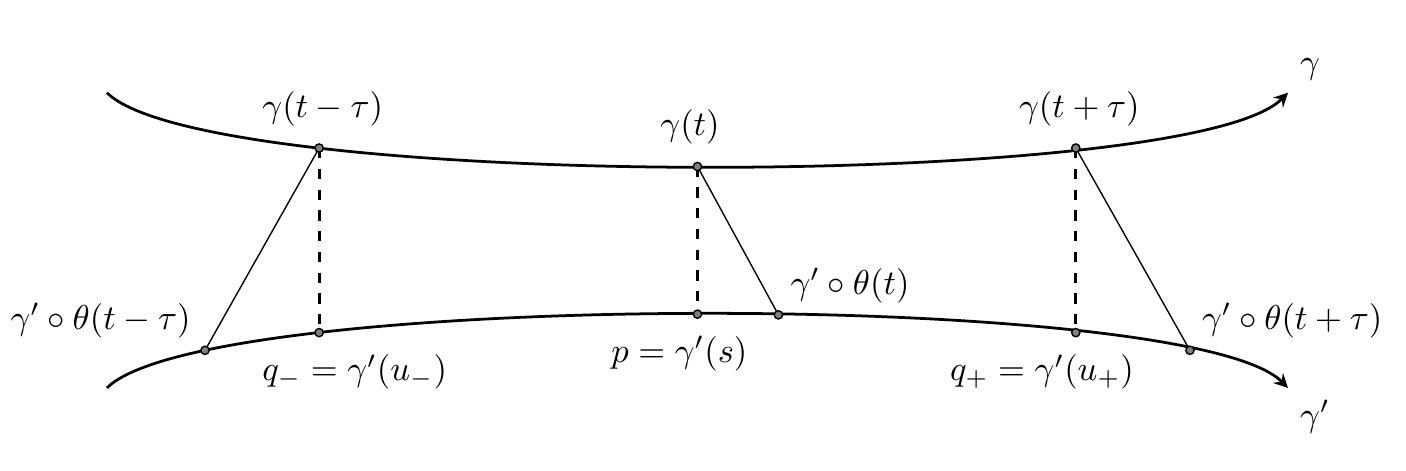}
		\caption{The paths $\gamma$ and $\gamma'$}
		\label{fig: following paths}
		\end{figure}

		\begin{itemize}
			\item $p=\gamma'(s)$ is a projection of $\gamma(t)$ on $\gamma'\left(\intval{\theta(t-\tau)}{\theta(t)}\right)$,
			\item $q_-=\gamma'(u_-)$ is a projection of $\gamma(t-\tau)$ on $\gamma'\left(\intval{\theta(t-\tau)}{s}\right)$,
			\item $q_+=\gamma'(u_+)$ is a projection of $\gamma(t+\tau)$ on $\gamma'\left(\intval{s}{\theta(t+\tau)}\right)$.
		\end{itemize}
		In particular, 
		\begin{displaymath}
			d \leq \dist{\gamma(t)}{p} < 10^{-2} \sigma 
			\text{ and } 
			\dist{\gamma(t\pm\tau)}{q_\pm} < 10^{-2} \sigma.
		\end{displaymath}
		
		\paragraph{Claim 1.} $\dist p{q_-} > \sigma/200$ and $\dist p{q_+} > \sigma/200$.
		Without loss of generality we can assume that $s \leq \theta(t)$.
		 Suppose, contrary to our claim, that $\dist p{q_-} \leq \sigma/200$.
		Hence the length of $\gamma'$ restricted to $\intval s{u_-}$ is at most $\sigma/200+l'$.
		However the distances $\dist {\gamma(t-\tau)}{q_-}$ and $\dist{\gamma(t)}p$ are bounded above by $\sigma/100$.
		Therefore there exists a path $\gamma_0$ joining $\gamma(t-\tau)$ to $\gamma(t)$ whose length is at most $\sigma/50+l'$ and which follows $\gamma'$ restricted to $\intval s{u_-}$.
		On the other hand $\dist{\theta(t-\tau)}{u_-}$ and $\dist {\theta(t)}s$ are bounded above by $\sigma/20 +l'$.
		It follows from Lemma~\ref{res: CH - extended following path} that $\gamma$ restricted to $\intval a{t-\tau}$ and $\gamma$ restricted to $\intval tb$ respectively follow $\gamma'$ restricted to $\intval{a'}{u_-}$ and $\gamma'$ restricted to $\intval s{b'}$.
		 According to Lemma~\ref{res: CH - concatenation of paths}, concatenating $\gamma$ restricted to $\intval a{t-\tau}$, $\gamma_0$ and $\gamma$ restricted to $\intval tb$ gives a path joining $x$ to $y$ which follows $\gamma'$.
		 In particular, it belongs to $\mathcal C$.
		 Nevertheless, its length is bounded above by
		 \begin{displaymath}
		 	L\left(\restriction \gamma{\intval a{t-\tau}}\right) + L\left(\restriction \gamma{\intval tb}\right) +\frac \sigma{50} + l'
			= L(\gamma) - \tau + \frac \sigma{50} +l'
			< L(\gamma) - l,
		 \end{displaymath}
		which contradicts the minimality of $\gamma$.
		With a similar argument, we prove that 
		\begin{displaymath}
			\dist{\gamma'\circ\theta(t)}{q_+} > \sigma/200+l'.
		\end{displaymath}
		However $\gamma'$ restricted to $\intval{\theta(t-\tau)}{\theta(t+\tau)}$ is a $(1,l')$-quasi-geodesic and $s \leq \theta(t)$.
		It follows that $\dist p{q_+} > \sigma/200$.

		\paragraph{Claim 2.} The length of $\gamma'$ restricted to $\intval {u_-}{u_+}$ is bounded above by 
		\begin{displaymath}
			L\left( \restriction \gamma{\intval{t-\tau}{t+\tau}}\right) - \dist{\gamma(t-\tau)}{q_-} - \dist{\gamma(t+\tau)}{q_+} - 2d + 9l' + 28 \delta.
		\end{displaymath}
		The points $p$ and $q_-$ are respective projections of $\gamma(t)$ and $\gamma(t-\tau)$ on the path $\gamma'$ restricted to $\intval{\theta(t-\tau)}{s}$.
		By projection on a quasi-convex, we have
		\begin{displaymath}
			\dist p{q_-} \leq \max\left\{\fantomB\dist{\gamma(t-\tau)}{\gamma(t)} - \dist {\gamma(t-\tau)}{q_-} -\dist{\gamma(t)}p +2\epsilon, \epsilon\right\}
		\end{displaymath}
		where $\epsilon = 2l' + 7\delta$.
		According to our previous claim, we necessarily have
		\begin{displaymath}
		\label{eq: CH - lower bound length - first half}
			\dist{\gamma(t-\tau)}{\gamma(t)} \geq  \dist {\gamma(t-\tau)}{q_-} + \dist{q_-}p + \dist p{\gamma(t)} - 4l' - 14 \delta.
		\end{displaymath}
		In the same way we obtain
		\begin{displaymath}
		\label{eq: CH - lower bound length - second half}
			\dist{\gamma(t+\tau)}{\gamma(t)} \geq  \dist {\gamma(t+\tau)}{q_+} + \dist{q_+}p + \dist p{\gamma(t)} - 4l' - 14 \delta.
		\end{displaymath}
		Combining these two inequalities, we see that the length of $\gamma$ restricted to $\intval{t-\tau}{t+\tau}$ is bounded below by
		\begin{eqnarray*}
			\dist {\gamma(t-\tau)}{q_-} + \dist{q_-}{q_+} +\dist{q_+}{\gamma(t+\tau)} + 2d- 8l' - 28\delta.
		\end{eqnarray*}
		However $\gamma'$ restricted to $\intval{t-\tau}{t+\tau}$ is a $(1,l')$-quasi-geodesic.
		Therefore,
		\begin{displaymath}
			\dist{q_-}{q_+}
			\geq L\left(\restriction{\gamma'}{\intval {u_-}{u_+}}\right) - l'.
		\end{displaymath}
		Combined with the previous inequality we obtain our second claim.
		
		\paragraph{}
		Recall that $\dist{\gamma(t\pm\tau)}{q_\pm} <\sigma/100$.
		Hence there exists a path $\gamma_1$ joining $\gamma(t-\tau)$ to $\gamma(t+\tau)$ which follows the restriction of $\gamma'$ to $\intval {u_-}{u_+}$  and whose length is at most 
		\begin{displaymath}
			\dist {\gamma(t-\tau)}{q_-} + L\left(\restriction{\gamma'}{\intval {u_-}{u_+}}\right) + \dist{q_+}{\gamma(t+\tau)} +l.
		\end{displaymath}
		According to our second claim
		\begin{displaymath}
			L\left(\gamma_1\right) \leq  L\left(\restriction{\gamma}{\intval {t-\tau}{t+\tau}}\right) -2d + 9l' + l+ 28 \delta .
		\end{displaymath}
 		On the other hand, by Lemma~\ref{res: CH - extended following path}, the restriction of $\gamma$ to $\intval a{t-\tau}$ (\resp $\intval {t+\tau}b$) follows the one of $\gamma'$ to $\intval{a'}{u_-}$ (\resp $\intval {u_+}{b'}$).
		By concatenating $\gamma$ restricted to $\intval a{t-\tau}$, $\gamma_1$ and $\gamma$ restricted to $\intval {t+\tau}b$ we obtain a path joining $x$ to $y$ which follows $\gamma'$ (see Lemma~\ref{res: CH - concatenation of paths}). 
		In particular, this is an element of $\mathcal C$.
		Moreover, its length is bounded above by $L(\gamma) - 2d +9l' + l + 28 \delta$.
		By construction, $\gamma$ minimizes up to $l$ the length of the elements of $\mathcal C$, hence $d \leq 5l' + l +14 \delta $.
	\end{proof}
	
	\begin{lemm}
	\label{res: CH - following path lemma 2}
		The path $\gamma$ $D$-follows $\gamma'$ with
		\begin{displaymath}
			D = \max\left\{\dist x{x'}, \dist y{y'},5l'+l+14\delta\right\} + 3l + 3l' + 8\delta < \sigma/100.
		\end{displaymath}
	\end{lemm}

	\begin{proof}
		We consider the subdivision $a = t_0 \leq \dots \leq t_m = b$ of $\intval ab$ such that for every $i \in \intvald 0{m-2}$, $\dist{t_{i+1}}{t_i}=2\tau$ and $2\tau\leq \dist {t_m}{t_{m-1}} \leq 4\tau$. 
		We claim that for every $i \in \intvald 0{m-1}$, the restriction of $\gamma$ to $\intval {t_i}{t_{i+1}}$ is a $(1,l)$-quasi-geodesic.
		According to Lemma~\ref{res: CH - following path - lemma 1}, there exists a subdivision $a'= s_0 \leq \dots \leq s_m = b'$ of $\intval {a'}{b'}$ satisfying the following properties.
		\begin{enumerate}
			\item For every $i \in \intvald 0m$, 
			\begin{displaymath}
				\dist{\gamma(t_i)}{\gamma'(s_i)} \leq \max \left\{\dist x{x'}, \dist y{y'}, 5l' +l+ 14 \delta \right\} < \sigma/200.
			\end{displaymath}
			\item For every $i \in \intvald 0{m-1}$, 
			\begin{displaymath}
				\dist {s_{i+1}}{s_i} \leq \dist {t_{i+1}}{t_i} + \frac\sigma{50} +l' \leq \frac \sigma3.
			\end{displaymath}
		\end{enumerate}
		Let $\eta \in (0, 10^{-10}\sigma)$.
		Let $i \in \intvald 0{m-1}$. 
		There exists a $(1,\eta)$-quasi-geodesic $c_i$ joining $\gamma(t_i)$ to $\gamma(t_{i+1})$.
		Thus $c_i$ and $\gamma'$ restricted to $\intval{t_i}{t_{i+1}}$ are respectively $(1,\eta)$- and $(1,l')$-quasi-geodesics contained in a ball of radius $\sigma$.
		By Lemma~\ref{res: CH - quasi-geodesics follow each other}, they follow each other.
		Consequently, the path $c$, obtained by concatenating the $c_i$'s, follows $\gamma'$ (see Lemma~\ref{res: CH - concatenation of paths}).
		In particular, $c$ belongs to $\mathcal C$.
		Therefore its length is bounded below by $L(\gamma) - l$.
		Since the $c_i$'s are $(1, \eta)$-quasi-geodesics one has
		\begin{displaymath}
			\sum_{i=0}^{m-1} L\left(\restriction{\gamma}{\intval{t_i}{t_{i+1}}}\right) \leq \sum_{i=0}^{m-1} \dist{\gamma(t_{i+1})}{\gamma(t_i)} + m\eta + l.
		\end{displaymath}
		It implies that for every $i \in \intvald 0{m-1}$, the length of $\gamma$ restricted to $\intval{t_i}{t_{i+1}}$ is bounded above by $\dist{\gamma(t_{i+1})}{\gamma(t_i)} + m\eta + l$.
		This inequality holds for every $\eta>0$.
		Hence the restriction of $\gamma$ to $\intval {t_i}{t_{i+1}}$ is a $(1,l)$-quasi-geodesic, which proves our claim.
		
		\paragraph{} As we did previously with the $c_i$'s, we can prove that for every $i \in \intvald 0{m-1}$, the paths $\gamma$ and $\gamma'$ respectively restricted to $\intval{t_i}{t_{i+1}}$ and $\intval{s_i}{s_{i+1}}$ $D$-follows each other where
		\begin{displaymath}
			D = \max\left\{\dist x{x'}, \dist y{y'},5l'+l+14\delta\right\} + 3l + 3l' + 8\delta.
		\end{displaymath}
		Hence by Lemma~\ref{res: CH - concatenation of paths} $\gamma$ $D$-follows $\gamma'$.
	\end{proof}
	
	We can now finish the proof by showing that $\gamma$ is a $\sigma/3$-local $(1,l)$-quasi-geodesic.
	With the notations of the previous lemma, $\gamma$ $D$-follows $\gamma'$.
	In particular, it provides a map $\mu : \intval ab \rightarrow \intval {a'}{b'}$ (possibly different from $\theta$) which satisfies the axioms of Definition~\ref{def: CH - following}.
	Let $s,t \in \intval ab$, $s\leq t$ such that $\dist st \leq \sigma/3$.
	Since $\gamma'$ is a $\sigma/3$-local $(1,l')$-quasi-geodesic, we have
	\begin{displaymath}
		\dist{\mu(s)}{\mu(t)}\leq \dist st + 4D + 2l' \leq \frac \sigma2.
	\end{displaymath}
	Therefore $\gamma'$ restricted to $\intval{\mu(s)}{\mu(t)}$ is contained in a ball $B$ of radius $\sigma$.
	By hyperbolicity, this path is a $(1,l'+2\delta)$-quasi-geodesic in $B$.
	Let $\eta \in (0, 10^{-10}\sigma)$ and $\gamma_2$ be a $(1,\eta)$-quasi-geodesic joining $\gamma(s)$ to $\gamma(t)$.
	Note the we can chose $B$ in such a way that is also contains $\gamma_2$.
	Recall that $\dist {\gamma(s)}{\gamma'\circ\mu(s)} < D$. 
	The same holds for $t$.
	Hence by Lemma~\ref{res: CH - quasi-geodesics follow each other}, $\gamma_2$ follows $\gamma'$ restricted to $\intval{\mu(s)}{\mu(t)}$.
	Now we consider the concatenation of $\gamma$ restricted to $\intval as$, $\gamma_2$ and $\gamma$ restricted to $\intval tb$, we obtain a path which follows $\gamma'$ joining $x$ to $y$ (see Lemma~\ref{res: CH - concatenation of paths}).
	Consequently, it belongs to $\mathcal C$ and its length is bounded below by $L(\gamma)-l$.
	However $\gamma_2$ is a $(1,\eta)$-quasi-geodesic, thus
	\begin{displaymath}
		L(\gamma) \leq L \left(\restriction\gamma{\intval as}\right) + \dist{\gamma(s)}{\gamma(t)} + L \left(\restriction\gamma{\intval tb}\right) + l+ \eta. 
	\end{displaymath}
	It implies that the length of $\gamma$ restricted to $\intval st$ is less than $\dist{\gamma(s)}{\gamma(t)} + \eta + l$. 
	This inequality holds for every $\eta>0$ and $s,t \in \intval ab$ such that $\dist st \leq \sigma/3$.
	Consequently, $\gamma$ is a $\sigma/3$-locally $(1,l)$-quasi-geodesic path.
\end{proof}

\subsection{The space of quasi-geodesic paths.}

	The goal of this section is to prove the stability of locally quasi-geodesic paths in $X$.
	More precisely prove the following result
	\begin{prop}
	\label{res: CH - stability of quasi-geodesics}
		Let $l \in (0, 10^{-5}\sigma)$.
		Let $x$, $y$ and $y'$ be three points of $X$ such that $\dist y{y'} < 10^{-5}\sigma$.
		Let $\gamma$ (\resp $\gamma'$) be a $\sigma/3$-locally $(1,l)$-quasi-geodesic path joining $x$ to $y$ (\resp to $y'$).
		Then the Hausdorff distance between $\gamma$ and $\gamma'$ is less than $\max\left\{\dist y{y'}, 2l+5\delta\right\}+ 6l + 8\delta$.
	\end{prop}
	
	Let $l \in (0,10^{-5}\sigma)$.
	We fix a base point $x_0 \in X$ and denote by $\Gamma$ the set of $\sigma/3$-local $(1,l)$-quasi-geodesics starting at $x_0$.
	We endow $\Gamma$ with the equivalence relation defined previously: two elements $\gamma$ and $\gamma'$ of $\Gamma$ are equivalent if they have the same terminal points and they follow each other.
	The set of equivalence classes is denoted by $\tilde X$.
	Given a path $\gamma \in \Gamma$ we write $[\gamma]$ for its equivalence class.
	By convention $\tilde x_0$ stands for the point of $\tilde X$ represented by the constant path equal to $x_0$.
	The map $p : \Gamma \rightarrow X$ is defined to send every path to its terminal point.
	It induces a map from $\tilde X$ to $X$ that we still denote by $p$.
	We are going to prove that $p: \tilde X \rightarrow X$ is a path-connected cover of $X$.
	To that end we need first to define a topology on $\tilde X$.
	
	\paragraph{Topology of $\tilde X$.}
	Let $\tilde x= [\gamma]$ be a point of $\tilde X$.
	Let $\epsilon \in (0, 10^{-5}\sigma)$.
	We define $U_{\tilde x,\epsilon}$ as the set of equivalence classes $\tilde x'= [\gamma'] \in \tilde X$ such that $\gamma'$ follows $\gamma$ and conversely and $\dist{p(\tilde x)}{p(\tilde x')} < \epsilon$.
	According to Proposition~\ref{res: CH - transitivity}, $U_{\tilde x,\epsilon}$ is well-defined and does not depend on the choice of the path $\gamma$ representing $\tilde x$.
	
	\begin{lemm}
	\label{res: CH - base of open sets}
		The collection $\{U_{\tilde x,\epsilon}\}$ where $\tilde x \in \tilde X$ and $\epsilon \in (0,10^{-5}\sigma)$ forms a base of open sets of $\tilde X$.
	\end{lemm}
	
	\begin{proof}
		First note that for every $\tilde x \in \tilde X$ and $\epsilon \in (0, 10^{-5}\sigma)$, $\tilde x$ belongs to $U_{\tilde x, \epsilon}$.
		Therefore $\{ U_{\tilde x, \epsilon}\}$ covers $\tilde X$.
		Let $\tilde y = [\nu]$ be a point of $U_{\tilde x,\epsilon} \cap U_{\tilde x', \epsilon'}$ where $\tilde x = [\gamma]$, $\tilde x' = [\gamma']$ are two points of $\tilde X$ and $\epsilon, \epsilon' \in (0,10^{-5}\sigma)$.
		By definition there exists $\eta  \in (0, 10^{-5}\sigma)$ such that $\dist{p(\tilde y)}{p(\tilde x)} + \eta < \epsilon$ and $\dist{p(\tilde y)}{p(\tilde x')} + \eta < \epsilon$.
		We are going to prove that $U_{\tilde y,\eta}$ is contained in $U_{\tilde x,\epsilon} \cap U_{\tilde x', \epsilon'}$.
		Let $\tilde z = [c]$ be a point of $U_{\tilde y, \eta}$.
		By definition $c$ follows $\nu$ and conversely.
		Moreover, $\dist{p(\tilde z)}{p(\tilde y)}<\eta$.
		However $\tilde y$ lies in $U_{\tilde x,\epsilon}$
		Therefore $\nu$ follows $\gamma$ and conversely.
		Moreover, $\dist{p(\tilde y)}{p(\tilde x)}<\epsilon$.
		By Proposition~\ref{res: CH - transitivity}, $c$ follows $\gamma$ and conversely.
		Furthermore 
		\begin{displaymath}
			\dist{p(\tilde z)}{p(\tilde x)} \leq \dist{p(\tilde z)}{p(\tilde y)} + \dist{p(\tilde y)}{p(\tilde x)} < \dist{p(\tilde y)}{p(\tilde x)} +\eta <\epsilon.
		\end{displaymath}
		Consequently, $\tilde z$ belongs to $U_{\tilde x,\epsilon}$.
		We prove in the same way that $\tilde z \in U_{\tilde x',\epsilon'}$, which establishes our claim.
	\end{proof}
	
	\begin{lemm}
	\label{res: CH - path connected}
		The space $\tilde X$ is path-connected.
	\end{lemm}
	
	\begin{proof}
		Let $\tilde x$ be a point of $\tilde X$. We choose $\gamma : \intval ab \rightarrow X$ an element of $\Gamma$ representing $\tilde x$.
		We define a path of $X$ joining $\tilde x_0$ to $\tilde x$ as follows.
		\begin{displaymath}
			\begin{array}{rccc}
				F : & \intval ab & \rightarrow & \tilde X, \\
				& t & \rightarrow & \left[ \restriction\gamma{\intval at}\right].
			\end{array}
		\end{displaymath}
		We claim that $F$ is continuous.
		Let $t \in \intval ab$ and $\epsilon \in (0,10^{-5}\sigma)$.
		Let $s \in \intval ab$ such that $\dist st < \epsilon$.
		Since $\gamma$ is parametrized by arclength $\dist{\gamma(s)}{\gamma(t)}< \epsilon$.
		By Lemmas~\ref{res: CH - concatenation of paths} and \ref{res: CH - extended following path} the paths $\gamma$ restricted to $\intval as$ and $\intval at$ follow each other. 
		Moreover, $\dist{p\circ F(s)}{p\circ F(t)} = \dist{\gamma(s)}{\gamma(t)} < \epsilon$.
		Hence $F(s)$ belongs to $U_{F(t),\epsilon}$ which proves that $F$ is continuous.
		Every point of $\tilde X$ can be connected to $\tilde x_0$ by a continuous path, thus $\tilde X$ is path-connected.
	\end{proof}
	
	\begin{lemm}
	\label{res: CH - p continuous}
		The map $p$ is continuous.
	\end{lemm}

	\begin{proof}
		Let $x \in X$ and $\epsilon \in  (0,10^{-5}\sigma)$.
		By construction $p^{-1}(B(x,\epsilon))$ is exactly the union of all $U_{\tilde x,\epsilon}$ where $\tilde x \in p^{-1}(x)$.
		Consequently, $p$ is continuous.
	\end{proof}

	\begin{lemm}
	\label{res: CH - p onto on open sets}
		Let $\tilde x \in \tilde X$ and $\epsilon \in (0,10^{-5}\sigma)$.
		Then $p$ induces a homeomorphism from $U_{\tilde x, \epsilon}$ onto $B(p(\tilde x), \epsilon)$.
	\end{lemm}
	
	\begin{proof}
	 	By construction, the image by $p$ of $U_{\tilde x,\epsilon}$ is contained in $B(p(\tilde x),\epsilon)$.
		We first prove the converse inclusion.
		To that end we denote by $\gamma$ a path of $\Gamma$ representing $\tilde x$.
		Let $y \in B(p(\tilde x),\epsilon)$.
		The path $\gamma$ is a $\sigma/3$-locally $(1, l)$-quasi-geodesic. 
		Moreover, $\dist {p(\tilde x)}y < \epsilon$.
		By Proposition~\ref{res: CH - following path} there exists a $\sigma/3$-locally $(1,l)$-quasi-geodesic $\nu$ joining $x_0$ to $y$ which follows $\gamma$.
		Applying Proposition~\ref{res: CH - local quasi-geodesics follow each other}, $\gamma$ also follows $\nu$.
		Therefore $\nu$ defines a point $\tilde y \in U_{\tilde x,\epsilon}$ whose image by $p$ is $y$.
		Consequently, $p$ maps $U_{\tilde x,\epsilon}$ onto $B(p(\tilde x), \epsilon)$.
		
		\paragraph{}
		Let $\tilde y = [\nu]$ and $\tilde y' = [\nu']$ be two points $U_{\tilde x, \epsilon}$ such that $p(\tilde y) = p(\tilde y')$.
		By construction $\nu$ (\resp $\nu'$) follows $\gamma$ and conversely.
		Moreover, $\dist{p(\tilde y)}{p(\tilde x)} <\epsilon$ and $\dist{p(\tilde y)}{p(\tilde x)} <\epsilon$.
		By Proposition~\ref{res: CH - transitivity} $\nu$ and $\nu'$ follow each other.
		By assumption they also have the same extremities.
		It follows that $\nu$ and $\nu'$ are equivalent i.e., $\tilde y = \tilde y'$.
		Thus the restriction of $p$ to $U_{\tilde x,\epsilon}$ is one-to-one.
		Consequently, $p$ induces a bijection from $U_{\tilde x,\epsilon}$ onto $B(x,\epsilon)$.
		To complete the proof we have to show that this bijection is a homeomorphism.
		Since $p$ is continuous, we only need to show that it is open.
		To that end it is sufficient to observe that for every $U_{\tilde y,\eta}$ contained in $U_{\tilde x,\epsilon}$, $p(U_{\tilde y,\eta})$ is an open set which is a consequence of the first part of the proof.
	\end{proof}
	
	\begin{prop}
	\label{res: CH - cover}
		The map $p : \tilde X \rightarrow X$ is a covering map.
	\end{prop}
	
	\begin{proof}
		Every point of $x$ of $X$ can be joined to $x_0$ by a $(1,l)$-quasi-geodesic.
		Therefore $p$ is onto.
		We already proved that $p$ is a continuous map which induces for every $\tilde x \in \tilde X$ and every $\epsilon \in (0,10^{-5}\sigma)$ a homeomorphism from $U_{\tilde x, \epsilon}$ onto $B(p(\tilde x), \epsilon)$ (Lemma~\ref{res: CH - p continuous} and Lemma~\ref{res: CH - p onto on open sets}).
		Therefore it is sufficient to prove that for all $\tilde x, \tilde x' \in \tilde X$ such that $p(\tilde x) = p(\tilde x')$ for every $\epsilon \in (0,10^{-5}\sigma)$, if $U_{\tilde x,\epsilon} \cap U_{\tilde x',\epsilon}$ is non-empty then $\tilde x = \tilde x'$.
		We denote by $\gamma$ and $\gamma'$ paths of $\Gamma$ respectively representing $\tilde x$ and $\tilde x'$.
		Assume that there exists $\tilde y = [\nu]$ a point in $U_{\tilde x,\epsilon} \cap U_{\tilde x',\epsilon}$.
		By definition $\nu$ follows $\gamma$ (\resp $\gamma'$) and conversely.
		Moreover, $\dist{p(\tilde y)}{p(\tilde x)}< \epsilon$ and $\dist{p(\tilde y)}{p(\tilde x')}< \epsilon$.
		According to Proposition~\ref{res: CH - transitivity} $\gamma$ and $\gamma'$ follows each other.
		In addition they have the same extremities. 
		Consequently, they are equivalent, hence $\tilde x=\tilde x'$.
	\end{proof}
	
	\begin{proof}[Proof of Proposition~\ref{res: CH - stability of quasi-geodesics}]
		Being a covering map, $p$ induces a one-to-one homomorphism $p_* : \pi_1(\tilde X) \rightarrow \pi_1(X)$.
		We claim that this map is actually onto.
		By assumption $X$ is $10^{-5}\sigma$-simply connected.
		Hence its fundamental group is normally generated by the set of homotopy classes of loops of diameter less than $10^{-5}\sigma$.
		Let $c$ be such a loop.
		An element $g \in \pi_1(X)$ of the conjugacy class represented by $c$ is characterized by a path $c_0$ joining $x_0$ to a point of $x$ of $c$.
		Since $p$ is a covering map, $c_0$ can be lifted as a path $\tilde c_0$ starting at $\tilde x_0$.
		We denote by $\tilde x$ the terminal point of $\tilde c_0$.
		Hence $p(\tilde x) = x$.
		However $p$ induces an homeomorphism from $U_{\tilde x,10^{-5}\sigma}$ onto $B(x,10^{-5}\sigma)$ (Lemma~\ref{res: CH - p onto on open sets}).
		Since $c$ lies in $B(x,10^{-5}\sigma)$, it can be lifted as a loop $\tilde c$ in $\tilde X$.
		The paths $\tilde c_0$ and $\tilde c$ defines an element $\tilde g \in \pi_1(\tilde X)$ which is sent on $g$ by $p_*$.
		Hence $p_*$ is onto and therefore an isomorphism.
		It follows that $p$ is a bijection.
		Let $\gamma$ and $\gamma'$ be two $\sigma/3$-locally $(1,l)$-quasi-geodesic starting at $x_0$.
		We assume that their respective endpoints $y$ and $y'$ satisfies $\dist y{y'} < 10^{-5}\sigma$.
		Hence $\tilde y=[\gamma]$ and $\tilde y'=[\gamma']$ are two points of $\tilde X$ such that $\dist {p(\tilde y)}{p(\tilde y')}<10^{-5}\sigma$.
		Since $p$ is a bijection, Lemma~\ref{res: CH - p onto on open sets} implies that $\tilde y'$ belongs to $U_{\tilde y,10^{-5}\sigma}$.
		In particular, $\gamma$ and $\gamma'$ follow each other.
		However by Lemma~\ref{res: CH - local quasi-geodesics follow each other} they actually $D$-follow each other where $D = \max\left\{\dist y{y'}, 2l+5\delta\right\}+ 6l + 8\delta$.
		Therefore the Hausdorff distance between them is less than $D$.
		Note that no constant in the proof depends on the choice of the base point $x_0$.
		Therefore the results holds for any two $\sigma/3$-locally $(1,l)$-quasi-geodesics with the same initial point.
	\end{proof}
	
\subsection{Global hyperbolicity of $X$}
	
	\begin{prop}
	\label{res: CH - thin triangles}
		Let $l \in (0,10^{-10}\sigma)$.
		Let $x$, $y$ and $y'$ be three points of $X$.
		Let $\gamma$ (\resp $\gamma'$, $\nu$) be a $(1,l)$-quasi-geodesic joining $x$ to $y$ (\resp $x$ to $y'$, $y$ to $y'$).
		Then $\nu$ is contained in the $\Delta$-neighborhood of $\gamma\cup\gamma'$ where $\Delta = 50l+96\delta$.
	\end{prop}
	
	\begin{proof}
		We denote by $\gamma : \intval ab \rightarrow X$ and $\gamma' : \intval{a'}{b'} \rightarrow X$ parametrization by arclength of $\gamma$ and $\gamma'$.
		We define $s$ to be the largest number in $\intval ab$ such that $d(\gamma(s),\gamma') \leq 4l+5\delta$.
		Assume first that $s=b$.
		Therefore $y= \gamma(b)$ is $(4l+5\delta)$-close to $\gamma'$.
		Applying Proposition~\ref{res: CH - stability of quasi-geodesics}, $\nu$ lies in the $(10l+13\delta)$-neighborhood of $\gamma'$.
		
		\paragraph{}Assume now that that $s <b$.
		It follows that $d(\gamma(s), \gamma')=4l+5\delta$.
		We write $p= \gamma(s)$.
		Let $p'=\gamma'(s')$ be a projection of $p$ on $\gamma'$.
		We put $t = \min\{b, s+\sigma/3\}$, $t' = \min\{ b', s'+\sigma/3\}$, $q=\gamma(t)$ and $q'= \gamma'(t')$.
		See Figure~\ref{fig: thin triangle}.
		\begin{figure}[htbp]
		\centering
			\includegraphics{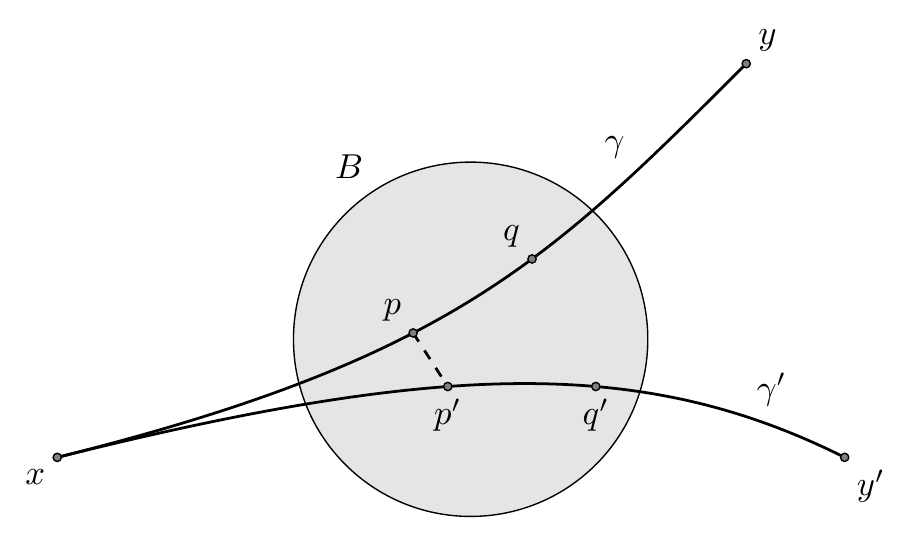}
		\caption{Thin ``geodesic'' triangle in $X$}
		\label{fig: thin triangle}
		\end{figure}
		The restriction of $\gamma'$ to $\intval {s'}{t'}$ and the point $p$ are simultaneously contained in a ball $B$ of radius $\sigma$.
		Hence by Proposition~\ref{res: quasi-convexity of quasi-geodesics} $\gro p{q'}{p'} \leq l + 2\delta$.
		On the other hand, by construction of $s$, for every $u \in \intval sb$,
		\begin{displaymath}
			\dist p{p'} = d(p,\gamma') \leq d(\gamma(u),\gamma') \leq \dist{\gamma(u)}{p'}.
		\end{displaymath}
		Hence $p$ is a projection of $p'$ on $\gamma$.
		Reasoning as previously, we get $\gro{p'}qp \leq l + 2\delta$.
		The four points $p$, $q$, $p'$ and $q'$ are contained in a ball of radius $\sigma$.
		Hence, by hyperbolicity,
		\begin{displaymath}
			\dist q{p'} + \dist{q'}p 
			\leq \max \left\{\fantomB \dist p{p'}+\dist q{q'},\dist pq+\dist{p'}{q'} \right\} + 2 \delta.
		\end{displaymath}
		Combined with the upper bound on the Gromov products, it leads to
		\begin{displaymath}
			\dist pq +\dist{p'}{q'} + 2\dist p{p'} \leq
			\max \left\{\fantomB \dist p{p'}+\dist q{q'},\dist pq+\dist{p'}{q'} \right\} + 4l +10 \delta.
		\end{displaymath}
		Recall that $\dist p{p'} = 4l+5\delta$. 
		Consequently, the maximum cannot be achieved by the second term. 
		Hence 
		\begin{equation}
		\label{eqn: CH - thin triangle equation}
			\dist q{q'} \geq \dist qp+ \dist p{p'} + \dist{p'}{q'} -4l-10\delta.
		\end{equation}
		Let $\eta \in (0,10^{-10}\sigma)$.
		We denote by $\nu_0$ be a $(1,\eta)$-quasi-geodesic joining $p$ to $p'$.
		We now define a path $\nu'$ by concatenating $\gamma$ restricted to $\intval sb$ (with the reverse parametrization), $\nu_0$ and $\gamma'$ restricted to $\intval {s'}{b'}$.
		Using (\ref{eqn: CH - thin triangle equation}), the length of the portion of $\nu'$ between $q$ and $q'$ is bounded above by
		\begin{displaymath}
			\dist qp+ \dist p{p'} + \dist{p'}{q'} + 2l + \eta \leq \dist q{q'} + 6l + \eta + 10 \delta.
		\end{displaymath}
 		Therefore it is a $(1,l')$-quasi-geodesic with $l' = 6l + \eta + 10 \delta$.
		On the other hand the portion of $\nu$ between $p$ and $y$ (\resp $p'$ and $y'$) is a $\sigma/3$-local $(1,l)$-quasi-geodesic.
		Note that if $t<b$ then $\dist st = \sigma/3$. 
		The same holds for $t'$.
		Hence $\nu'$ is a $\sigma/3$-local $(1,l')$-quasi-geodesic joining $y$ and $y'$.
		So is also $\nu$.
		Hence by Proposition~\ref{res: CH - stability of quasi-geodesics}, the Hausdorff distance between $\nu$ and $\nu'$ is less than $48l+8\eta+93\delta$.
		However, by construction $\nu'$ is contained in the $(2l+3\delta+\eta)$-neighborhood of $\gamma \cup \gamma'$.
		Hence $\nu$ lies in the $(50l+9\eta+96\delta)$-neighborhood of $\gamma\cup\gamma'$.
		This estimate holds for every $\eta>0$ which leads to the result.
	\end{proof}
	
	\begin{coro}
	\label{res: CH - dist vs gromov product}
		Let $x$, $y$ and $y'$ be three points of $X$.
		Let $l \in(0,10^{-10}\sigma)$.
		Let $\nu$ be a $(1,l)$-quasi-geodesic joining $y$ to $y'$.
		Then
		\begin{displaymath}
			\gro y{y'}x - \frac 12l \leq d(x,\nu) \leq \gro y{y'}x +101l+192\delta.
		\end{displaymath}
	\end{coro}
	
	\begin{proof}
		Let $\nu : \intval ab \rightarrow X$ be a parametrization by arclength of $\nu$.
		For every $t \in \intval ab$, the triangle inequality gives
		\begin{displaymath}
			\gro y{y'}x \leq \dist x{\nu(t)} + \gro y{y'}{\nu(t)} \leq \dist x{\nu(t)} + \frac 12 l.
		\end{displaymath}
		Hence $\gro y{y'}x \leq d(x,\nu) + l/2$, which gives the first inequality.
		Let $\gamma$ and $\gamma'$ be $(1,l)$-quasi-geodesics respectively joining $x$ to $y$ and $x$ to $y'$.
		By Proposition~\ref{res: CH - thin triangles}, $\nu$ lies in the $(50l+96\delta)$-neighborhood of $\gamma \cup \gamma'$.
		Therefore there exists $t\in \intval ab$ such  that $d(\nu(t),\gamma) \leq 50l+96\delta$ and $d(\nu(t),\gamma')\leq 50l+96\delta$.
		We denote by $u$ and $u'$ respective projections of $p$ on $\gamma$ and $\gamma'$.
		The triangle inequality leads to 
		\begin{displaymath}
			\gro y{y'}x \geq \dist {\nu(t)}x - \dist {\nu(t)}u - \dist{\nu(t)}{u'} - \gro xyu - \gro x{y'}{u'}.
		\end{displaymath}
		Thus $\gro y{y'}x \geq d(x,\nu) - 101l-192\delta$.
	\end{proof}
	
	\begin{coro}
	\label{res: CH - X globally hyperbolic}
		The space $X$ is $300\delta$-hyperbolic.
	\end{coro}

	\begin{proof}
		Let $x$, $y$, $y'$ and $t$ be four points of $X$.
		Let $l \in(0,10^{-10}\sigma)$.
		We denote by $\gamma$ (\resp $\gamma'$, $\nu$) a $(1,l)$-quasi-geodesic joining $x$ to $y$ (\resp $x$ to $y'$, $y$ to $y'$).
		Let $p$ be a projection of $t$ on $\nu$.
		By Proposition~\ref{res: CH - thin triangles}, $p$ lies in the $(50l+96\delta)$-neighborhood of $\gamma \cup \gamma'$.
		Thus the triangle inequality gives
		\begin{displaymath}
			d(t,\nu)  = \dist tp \geq \min\left\{\fantomB d(t,\gamma), d(t,\gamma')\right\} - 50l-96\delta.
		\end{displaymath}
		Combined with Corollary~\ref{res: CH - dist vs gromov product} it yields
		\begin{displaymath}
			\gro y{y'}t \geq \min\left\{\fantomB \gro yxt, \gro {y'}xt \right\} - 152l -288\delta.
		\end{displaymath}
		This inequality holds for every $l>0$, thus $\gro y{y'}t \geq \min\left\{ \gro yxt, \gro {y'}xt \right\} - 288\delta$.
		Hence $X$ is $288\delta$-hyperbolic.
	\end{proof}

\makebiblio

\end{document}